\documentclass[12pt,authoryear,preprint]{elsarticle}
\usepackage[utf8x]{inputenc}
\usepackage{pgf}
\usepackage{amsmath}   
\usepackage{amsmath}
\usepackage{amsmath,amssymb}

\usepackage{multirow}
\usepackage[]{algorithm2e}
\usepackage{rotating}
\usepackage{enumerate}
\usepackage{array}
\usepackage[labelformat=simple]{subfig}
\usepackage{float}
\usepackage{natbib}
\usepackage{color}
\usepackage{colortbl}
\usepackage{wasysym}
\usepackage{mathrsfs}
\usepackage{enumerate}
\usepackage{pxfonts}
 \usepackage{amsfonts}
\usepackage{dsfont}
\usepackage{array}
\usepackage{lscape}
\usepackage[pdftex,                %
    bookmarks         = true,
    bookmarksnumbered = true,
    pdfpagemode       = None,
    pdfstartview      = FitH,
    pdfpagelayout     = SinglePage,
    colorlinks        = true,
    urlcolor          = magenta,
    pdfborder         = {0 0 0}
    ]{hyperref}%

\usepackage{tabularx,ragged2e,booktabs,caption}
\newcolumntype{Y}{>{\centering\arraybackslash}X}
\newtheorem{The}{Theorem}[section]
\newtheorem{Pro}[The]{Proposition}
\newtheorem{Rq}[The]{Remark}
\newtheorem{Ex}[The]{Example}
\newtheorem{Def}[The]{Definition}
\newtheorem{Cor}[The]{Corollary}
\newtheorem{Lem}[The]{Lemma}

\DeclareMathOperator{\trace}{trace}
\DeclareMathOperator{\vech}{vech}
\DeclareMathOperator{\Bias}{Bias}
\DeclareMathOperator{\Var}{Var}
\DeclareMathOperator{\MSE}{MSE}
\DeclareMathOperator{\MISE}{MISE}
\DeclareMathOperator{\AMISE}{AMISE}
\DeclareMathOperator{\Cov}{Cov}

\usepackage{color}
\definecolor{darkred}{rgb}{.7, 0, 0}

\numberwithin{equation}{section}
\numberwithin{table}{section}
\numberwithin{figure}{section}
\usepackage{geometry}

\journal{  publication.}
\date{January 29, 2015}
\begin{document}

\begin{frontmatter}

\title{On multivariate associated kernels for smoothing \\general density functions}
\author[rvt]{C\'elestin C. Kokonendji}
\ead{celestin.kokonendji@univ-fcomte.fr}
\author[rvt]{Sobom M. Som\'e\corref{cor1}}
\ead{sobom.some@univ-fcomte.fr}
\cortext[cor1]{\textit{Corresponding author:} Universit\'e de Franche-Comt\'e, UFR Sciences et Techniques, Laboratoire de Math\'ematiques de Besan\c{c}on -- UMR 6623 CNRS-UFC, 16 route de Gray, 25030 Besan\c{c}on cedex, France; Tel. +33 381 666 398; Fax +33 381 666 623.}

\address[rvt]{University of Franche-Comt\'e, Besan\c{c}on, France}
\begin{abstract}
Multivariate associated kernel estimators, which depend on both target point and bandwidth matrix, are appropriate for partially or totally bounded distributions and generalize the classical ones as Gaussian. Previous studies on multivariate associated kernels have been restricted to product of univariate associated kernels, also considered having diagonal bandwidth matrices. However, it is shown in classical cases that for certain forms of target density such as multimodal, the use of full bandwidth matrices offers the potential for significantly improved density estimation. In this paper, general associated kernel estimators with correlation structure are introduced. Properties of these estimators are presented; in particular, the boundary bias is investigated. Then, the generalized bivariate beta kernels are handled with more details. The associated kernel with a correlation structure is built with a variant of the mode-dispersion method and two families of bandwidth matrices are discussed under the criterion of cross-validation. Several simulation studies are done. In the particular situation of bivariate beta kernels, it is therefore pointed out the very good performance of associated kernel estimators with correlation structure compared to the diagonal case. Finally, an illustration on real dataset of paired rates in a framework of political elections is presented.
\end{abstract}

\end{frontmatter}

\section{Introduction}

Nonparametric estimation of unknown densities on partially or totally bounded supports, with or without correlation in its multivariate components, is a recurrent practical problem. Because of symmetry, the multivariate classical or symmetric kernels, not depending on any parameter, are not appropriate for these densities. In fact, these estimators give weights outside the support causing a bias in boundary regions. In order to reduce the boundary problem with multivariate symmetric kernels as Gaussian, \cite{S02} and recently  \cite{ZAK14} have proposed adaptive full bandwidth matrix selection; but the bias does not disappear completely. \cite{{C99},{C00}} is one of the pioneers who has proposed, in univariate continuous case, some asymmetric kernels (i.e. beta and gamma) whose supports coincide with those of the densities to be estimated. Also recently,  \cite{L13} investigated several families of these univariate continuous kernels that he called univariate associated kernels; see also \cite{KSKZ07}, \cite{KSK11},  \cite{{ZAK12},{ZAK13a}} for univariate discrete situations. This procedure cancels of course the boundary bias; however, it creates a quantity in the bias of the estimator which needs reduction; see, for instance, \cite{MS13},  \cite{HiruSakudo13} and \cite{IgarashKakizawa2015}. 

Several approaches on multivariate kernel estimation have been proposed for various processings.  \cite{GP-C-GM13} used product of kernels for estimating the different nature of both directional and linear components of a random vector. A classical kernel density estimation on the rotation group appropriate for crystallographic texture analysis has been investigated by \cite{H13}. Symmetric kernel smoothers with univariate local bandwidth have been studied by  \cite{GM-LL-13} for semiparametric mixed effect models. \cite{GGS13} presented frontier estimation with classical kernel regression on high order moments. In discrete case, \citet{AJG76} provided kernel estimators for binary data while \cite{RL04} proposed the product of them with classical continuous one for smoothing regression on both categorical and continuous data. In the same spirit of \cite{RL04}, \cite{BR10} considered some products of different univariate associated kernels in continuous case; i.e. the bandwidth matrix obtained is diagonal. In the classical kernels case, \cite{CD11a} and \cite{CD11} have shown the importance of full bandwidth matrices for certain  target densities. See also \cite{HAZMR} for a support with arbitrary shape.

The main goal of this work is to introduce the multivariate associated kernels with the most general bandwidth matrix. In other words, the support of the suggested associated kernels coincides to the support of the densities to be estimated; also, the full bandwidth matrices take into account different correlation structures in the sample. Note that, a full bandwidth matrix significantly improves some complex target densities (e.g. multimodal); see \cite{S02}. In high dimensions, the computational choice of this full bandwidth matrix needs some special techniques. We can refer to \cite{{CD10}, {CD11a}} and \cite{CD11} for classical (symmetric) kernels. For illustrations in the present paper, we then focus on a bivariate case as beta kernel with correlation structure introduced by \cite{S66}; see also \cite{L96}. A motivation to investigate the smoothing of these densities on $\left[0, 1\right]\times\left[0, 1\right]$ comes from a joint distribution of two comparable proportions. Many datasets in $\left[0, 1\right]\times\left[0, 1\right]$ can be found in statistical problems, for example, for  comparing two rates. We shall examine the theoretical bias reduction and practical performances of the full bandwidth selection and two others bandwidth matrix parametrization using least squares or unbiased cross validation; see, e.g.,  \cite{WJ93}.

The rest of the paper is organized as follows. Section~\ref{sec:MultiAssKern} gives a complete definition of  multivariate associated kernels which includes both the product and the classical symmetric ones. A method to construct any multivariate associated kernel from a parametric probability density function (pdf) is then provided. Some pointwise properties of the corresponding estimator are investigated, in particular the convergence in the sense of mean integrated squared error (MISE) and an algorithm of the bias reduction. Section~\ref{sec:BivBetaKern} provides a particular study of a bivariate beta kernel with a correlation structure introduced by \cite{S66}. Also, some algorithms for the choice of the optimal bandwidth matrix by unbiased cross validation method are presented. This is followed, in Section~\ref{sec:Simulations}, by simulation studies and a real data analysis of electoral behaviour of a population with regard to a candidate. Especially, the role of forms of bandwidth matrices is explored in details. Section~\ref{sec:Conclusion} concludes with summary and final remarks, while the proof of a proposition is deferred to the appendix of Section~\ref{sec:Appendix}.

\section{Multivariate associated kernel estimators}
\label{sec:MultiAssKern}

Let $\mathbf{X}_{1}, \ldots, \mathbf{X}_{n}$ be independent and identically distributed (iid) random vectors with an unknown density function $f$ on $\mathbb{T}_{d}$, a subset of $\mathbb{R}^{d}$ ($d \geq 1$). As frequently observed in practice, the subset $\mathbb{T}_{d}$ might be unbounded, partially bounded or totally bounded as:
\begin{equation}\label{Examplesupport}
 \mathbb{T}_{d} =  \mathbb{R}^{d_{\infty}} \times \left[z, \infty\right)^{d_z} \times \left[u, w\right]^{d_{uw}}
 \end{equation} for given reals $u<w$ and $z$ with nonnegative values of $d_{\infty}$, $d_z$ and $d_{uw}$ in $\{0,1,\ldots,d\}$ such that $d = d_{\infty} + d_z + d_{uw}$. A multivariate  associated kernel estimator $ \widehat{f}_{n}$ of $f$ is simply defined by
\begin{equation}\label{Asskern}
 \widehat{f}_{n}(\mathbf{x}) = \dfrac{1}{n} \displaystyle \sum_{i=1}^{n} K_{\mathbf{x}, \mathbf{H}}(\mathbf{X}_{i}),~~~ \forall  \mathbf{x} \in  \mathbb{T}_{d} \subseteq \mathbb{R}^{d},
\end{equation}where $\mathbf{H}$ is a $d \times d$ bandwidth matrix (i.e. symmetric and positive definite) such that $\mathbf{H} \equiv \mathbf{H}_{n} \rightarrow \mathbf{0}_d$ (the $d \times d$ null matrix) as $n \rightarrow \infty$, and  $K_{\mathbf{x}, \mathbf{H}}(\cdot)$ is the so-called associated kernel, parametrized by $\mathbf{x}$ and $\mathbf{H}$, and precisely defined as follows.

\begin{Def}\label{Defkern}
Let $\mathbb{T}_{d}\left(\subseteq \mathbb{R}^{d}\right)$ be the support of the pdf to be estimated,  $\mathbf{x}  \in \mathbb{T}_{d}$ a target vector and   $\mathbf{H}$ a bandwidth matrix. A parametrized pdf $K_{\mathbf{x}, \mathbf{H}}(\cdot)$ of support $\mathbb{S}_{\mathbf{x}, \mathbf{H}}  \left(\subseteq \mathbb{R}^{d}\right)$
 is called ``multivariate (or general) associated kernel'' if the following conditions are satisfied:
\begin{align}
  &\mathbf{x} \in \mathbb{S}_{\mathbf{x}, \mathbf{H}},  \label{NoyAss1}\\
  &\mathbb{E}\left(\mathcal{Z}_{\mathbf{x}, \mathbf{H}}\right) = \mathbf{x} +  \mathbf{a}(\mathbf{x}, \mathbf{H}), \label{NoyAss2} \\
  & \rm{Cov}\left(\mathcal{Z}_{\mathbf{x}, \mathbf{H}}\right) =   \mathbf{B}(\mathbf{x}, \mathbf{H}), \label{NoyAss3}
\end{align}
where $\mathcal{Z}_{\mathbf{x}, \mathbf{H}}$ denotes the random vector with pdf  $K_{\mathbf{x},  \mathbf{H}}$ and both $\mathbf{a}(\mathbf{x},  \mathbf{H}) = \left(a_{1}(\mathbf{x},  \mathbf{H}), \ldots, a_{d}(\mathbf{x},  \mathbf{H})\right)^{\top}$ and  $\mathbf{B}(\mathbf{x},  \mathbf{H}) = \left(b_{ij}(\mathbf{x},  \mathbf{H})\right)_{i,j = 1, \ldots, d}$ tend, respectively, to the null vector $\mathbf{0}$ and the null matrix $\mathbf{0}_d$ as $\mathbf{H}$ goes to $\mathbf{0}_d$.
\end{Def}

\begin{Rq}\label{Remark1}
\begin{enumerate}[(i)]
\item The function $K_{\mathbf{x}, \mathbf{H}}(\cdot)$  is not necessary symmetric and is intrinsically linked to $\mathbf{x}$ and  $\mathbf{H}$.
\item The support $\mathbb{S}_{\mathbf{x},\mathbf{H}}$ is not necessary symmetric around of $\mathbf{x}$; it can  depend or not on $\mathbf{x}$ and $\mathbf{H}$.
\item The condition (\ref{NoyAss1}) can be viewed as $\cup_{\mathbf{x}\in\mathbb{T}_{d}}\mathbb{S}_{\mathbf{x},\mathbf{H}}\supseteq\mathbb{T}_{d}$ and it implies that the associated kernel takes into account the support $\mathbb{T}_{d}$ of the density $f$, to be estimated.
\item If $\cup_{\mathbf{x}\in\mathbb{T}_{d}}\mathbb{S}_{\mathbf{x},\mathbf{H}}$ does not contain $\mathbb{T}_{d}$ then this is the well-known problem of boundary bias.
\item Both conditions (\ref{NoyAss2}) and (\ref{NoyAss3}) indicate that the associated kernel is more and more concentrated around of $\mathbf{x}$ as $\mathbf{H}$ goes to $\mathbf{0}_d$. This highlights the peculiarity of associated kernel which can change its shape according to the target position.
\item The form of orientation of the kernel is controlled by the parametrization of bandwidth matrix $\mathbf{H}$; i.e. a full bandwidth matrix allows any orientation of the kernel and therefore any correlation structure.
\end{enumerate}
\end{Rq}

The following two examples provide well-known and also interesting particular cases of multivariate associated kernel estimators. The first can be seen as an interpretation of  associated kernels through symmetric kernels. The second deals on associated kernels without correlation structure.

\begin{Ex}\label{par:Exple1}
(Classical kernels) The kernel estimator $\widehat{f}_{n}$ of the density $f$, appropriate for unbounded supports in particular $\mathbb{R}^{d}$,  is usually defined by:
\begin{equation}
  \widehat{f}_{n}(\mathbf{x}) = \dfrac{1}{n } \displaystyle \sum_{i=1}^{n}K_{ \mathbf{H}}(\mathbf{x} - \mathbf{X}_{i}),~~ \forall  \mathbf{x} \in  \mathbb{T}_{d} = \mathbb{R}^{d}, \label{classkern}
\end{equation}where $\mathbf{x}$ is a target, $\mathbf{H}$ a bandwidth matrix and $K_{\mathbf{H}}(\mathbf{y}) = (1/\det\mathbf{H})K\left(\mathbf{H}^{-1}\mathbf{y}\right)$ or sometimes $K_{\mathbf{H}}(\mathbf{y}) = (1/\det\mathbf{H})^{1/2}K\left(\mathbf{H}^{-1/2}\mathbf{y}\right)$ for all $\mathbf{y} \in \mathbb{R}^{d}$. The function $K$ is the multivariate kernel assumed to be spherically symmetric and it does not depend on any parameter in particular $\mathbf{x}$ and $\mathbf{H}$. The kernel function has also mean vector and covariance matrix respectively equal to zero and $\mathbf{\Sigma}$; in general, the covariance matrix is the identity matrix: $\mathbf{\Sigma} = \mathbf{I}$. This function $K$ is here called classical kernel.
\end{Ex}

The following result connects a classical kernel to its corresponding symmetric or classical (multivariate) associated kernel.

\begin{Pro} \label{Pro1}
Let $\mathbb{R}^{d} = \mathbb{T}_{d}$ be the support of the density to be estimated. Let  $K$ be a classical kernel with support $\mathbb{S}_{d} \subseteq \mathbb{R}^{d}$, mean vector $\mathbf{0}$ and covariance matrix  $\mathbf{\Sigma}$. Given a target vector $\mathbf{x}\in\mathbb{R}^{d}$ and a bandwidth matrix $\mathbf{H}$, then the classical kernel induces the so-called classical (multivariate) associated kernel: (i)
\begin{eqnarray}\label{classcical}
K_{\mathbf{x}, \mathbf{H}}(\cdot) = \dfrac{1}{\det\mathbf{H}} K\left\{\mathbf{H}^{-1}(\mathbf{x} - \cdot)\right\}
\end{eqnarray} on $\mathbb{S}_{\mathbf{x}, \mathbf{H}} = \mathbf{x} - \mathbf{H}\mathbb{S}_{d}$ with $\mathbb{E}\left(\mathcal{Z}_{\mathbf{x}, \bf{H}}\right) = \mathbf{x}$ (i.e. $\mathbf{a}(\mathbf{x}, \mathbf{H}) = \mathbf{0}$) and $\rm{Cov}\left(\mathcal{Z}_{\mathbf{x}, \mathbf{H}}\right) =  \mathbf{H}\mathbf{\Sigma}\mathbf{H}$; (ii)
\begin{eqnarray*}
K_{\mathbf{x}, \mathbf{H}}(\cdot) = \dfrac{1}{(\det\mathbf{H})^{1/2}} K\left\{\mathbf{H}^{-1/2}(\mathbf{x} - \cdot)\right\}
\end{eqnarray*} on $\mathbb{S}_{\mathbf{x}, \mathbf{H}} = \mathbf{x} - \mathbf{H}^{1/2}\mathbb{S}_{d}$ with $\mathbb{E}\left(\mathcal{Z}_{\mathbf{x}, \bf{H}}\right) = \mathbf{x}$ (i.e. $\mathbf{a}(\mathbf{x}, \mathbf{H}) = \mathbf{0}$) and $\rm{Cov}\left(\mathcal{Z}_{\mathbf{x}, \mathbf{H}}\right) =  \mathbf{H}^{1/2}\mathbf{\Sigma}\mathbf{H}^{1/2}$.
\end{Pro}

\textbf{Proof.} We only proof (i) because (ii) is similar. From (\ref{Asskern}) and (\ref{classkern}) with $K_{\mathbf{H}}(\mathbf{y}) = (1/\det\mathbf{H})K\left(\mathbf{H}^{-1}\mathbf{y}\right)$, we easily deduce the expression (\ref{classcical}). From (\ref{classcical}) and Definition \ref{Defkern}, for a fixed $\mathbf{x} \in \mathbb{T}_{d}=\mathbb{R}^{d}$ and for all $\mathbf{t} \in \mathbb{T}_{d}=\mathbb{R}^{d}$, there exists $\mathbf{u} \in \mathbb{S}_{d}$ such that $\mathbf{u} = \mathbf{H}^{-1}(\mathbf{x} - \mathbf{t})$ and therefore $\mathbf{t} = \mathbf{x} - \mathbf{H}\mathbf{u}$. This implies, from (\ref{NoyAss1}), that $\mathbb{S}_{\mathbf{x}, \mathbf{H}} = \mathbf{x} - \mathbf{H}\mathbb{S}_{d}$. The last two results are simply derived from calculating the covariance matrix and the mean vector of $\mathcal{Z}_{\mathbf{x}, \bf{H}}$ (the random vector of pdf $K_{\mathbf{x}, \mathbf{H}}$) by making the previous substitution $\mathbf{u} = \mathbf{H}^{-1}(\mathbf{x} - \mathbf{t})$.$\blacksquare$\

It is known that the choice of classical kernels is not important; nevertheless, the best classical kernel is the \cite{E69} one in the sense of MISE with bounded support $\mathbb{S}_{d}$. The most popular is the Gaussian kernel with $\mathbb{S}_{d} = \mathbb{R}^{d}=\mathbb{T}_{d}$, $\mathbf{\Sigma} = \mathbf{I}$ and therefore $\mbox{Cov}\left(\mathcal{Z}_{\mathbf{x}, \mathbf{H}}\right) =  \mathbf{H}^{2}$; see  \cite{CD10} and \cite{ZAK14}. An interpretation of any classical multivariate associated kernel $K_{\mathbf{x},\mathbf{H}}(\cdot)$ can be presented as follows: through the symmetry property of the classical kernel, the mean $\mathbb{E}\left(\mathcal{Z}_{\mathbf{x}, \mathbf{H}}\right) = \mathbf{x}$ coincides with the mode which is the target $\mathbf{x}$; separately and in contrario to the general case (\ref{NoyAss3}), the dispersion measure around of the target $\mathbf{x}$, $\mbox{Cov}\left(\mathcal{Z}_{\mathbf{x},\mathbf{H}}\right)= \mathbf{H}\mathbf{\Sigma}\mathbf{H}$ which does not here depend on $\mathbf{x}$, serves essentially to the smoothing parameters or to the bandwidth matrix. This is the basic concept of general associated kernels and it is a different approach with respect to the convolution point of view. Note that the bandwidth matrix is similar to the dispersion matrix, which is symmetric and positive definite; see for instance \cite{J13}. For univariate dispersion parameter, we can refer to~\cite{J97} and~\cite{JK11} for different uses.

\begin{Ex}\label{par:Exple2}
(Multiple kernels) The product kernel estimator  introduced by~\cite{BR10} can be defined as a product of univariate associated kernel estimators of~\cite{L13}.  We here call it ``multiple associated kernel estimator'' $\widehat{f}_{n}$ of the density $f$:
\begin{equation}
  \widehat{f}_{n}(\mathbf{x})  = \dfrac{1}{n } \displaystyle \sum_{i=1}^{n} \prod_{j=1}^{d}K^{[j]}_{x_{j}, h_{jj}}(X_{ij}),~~ \forall  x_{j} \in  \mathbb{T}^{[j]}_{1} \subseteq \mathbb{R}, \label{prodkern}
\end{equation}where $\mathbb{T}^{[j]}_{1}$ is the support of  univariate margin of $f$ for $j = 1, \ldots, d$, $\mathbf{x} = (x_{1}, \ldots, x_{d})^{\top} \in \displaystyle \times_{j=1}^{d}\mathbb{T}^{[j]}_{1}$, $\mathbf{X}_{i} = (X_{i1}, \ldots, X_{id})^{\top}$ for $i = 1, \ldots, n$, and $h_{11}, \ldots, h_{dd}$ are the univariate  bandwidth parameters. The function $K_{x_{j}, h_{jj}}^{[j]}$ is the $j$th univariate   associated kernel on the support $\mathbb{S}_{x_{j}, h_{jj}} \subseteq \mathbb{R}$. In principle, this estimator is more appropriate for bounded or partially bounded distributions without correlation in its components. A particular multiple associated kernel estimator is obtained by using univariate classical kernels. 
\end{Ex}

In the following proposition, we point out that all multiple associated kernels are multivariate associated kernels without correlation structure in the bandwidth matrix.

\begin{Pro} \label{Proprod}
Let $\times_{j=1}^{d}\mathbb{T}^{[j]}_{1} = \mathbb{T}_{d}$ be the support of the density $f$ to be estimated with $\mathbb{T}^{[j]}_{1} (\subseteq \mathbb{R})$ the supports of univariate margins of $f$. Let $\mathbf{x} = (x_{1}, \ldots, x_{d})^{\top} \in \times_{j=1}^{d}\mathbb{T}^{[j]}_{1}$  and $\mathbf{H} = \mathbf{Diag} (h_{11}, \ldots, h_{dd})$ with $h_{jj} > 0$. Let  $K^{[j]}_{x_{j}, h_{jj}}$ be a univariate associated kernel (see Definition \ref{Defkern} for $d = 1$) with its corresponding random variable $\mathcal{Z}^{[j]}_{x_{j}, h_{jj}}$  on $\mathbb{S}_{x_{j}, h_{jj}} (\subseteq \mathbb{R})$ for all $j = 1, \ldots, d$. Then,  the multiple associated kernel is also a multivariate associated kernel:
\begin{eqnarray}\label{prodkern1}
  K_{\mathbf{x}, \mathbf{H}}(\cdot) = \prod_{j=1}^{d}K^{[j]}_{x_{j}, h_{jj}}(\cdot)
\end{eqnarray}on $\mathbb{S}_{\mathbf{x}, \mathbf{H}} =  \displaystyle \times_{j=1}^{d}\mathbb{S}_{x_{j}, h_{jj}}$ with $\mathbb{E}\left(\mathcal{Z}_{\mathbf{x}, \bf{H}}\right) = \left(x_{1} + a_{1}(x_{1}, h_{11}), \ldots, x_{d} + a_{d}(x_{d}, h_{dd})\right)^{\top}$ and $\Cov\left(\mathcal{Z}_{\mathbf{x}, \mathbf{H}}\right)$ = $ \mathbf{Diag}\left(b_{jj}(x_{j}, h_{jj})\right)_{j = 1, \ldots, d}$. In other words, the random variables  $\mathcal{Z}^{[j]}_{x_{j}, h_{jj}}$ are independent components of the random vector $\mathcal{Z}_{\mathbf{x}, \bf{H}}$.
\end{Pro}

\textbf{Proof.} From (\ref{Asskern}) and (\ref{prodkern}), the expression (\ref{prodkern1}) is easily deduced. The remainder results are obtained directly by calculating the mean vector  and covariance matrix of $(\mathcal{Z}^{[1]}_{x_{1}, h_{11}},  \ldots, \mathcal{Z}^{[d]}_{x_{d}, h_{dd}})^{\top} = \mathcal{Z}_{\mathbf{x}, \mathbf{H}} $  which is the random vector of the pdf (\ref{prodkern1}).$\blacksquare$

The multiple associated kernels have been illustrated in bivariate case by \cite{BR10}. For simulation studies, the authors used two independent univariate beta kernels and also two independent univariate gamma kernels. It is easy to generalize the investigation from two to more independent univariate associated kernels.
  
If we have an associated kernel, an estimator can be easily deduced as in (\ref{Asskern}). Otherwise, a construction of associated kernels is possible by using an appropriate pdf. The pdf used must have at least the same numbers of parameters than the number of components in the couple $\left(\mathbf{x}, \mathbf{H}\right)$ as parameters of the expected associated kernel. The rest of this section is devoted to a construction of the multivariate associated kernels and then to some properties of the corresponding estimators.

\subsection{Construction of general associated kernels}
\label{ssec:Construction}

In order to build a multivariate associated kernel $K_{\mathbf{x}, \mathbf{H}}(\cdot)$, we have to evaluate the dimensions of  $\mathbf{x}$ and $\mathbf{H}$. We always have $d$ components for the target vector $\mathbf{x}\in\mathbb{T}_d$ which is completely separate from the bandwidth matrix $\mathbf{H}$ in the classical multivariate associated kernel; but, in general, $\mathbf{x}$ is intrinsically linked to $\mathbf{H}$. 
\begin{table}[!htbp]
\begin{center}
\begin{tabular}{lccc}
\hline \hline
$\mathbf{H}$&   General    & Multiple  &Classical \\
\hline
Full &$d(d+3)/2$ &$2d$& $d(d+3)/2$ \\
Scott  &$d+1$ &$d+1$ & $d+1$   \\
Diagonal &$2d$  &$2d$ &$2d$ \\
 \hline \hline
\end{tabular}
 \caption{Numbers $k_{d}$ of parameters according to the form of bandwidth matrices for general, multiple and classical asociated kernels.} \label{parameterstable}
\end{center}
\end{table}

Table~\ref{parameterstable} gives the exact or minimal numbers $k_{d}$ ($> d$) of parameters in $\left(\mathbf{x}, \mathbf{H}\right)$ according to different forms of $\mathbf{H}$. The bandwidth matrix $\mathbf{H} = (h_{ij})_{i,j=1, \ldots, d}$ is said \textit{full} (i.e. with complete structure of correlations) and admits $d(d+1)/2$ parameters. It is said \textit{diagonal} if $\mathbf{H} = \mathbf{Diag}(h_{11}, \ldots, h_{dd})$, i.e. without correlation, and has only $d$  parameters. We denote by the \textit{Scott} bandwidth matrix the form $\mathbf{H} = h\mathbf{H}_{0}$ with only one parameter $h > 0$ and fixed $\mathbf{H}_{0} = (h_{ij,0})_{i,j = 1, \ldots, d}$; see \citet[page 154]{S92}. In practice, the matrix $\mathbf{H}_{0}$ can be fixed empirically from the data. Although $k_{d}$ is the same for classical and general associated kernels in Table~\ref{parameterstable}, the differences arise because of separation or not between $\mathbf{x}$ and $\mathbf{H}$ and, also, the presented $k_{d}$ are exact numbers for classical and minimal numbers for both general and multiple associated kernels.
It becomes clear that any pdf, with at least $k_{d}$ parameters and having a unique mode and a dispersion matrix, can lead to an associated kernel. Now, we introduce the notion of \textit{type of kernel} which is necessary for the construction from a given pdf.

\begin{Def} \label{Def1}
 A type of kernel $K_{\theta}$  is a parametrized pdf  with support $\mathbb{S}_{\theta} \subseteq \mathbb{R}^{d}$, $\theta \in  \Theta \subseteq \mathbb{R}^{k_{d}}$, such that $K_{\theta}$  is  squared integrable, unimodal with mode $\mathbf{m} \in \mathbb{R}^{d}$ and admitting a $d\times d$ dispersion matrix $\mathbf{D}$ (which is symmetric and positive definite); i.e. $\theta = \theta(\mathbf{m}, \mathbf{D})$ a $k_{d} \times 1$ vector with $k_{d}$ given in Table \ref{parameterstable} and the $d$ first coordinates of $\theta$ corresponds to those of $\mathbf{m}$.
\end{Def}

Let us denote by $\mathbf{D}_{\mathbf{x_0}} = \int_{\mathbb{R}^{d}}(\mathbf{x}-\mathbf{x_0})(\mathbf{x}-\mathbf{x_0})^{\top}K_{\theta}(d\mathbf{x})$ the dispersion matrix of $K_{\theta}$ around the fixed vector $\mathbf{x_0}$.
Here is a series of facts to have in mind.

\begin{Lem}\label{RmkTypKern}
Let $K_{\theta}$ be a type of kernel on $\mathbb{S}_{\theta} \subseteq \mathbb{R}^{d}$. The three following assertions are satisfied:
\begin{enumerate}[(i)]
\item the mode vector  $\mathbf{m}$ of $K_{\theta}$ always belongs in  $\mathbb{S}_{\theta}$;
\item $K_{\theta}(\mathbf{m}) \geq K_{\theta}(\boldsymbol{\mu} )$ where $\boldsymbol{\mu} $ is the mean vector of $K_{\theta}$; 

\item if $\mathbf{D}_{\mathbf{m}}$  tends to the null matrix, then  $\mathbf{D}_{\boldsymbol{\mu}} $ also goes  to the null matrix.
\end{enumerate}
\end{Lem}

\textbf{Proof.} (i) and (ii) are trivial. As for (iii), it is easy to check that
$\mathbf{D}_{\mathbf{m}} = \mathbf{D}_{\boldsymbol{\mu}}  + (\boldsymbol{\mu}-\mathbf{m})(\boldsymbol{\mu}-\mathbf{m})^{\top}$. Thus,  $\mathbf{D}_{\mathbf{m}}$ tends to the null matrix means $K_{\theta}$ goes to the Dirac probability in the sense of distribution; then, $\boldsymbol{\mu}-\mathbf{m}$ goes to the null vector and therefore $\mathbf{D}_{\boldsymbol{\mu}}$ also goes to the null matrix.$\blacksquare$

Without loss of generality, we only present a construction of general (or multivariate) associated kernels excluding both multiple and classical ones. In fact, \cite{L13} built some univariate associated kernels that can be used in multiple associated kernels. For this, he proposed a ``mode-dispersion principle'' saying that it must put the mode on the target and the dispersion parameter on the bandwidth. In the same spirit, we here propose a construction of general associated kernels using the {\it multivariate mode-dispersion} method given in (\ref{eqtheta}) below. 

Indeed, since $\theta = \theta(\mathbf{m}, \mathbf{D})$ is a $k_{d}\times 1$ vector, we must vectorize the couple $(\mathbf{x}, \mathbf{H})$ where $\mathbf{x}\in\mathbb{T}_d$ is the target and $\mathbf{H} =  (h_{ij})_{i,j = 1, \ldots, d}$ is the bandwidth matrix. According to the symmetry of $\mathbf{H}$, we use the so-called half vectorization of   $(\mathbf{x}, \mathbf{H})$. That  is defined as $\vech(\mathbf{x}, \mathbf{H}) = (\mathbf{x}, \vech\mathbf{H})$, where 
\begin{equation}
\vech\mathbf{H} = \left(h_{11}, \ldots, h_{1d},  h_{22}, \ldots, h_{2d}, \ldots, 
 h_{(d-1)(d-1)},  h_{(d-1)d}, h_{dd}\right)^{\top} 
\end{equation}is a $\left[d(d+1)/2\right]\times 1$ vector obtained by stacking the columns of the lower triangular of $\mathbf{H}$; see, e.g., \cite{HS79}.
 Making general associated kernels from a type of kernel $K_{\theta}$ on $\mathbb{S}_{\theta}$ with $\theta = \theta(\mathbf{m}, \mathbf{D})$ by multivariate mode-dispersion method requires solving the system of the equations
\begin{equation}
\left(\theta(\mathbf{m}, \mathbf{D})\right)^{\top} = (\mathbf{x}, \vech\mathbf{H})^{\top}. \label{eqtheta}
\end{equation}The solution of (\ref{eqtheta}), if there exists, is a $k_{d} \times 1$ vector denoted by $\theta(\mathbf{x}, \mathbf{H}) := \theta\left(\mathbf{m}(\mathbf{x}, \mathbf{H}), \mathbf{D}(\mathbf{x}, \mathbf{H})\right)$. For $d = 1$, the system (\ref{eqtheta}) provides the result in \cite{L13}. A light version will be presented for bivariate case ($d = 2$) in Section~\ref{ssec:BivbetaSarmanov}. For classical associated kernels, the system (\ref{eqtheta}) means to solve directly $(\mathbf{m}, \mathbf{D}) = (\mathbf{x}, \mathbf{H})$.  More generally, the solution of (\ref{eqtheta})
depends on the flexibility of the type of kernel $K_{\theta}$ and leads to the corresponding  associated kernel denoted $K_{\theta(\mathbf{x}, \mathbf{H})}$. The constructed associated kernel satisfies Definition~\ref{Defkern} of multivariate associated kernel:

\begin{Pro} \label{ProTyp-Ass}
The associated kernel function $K_{\theta(\mathbf{x}, \mathbf{H})}(\cdot)$, obtained from (\ref{eqtheta}) and having support $\mathbb{S}_{\theta(\mathbf{x}, \mathbf{H})}$, is such that:
\begin{align}
  &\mathbf{x} \in \mathbb{S}_{\theta(\mathbf{x}, \mathbf{H})}, \label{ProTyp-Ass1} \\
  &\mathbb{E}\left(\mathcal{Z}_{\theta(\mathbf{x}, \mathbf{H})}\right) - \mathbf{x}  =  \mathbf{a}_{\theta}(\mathbf{x}, \mathbf{H}),\label{ProTyp-Ass2}  \\
  & \rm{Cov}\left(\mathcal{Z}_{\theta(\mathbf{x}, \mathbf{H})}\right) =   \mathbf{B}_{\theta}(\mathbf{x}, \mathbf{H}), \label{ProTyp-Ass3}
\end{align}
where $\mathcal{Z}_{\theta(\mathbf{x}, \mathbf{H})}$ is the random vector with pdf $K_{\theta(\mathbf{x}, \mathbf{H})}$ and both $\mathbf{a}_{\theta}(\mathbf{x}, \mathbf{H}) = \left(a_{\theta1}(\mathbf{x},  \mathbf{H}), \ldots, a_{\theta d}(\mathbf{x},  \mathbf{H})\right)^{\top}$ and  $\mathbf{B}_{\theta}(\mathbf{x},  \mathbf{H}) = \left(b_{\theta ij}(\mathbf{x}, \mathbf{H})\right)_{i,j = 1, \ldots, d}$ tend, respectively, to the null vector $\mathbf{0}$ and the null matrix $\mathbf{0}_d$ as  $\mathbf{H}$ goes to $\mathbf{0}_d$.
\end{Pro}

\textbf{Proof.} The multivariate mode-dispersion method (\ref{eqtheta}) implies  $\theta(\mathbf{x}, \mathbf{H}) = \theta\left(\mathbf{m}(\mathbf{x}, \mathbf{H}), \mathbf{D}(\mathbf{x}, \mathbf{H})\right)$ which leads to $\mathbb{S}_{\theta(\mathbf{x}, \mathbf{H})} := \mathbb{S}_{\theta\left(\mathbf{m}(\mathbf{x}, \mathbf{H}), \mathbf{D}(\mathbf{x}, \mathbf{H})\right)}$. Since $K_{\theta}$ is unimodal of mode $\mathbf{m} \in \mathbb{S}_{\theta}$ (see Part~(i) of Lemma~\ref{RmkTypKern}), we obviously have (\ref{ProTyp-Ass1}), and then (\ref{NoyAss1}) is checked, because $\mathbf{m}$ is identified to $\mathbf{x}$ in (\ref{eqtheta}). 
Let $\mathcal{Z}_{\theta(\mathbf{m}, \mathbf{D})}$ be the random vector with unimodal pdf as the type of kernel $K_{\theta(\mathbf{m}, \mathbf{D})}$. From Part~(ii) of Lemma~\ref{RmkTypKern} we can write $$\mathbb{E}\left(\mathcal{Z}_{\theta(\mathbf{m}, \mathbf{D})}\right) = \mathbf{m} + \boldsymbol{\tau}(\mathbf{m}, \mathbf{D}),$$ where $\boldsymbol{\tau}(\mathbf{m}, \mathbf{D})$ is the difference between the mean vector $\mathbb{E}\left(\mathcal{Z}_{\theta(\mathbf{m}, \mathbf{D})}\right)$ and the mode vector $\mathbf{m}$ of $\mathcal{Z}_{\theta(\mathbf{m}, \mathbf{D})}$. Thus, from the mode-dispersion method (\ref{eqtheta}), we have  $\mathbf{m}  = \mathbf{x}$ and $\boldsymbol{\tau}(\mathbf{m}, \mathbf{D}) = \boldsymbol{\tau}(\mathbf{m}(\mathbf{x}, \mathbf{H}), \mathbf{D}(\mathbf{x}, \mathbf{H}))$; taking $\mathbf{a}_{\theta}(\mathbf{x}, \mathbf{H}) = \boldsymbol{\tau}(\mathbf{m}(\mathbf{x}, \mathbf{H}), \mathbf{D}(\mathbf{x}, \mathbf{H}))$, Part~(iii) of Lemma~\ref{RmkTypKern} leads to the second result (\ref{ProTyp-Ass2}), and therefore (\ref{NoyAss2}) is verified. Also, since $K_{\theta(\mathbf{m}, \mathbf{D})}$ admits a moment of second order, the covariance matrix of $K_{\theta(\mathbf{m}, \mathbf{D})}$ exists and that can be written as 
  $\rm{Cov}\left(\mathcal{Z}_{\theta(\mathbf{m}, \mathbf{D})}\right) = \mathbf{B}_{\theta}\left(\mathbf{m}, \mathbf{D}\right)$; solving (\ref{eqtheta}) and then taking $\mathbf{B}_{\theta}(\mathbf{x},  \mathbf{H}) := \mathbf{B}_{\theta}\left(\mathbf{m}(\mathbf{x}, \mathbf{H}), \mathbf{D}(\mathbf{x}, \mathbf{H})\right)$, the last result (\ref{ProTyp-Ass3}) holds in the sense of (\ref{NoyAss3}) using again Part~(iii) of Lemma~\ref{RmkTypKern}.$\blacksquare$
  
In practice, both characteristics $\mathbf{a}_{\theta}(\mathbf{x}, \mathbf{H})$ and $\mathbf{B}_{\theta}(\mathbf{x}, \mathbf{H})$ are derived from the calculation of the mean vector and covariance matrix of $\mathcal{Z}_{\theta(\mathbf{x}, \mathbf{H})}$ in terms of the mode vector $\mathbf{m}(\mathbf{x}, \mathbf{H})$ and the dispersion matrix $\mathbf{D}(\mathbf{x}, \mathbf{H})$. 
In a general way, a given $K_{\mathbf{x}, \mathbf{H}}$ or the constructed associated kernel $K_{\theta(\mathbf{x}, \mathbf{H})}$  in Proposition~\ref{ProTyp-Ass} (that we will call standard version) creates a quantity in the bias of the kernel density estimation. In order to eliminate this quantity in the larger part of the support $\mathbb{T}_{d}$ of the density to be estimated, we will also study  a modified version of the standard one. The following two subsections investigate some properties of these estimators.

\subsection{Standard version of the estimator}
\label{ssec:Standard estimator}

Here, we give some properties of the estimator $\widehat{f}_{n}$ of $f$ through a given associated kernel presented in Definition~\ref{Defkern} or the constructed associated kernel in Proposition~\ref{ProTyp-Ass}; i.e. $K_{\theta(\mathbf{x}, \mathbf{H})}(\cdot)\equiv K_{\mathbf{x}, \mathbf{H}}(\cdot)$. For a given bandwidth matrix $\mathbf{H}$, similarly to (\ref{Asskern}) we consider
\begin{equation}
 \widehat{f}_{n}(\mathbf{x})= \dfrac{1}{n}\displaystyle\sum_{i=1}^{n}K_{\theta(\mathbf{x}, \mathbf{H})}(\mathbf{X}_{i}),\quad\forall \mathbf{x}\in\mathbb{T}_{d}. \label{AssCKernest}
\end{equation}

\begin{Pro} \label{Pro2}
For given $\mathbf{x} \in \mathbb{T}_{d}$,
\begin{equation}
 \mathbb{E}\left\{\widehat{f}_{n}(\mathbf{x})\right\} = \mathbb{E}\left\{ f\left(\mathcal{Z}_{\theta(\mathbf{x}, \mathbf{H})}\right)\right\} \label{EstversKern}
  \end{equation}and $\widehat{f}_{n}(\mathbf{x}) \geq 0$. Furthermore,  one has
\begin{align}
  &\int_{\mathbb{T}_{d}} \widehat{f}_{n}(\mathbf{x}) d\bold{x} = \Lambda, \label{est2}
\end{align}where the total mass $\Lambda = \Lambda(n; \mathbf{H}, K_{\theta})$ is a positive real and, it is also a finite constant if $\int_{ \mathbb{T}_{d}} K_{\theta(\mathbf{x}, \mathbf{H})} (\mathbf{t})d\bold{x} <\infty$ for all $\mathbf{t} \in \mathbb{T}_{d}$.
\end{Pro}

\textbf{Proof.} The first result (\ref{EstversKern}) is straightforwardly obtained from (\ref{AssCKernest}) as follows:
\begin{equation*}
 \mathbb{E}\left\{\widehat{f}_{n}(\mathbf{x})\right\} = \int_{\mathbb{S}_{\theta(\mathbf{x}, \mathbf{H})} \cap \mathbb{T}_{d}} K_{\theta(\mathbf{x}, \mathbf{H})} (\mathbf{t})f(\mathbf{t})d\bold{t} =  \mathbb{E}\left\{ f\left(\mathcal{Z}_{\theta(\mathbf{x}, \mathbf{H})}\right)\right\}.
\end{equation*}Also, the estimates $\widehat{f}_{n}(\mathbf{x}) \geq 0$ and the total mass $\Lambda > 0$ stem immediately from the fact that $K_{\theta(\mathbf{x}, \mathbf{H})}(\cdot)$ is a pdf.
Finally, the values of $\mathbf{X}_{i}$ belonging to the set $\mathbb{T}_{d}$, we have (\ref{est2}) as
\begin{equation*}
\int_{\mathbb{T}_{d}}\widehat{f}_{n}(\mathbf{x})d\bold{x} = \frac{1}{n} \displaystyle \sum_{i = 1}^{n} \int_{\mathbb{T}_{d}}K_{\theta(\mathbf{x}, \mathbf{H})} (\mathbf{X}_{i})d\bold{x} < \infty,
\end{equation*}since the integration vector of $\int_{\mathbb{T}_{d}}K_{\theta(\mathbf{x}, \mathbf{H})} (\mathbf{X}_{i})d\bold{x}$ is on the target $\mathbf{x}$ which is a parameter of $K_{\theta(\mathbf{x}, \mathbf{H})}(\cdot)$.$\blacksquare$

From the above Proposition~\ref{Pro2}, the total mass $\Lambda$ of $\widehat{f}_{n}$ by a non-classical associated kernel generally fails to be equal to 1; an illustration for $d=2$ is given in Table~\ref{Tablnormconst} below. Hence, non-classical associated kernel estimators $\widehat{f}_{n}$ are improper density estimates or as kind of ``balloon estimators''; see \cite{S02} and \cite{ZAK14}. The fact that $\Lambda = \Lambda(n; \mathbf{H}, K_{\theta})$ is close to 1 can find a statistical explanation in both Examples 1 and 2 of \cite{RT96}, or by showing  $\mathbb{E}\left(\int_{\mathbb{T}_{d}}\widehat{f}_{n}(\mathbf{x})d\bold{x}\right) = 1$ which does not depend on $n$: $\mathbb{E}\left(\int_{\mathbb{T}_{d}}\widehat{f}_{n}(\mathbf{x})d\bold{x}\right)  = \mathbb{E}\left\{ \phi(\mathbf{X}_{1}\mathds{1}_{\mathbb{S}_{\theta(\mathbf{x}, \mathbf{H})} \cap \mathbb{T}_{d}}; \mathbf{H}, K_{\theta})\right\}$ with $\phi(\mathbf{t}; \mathbf{H}, K_{\theta}) = \int_{ \mathbb{T}_{d}} K_{\theta(\mathbf{x}, \mathbf{H})} (\mathbf{t})d\bold{x} <\infty$ for all $\mathbf{t} \in \mathbb{T}_{d}$. Without loss of generality, we study $\mathbf{x} \mapsto\widehat{f}_{n}(\mathbf{x})$ up to normalizing constant which is used at the end of the density estimation process.

\begin{Pro} \label{Pro3}
Let $\mathbf{x} \in \mathbb{T}_{d}$ be a target and a bandwidth matrix $\mathbf{H} \equiv \mathbf{H}_{n}\to \mathbf{0}_d$ as $n\to \infty$. Assume $f$ in the class $\mathscr{C}^{2}(\mathbb{T}_{d})$, then
\begin{eqnarray}
 \Bias\left\{\widehat{f}_{n}(\mathbf{x})\right\}
&=& \mathbf{a}^{\top}_{\theta}(\mathbf{x}, \mathbf{H})\nabla f(\mathbf{x}) + \dfrac{1}{2}\trace\left[\left\{\mathbf{a}_{\theta}(\mathbf{x}, \mathbf{H})\mathbf{a}^{\top}_{\theta}(\mathbf{x}, \mathbf{H}) + \mathbf{B}_{\theta}(\mathbf{x}, \mathbf{H}) \right\}\nabla^{2}f\left(\mathbf{x}\right)\right]   \nonumber \\
&&+ \oldstylenums{0}\left\{\trace\left(\mathbf{H}^{2}\right)\right\}.  \label{biaisg}
\end{eqnarray}
Furthermore, if  $f$ is bounded on  $\mathbb{T}_{d}$ then there exists the largest positive real number $r_{2} = r_{2}\left(K_{\theta}\right)$  such that
 $\left\|K_{\theta(\mathbf{x}, \mathbf{H})}\right\|^{2}_{2}  \lesssim  c_{2}(\mathbf{x})(\det\mathbf{H})^{-r_{2}}$,  $0 \leq c_{2}(\mathbf{x}) \leq \infty$ and
\begin{equation}
  \Var\left\{\widehat{f}_{n}(\mathbf{x})\right\}  = \frac{1}{n} \left\|K_{\theta(\mathbf{x}, \mathbf{H})}\right\|^{2}_{2}f(\mathbf{x}) +  \oldstylenums{0}\left(n^{-1}\;(\det\mathbf{H})^{-r_{2}}\right), \label{var}
\end{equation}
with $\left\|K_{\theta(\mathbf{x}, \mathbf{H})}\right\|^{2}_{2} := \int_{\mathbb{S}_{\theta(\mathbf{x}, \mathbf{H})}} K_{\theta(\mathbf{x}, \mathbf{H})}^{2}(\mathbf{u})d\bold{u}$ and where ``$\lesssim$'' stands for ``$\leq$ and then approximation as $n \to \infty$''.
\end{Pro}
\textbf{Proof.} See Appendix in Section~\ref{sec:Appendix}.

The univariate case ($d = 1$) of Proposition~\ref{Pro3} can be found in \cite{L13}; see \cite{{C99},{C00}} for both beta and gamma kernels with $r_{2} = 1/2$ and, also, \cite{HiruSakudo13} for other values of $r_{2}$. In the situation of multiple associated kernels (\ref{prodkern1}), in contrario to (\ref{biaisg}) the general representation (\ref{var}) is simply expressed in terms of univariate associated kernels as follows:

\begin{Cor}\label{Corollary}
Under (\ref{prodkern1}) with $r_{2,j} =  r_{2}\left(K_{\theta_{j}}^{[j]}\right)$ and where $K_{\theta_{j}}^{[j]} \equiv K^{[j]}$ refers to the $j$th univariate type of kernel, then (\ref{var}) is equivalent to 
$$ \Var\left\{\widehat{f}_{n}(\mathbf{x})\right\}  = \frac{1}{n} \left(\prod_{j=1}^{d}\left\|K^{[j]}_{\theta_{j}(x_{j}, h_{jj})}\right\|^{2}_{2}\right)f(\mathbf{x}) +  \oldstylenums{0}\left(n^{-1}\; \prod_{j=1}^{d}h_{jj}^{-r_{2,j}}\right).$$
\end{Cor}

\textbf{Proof.}  Easy.$\blacksquare$

Now, we recall  natural measures for assessing the similarity of the multivariate associated kernel estimator $\widehat{f}_{n}$ according to the true pdf $f$, to be estimated. Since the pointwise measure is the mean squared error (MSE) and expressed by
\begin{eqnarray}
 \MSE\left(\mathbf{x}\right) = \mbox{Bias}^{2}\left\{\widehat{f}_{n}(\mathbf{x})\right\} +    \mbox{Var}\left\{\widehat{f}_{n}(\mathbf{x})\right\}, \label{MSE}
\end{eqnarray}
the integrated form of MSE on $\mathbb{T}_{d}$ is $\MISE\left(\widehat{f}_{n}\right) =  \int_{\mathbb{T}_{d}} \MSE\left(\mathbf{x}\right)d\bold{x}$ and its approximate expression satisfies 
\begin{eqnarray}
\label{amiseg}
\AMISE\left(\widehat{f}_{n}\right)  &=&  \int_{\mathbb{T}_{d}}\left(\left[\mathbf{a}^{\top}_{\theta}(\mathbf{x}, \mathbf{H})\nabla f\left(\mathbf{x}\right) + \dfrac{1}{2}\trace\left\{\left(\mathbf{a}_{\theta}(\mathbf{x}, \mathbf{H})\mathbf{a}^{\top}_{\theta}(\mathbf{x}, \mathbf{H}) + \mathbf{B}_{\theta}(\mathbf{x}, \mathbf{H}) \right)\nabla^{2}f\left(\mathbf{x}\right)\right\}\right]^{2} \right. \nonumber \\ 
&& + \left.\frac{1}{n} \left\|K_{\theta(\mathbf{x}, \mathbf{H})}\right\|^{2}_{2}f(\mathbf{x})\right)d\bold{x}.
\end{eqnarray} 
In the general case of associated kernels with correlation structure, an optimal bandwidth matrix that minimizes the AMISE (\ref{amiseg}) still remains challenging problem, even in the bivariate case. However, the particular case of diagonal bandwidth matrix is $\mathbf{H}_{{\rm opt, diag}} = \oldstylenums{O}\left(n^{-2/(d + 4)}\right)$ for multiple gamma kernels; see \cite{BR10} for further details. So, the next result only gives the optimal bandwidth matrix for the Scott form which still has a correlation structure.

In fact, let us consider $\mathbf{H} = h\mathbf{H}_{0}$ be a Scott bandwidth matrix with fixed matrix $\mathbf{H}_{0}$ and positive $h \equiv h_{n}\to 0$ as $n\to\infty$. From Definition~\ref{Defkern}, Proposition~\ref{Pro1} and Proposition~\ref{ProTyp-Ass}, it follows that there exists a $d \times 1$ vector  $\mathbf{c}^{*}(\mathbf{x}, \mathbf{H}_{0}) = (c^{*}_{1}(\mathbf{x}, \mathbf{H}_{0}), \ldots, c^{*}_{d}(\mathbf{x}, \mathbf{H}_{0}))^{\top}$ and a $d \times d$ matrix $\mathbf{C}^{**}(\mathbf{x}, \mathbf{H}_{0}) = (c^{**}_{ij}(\mathbf{x}, \mathbf{H}_{0}))_{i, j = 1, \ldots, d}$ of finite constants connected respectively to $\mathbf{a}_{\theta}(\mathbf{x}, \mathbf{H}) = (a_{\theta 1}(\mathbf{x}, \mathbf{H}), \ldots, a_{\theta d}(\mathbf{x}, \mathbf{H}))^{\top}$  and $\mathbf{B}_{\theta}(\mathbf{x}, \mathbf{H}) = (b_{\theta ij}(\mathbf{x}, \mathbf{H}))_{i, j = 1, \ldots, d}$ such that  for all $i, j = 1, \ldots, d$, $a_{\theta i}(\mathbf{x}, \mathbf{H}) \leq hc^{*}_{i}(\mathbf{x}, \mathbf{H}_{0})$ and $b_{\theta ij}(\mathbf{x}, \mathbf{H}) \leq h^{2}c^{**}_{ij}(\mathbf{x}, \mathbf{H}_{0})$. Using Proposition~\ref{Pro3} and assuming  $f$ such that  its all first and second partial derivatives are bounded, one has:
\begin{equation*}
 \Bias\left\{\widehat{f}_{n}(\mathbf{x})\right\}\leq hs_{1}\left(\mathbf{x}\right) \textrm{ and } \mbox{Var}\left\{\widehat{f}_{n}(\mathbf{x})\right\}\leq n^{-1}h^{-dr_{2}}s_{2}\left(\mathbf{x}\right),
\end{equation*}
where $s_{1}(\mathbf{x}) \geq  \mathbf{c}^{*\top}(\mathbf{x}, \mathbf{H}_{0})\nabla f\left(\mathbf{x}\right) + (h/2)\trace\left[\left\{\mathbf{c}^{*\top}(\mathbf{x}, \mathbf{H}_{0})\mathbf{c}^{*}(\mathbf{x}, \mathbf{H}_{0}) + \mathbf{C}^{**}(\mathbf{x}, \mathbf{H}_{0}) \right\}\nabla^{2}f\left(\mathbf{x}\right)\right]$ and $s_{2}(\mathbf{x})$ are positive scalars. From (\ref{MSE}), it follows that $ \MSE\left(\mathbf{x}\right)\leq h^{2}s_{1}^{2}\left(\mathbf{x}\right) + n^{-1}h^{-dr_{2}}s_{2}\left(\mathbf{x}\right)$.
By integration of $\MSE\left(\mathbf{x}\right)$, one obtains
\begin{equation}\label{AMISEH0}
 \AMISE(\widehat{f}_{n,h\mathbf{H}_{0},K_{\theta},f})\leq h^{2}t_{1}+n^{-1}h^{-dr_{2}}t_{2},
\end{equation} with $t_{1}$ and $t_{2}$ the anti-derivatives of respectively  $s_{1}^{2}\left(\mathbf{x}\right)$ and $s_{2}\left(\mathbf{x}\right)$ on $\mathbb{T}_{d}$. Taking the derivative of the second member of the inequality (\ref{AMISEH0}) equal to $0$ leads to the following proposition.

\begin{Pro} \label{Pro4}
Under assumption of Proposition~\ref{Pro3}, let $\mathbf{H} = h\mathbf{H}_{0}$ be a Scott bandwidth matrix with fixed matrix $\mathbf{H}_{0}$ and positive $h \equiv h_{n}\to 0$ as $n\to\infty$ and such that the right member of (\ref{AMISEH0}) satisfies
\begin{equation}\label{RightAMISEH0}
0<h^{2}t_{1}+n^{-1}h^{-dr_{2}}t_{2}< \infty,
\end{equation} where $r_2=r_{2}\left(K_{\theta}\right)$ given in (\ref{var}).
 Then, the optimal  bandwidth matrix $ \mathbf{H}_{{\rm opt, Scott}}$ minimizing the $\AMISE$ in (\ref{amiseg}) is
\begin{equation}\label{OptimalH0}
 \mathbf{H}_{{\rm opt, Scott}} = C n^{-1/(dr_{2} + 2)}\mathbf{H}_{0},
\end{equation}
where C is a positive constant.
\end{Pro}
Note that we can use the Scott bandwidth matrix $\mathbf{H} = h\mathbf{H}_{0}$ if (\ref{RightAMISEH0}) holds. However, in practice, we cannot check (\ref{RightAMISEH0}) because $f$ is unknown. But, if the quantity  $h^{2}t_{1}+n^{-1}h^{-dr_{2}}t_{2}$ of (\ref{RightAMISEH0})  becomes 0 (resp. $\infty$) then one observes an undersmoothing (resp. oversmoothing). Finally, the practical choice of $\mathbf{H}_{0}$ in (\ref{OptimalH0}) can be the sample covariance matrix.

Unlike to classical associated kernels for $\mathbb{T}_{d}=\mathbb{R}^{d}$ (Example \ref{par:Exple1}), the choice of non-classical associated kernels is very important for the support $\mathbb{T}_{d}\left(\subseteq \mathbb{R}^{d}\right)$ of the pdf $f$ to be estimated; see Parts (i)-(iv) of Remark~\ref{Remark1}. Also, from Parts (v)-(vi) of Remark~\ref{Remark1}, different positions of the target $\mathbf{x} \in \mathbb{T}_{d}$ and correlation structure need a suitable general associated kernel. Nevertheless, the selection of bandwidth matrix remains crucial when the general associated kernel is chosen; see, e.g. \cite{CD11a} and \cite{CD11}. Here, we consider the multivariate least squares cross validation (LSCV) method to select the bandwidth matrix. From (\ref{AssCKernest}), the LSCV method is based on the minimization of the integrated squared error (ISE) which can be written as
$$\mbox{ISE}(\mathbf{H})=\displaystyle\int_{\mathbb{T}_{d}}\widehat{f_{n}}^{2}(\mathbf{x})d\bold{x}-2\displaystyle\int_{\mathbb{T}_{d}}\widehat{f}_{n}(\mathbf{x})f(\mathbf{x})d\bold{x}+
\displaystyle\int_{\mathbb{T}_{d}}f^{2}(\mathbf{x})d\bold{x}.$$
Minimizing this $\mbox{ISE}(\mathbf{H})$ means to minimize the two first terms. However, we need to estimate the second term since it depends on the unknown pdf $f$. The LSCV estimator of $\mbox{ISE}(\mathbf{H})-\displaystyle\int_{\mathbb{T}_{d}}f^{2}(\mathbf{x})d\bold{x}$ is
\begin{equation*}
  \mbox{LSCV}(\mathbf{H})  =  \int_{\mathbb{T}_{d}}\left\{\widehat{f}_{n}(\mathbf{x})\right\}^{2}d\bold{x} - \dfrac{2}{n} \displaystyle\sum_{i = 1}^{n}\widehat{f}_{n, -i}(\mathbf{X}_{i}),
\end{equation*}
where $\widehat{f}_{n, -i}\left(\mathbf{X}_{i}\right) = (n - 1)^{-1}\displaystyle\sum_{j \neq i} K_{\theta(\mathbf{X}_{i}, \mathbf{H})}(\mathbf{X}_{j})$ is being computed as $\widehat{f}_{n}(\mathbf{X}_{i})$ excluding the observation $\mathbf{X}_{i}$. The bandwidth matrix obtained by the LSCV rule selection is defined as follows:
\begin{equation} \label{Hcv}
  \mathbf{\widehat{H}}  =  \displaystyle \mbox{arg }  \underset{\mathbf{H} \in \mathcal{H} }{ \min} \mbox{ LSCV}(\mathbf{H}),
\end{equation}where $\mathcal{H}$ is the set of all positive definite full bandwidth matrices. The LSCV rule is the same for the Scott and diagonal bandwidth matrices where $\mathcal{S}$ and $\mathcal{D}$ are their respective sets. The difficulty comes from the level of dimension of $\mathbf{H}$ which is, respectively, $d(d+1)/2$, $1$ and $d$ for $\mathcal{H}$, $\mathcal{S}$ and $\mathcal{D}$. Thus, for high dimension $d>2$, the set of the Scott bandwidth matrices might be a good compromise between computational problems and correlation structures in the sample.

\subsection{Modified version of the estimator}
\label{ssec:Modifiedversion}

Following \cite{{C99},{C00}}, \cite{L13}, \cite{MS13},  \cite{HiruSakudo13} and \cite{IgarashKakizawa2015} in univariate case, a second version of the estimator (\ref{AssCKernest}) is sometimes necessary. Indeed, the presence of the non null term $\mathbf{a}_{\theta}(\mathbf{x},\mathbf{H})$ with the gradient $\nabla f(\mathbf{x})$ in (\ref{biaisg}) increases the pointwise bias of $\widehat{f}_{n}\left(\mathbf{x}\right)$. Thus, we propose below an algorithm for eliminating the term of gradient in the largest region of $\mathbb{T}_{d}$. Since \cite{BR10} has shown the results for multiple associated kernels, we here investigate the case of general associated kernels with $d>1$.

The algorithm of bias reduction has two steps. The first step consists to define both inside and boundary regions. The second one deals on the modified associated kernel which leads to the bias reduction in the interior domain.

\textbf{First step}. Partitioning $\mathbb{T}_{d}$ into two regions of order $\bf{\alpha}({\mathbf{H}})$ which is  a $d \times 1$ vector,  and where  $\mathbf{\alpha}({\mathbf{H}})$ tends to the null vector $\mathbf{0}$ as  $\mathbf{H}$ goes to the null matrix $\mathbf{0}_d$:
    \begin{enumerate}[a.]
\item \textit{interior region} is the largest one inside the interior of $\mathbb{T}_{d}$ in order to contain at least 95 percent of observations, and it is denoted by $\mathbb{T}^{\alpha(\mathbf{H}),I}_{d}$;
       \item \textit{boundary regions} representing the complementary of $\mathbb{T}^{\alpha(\mathbf{H}),I}_{d}$ in $\mathbb{T}_{d}$, and it is denoted by $\mathbb{T}_{d}^{\alpha(\mathbf{H}),B}$ which could be empty; recall that $\mathbb{T}_{d} = \mathbb{T}^{\alpha(\mathbf{H}),I}_{d} \cup \mathbb{T}^{\alpha(\mathbf{H}),B}_{d}$ and 
   $\mathbb{T}^{\alpha(\mathbf{H}),I}_{d} \cap \mathbb{T}^{\alpha(\mathbf{H}),B}_{d}=\emptyset$.
\end{enumerate}
Since $\mathbb{T}_{d}\left(\subseteq \mathbb{R}^{d}\right)$ might have each one of its $d$ convex components as unbounded, partially bounded or totally bounded interval as in (\ref{Examplesupport}), there is only one $\mathbb{T}^{\alpha(\mathbf{H}),I}_{d}$ but 
\begin{equation}
N_{d} = 1^{d_\infty}2^{d_z}3^{d_{uw}} - 1 \label{bounregions}
\end{equation}
boundary subregions of $\mathbb{T}_{d}^{\alpha(\mathbf{H}),B}$. In the below Section~\ref{sssec:Modversionbeta} an illustration is provided for $d = 2 $ with $\mathbb{T}_{2} = \left[0, 1\right]\times\left[0, 1\right]$ and, therefore, the number of boundary subregions (\ref{bounregions}) is $N_{2} = 3^2 - 1 = 8$.

\textbf{Second step}. Changing the general associated kernel $K_{\theta(\mathbf{x}, \mathbf{H})}$ into its modified version $K_{\widetilde{\theta}(\mathbf{x}, \mathbf{H})}$; that leads to replace the couple $\left(\mathbf{a}_{\theta}(\mathbf{x},\mathbf{H}), \mathbf{B}_{\theta}(\mathbf{x},\mathbf{H})\right)$ into $\left(\mathbf{a}_{\widetilde{\theta}}(\mathbf{x},\mathbf{H}),  \mathbf{B}_{\widetilde{\theta}}(\mathbf{x},\mathbf{H})\right)$ with $\mathbf{a}_{\widetilde{\theta}}(\mathbf{x},\mathbf{H}) =  \mathbf{a}_{\widetilde{\theta}_B}(\mathbf{x},\mathbf{H})\mathds{1}_{\mathbb{T}_{d}^{\alpha(\mathbf{H}),B}}(\mathbf{x})$ because  $\mathbf{a}_{\widetilde{\theta}_I}(\mathbf{x},\mathbf{H})\mathds{1}_{\mathbb{T}_{d}^{\alpha(\mathbf{H}),I}}(\mathbf{x}) = \mathbf{0}$ in the interior, and $\mathbf{B}_{\widetilde{\theta}}(\mathbf{x},\mathbf{H}) = \mathbf{B}_{\widetilde{\theta}_I}(\mathbf{x},\mathbf{H})\mathds{1}_{\mathbb{T}_{d}^{\alpha(\mathbf{H}),I}}(\mathbf{x}) +  \mathbf{B}_{\widetilde{\theta}_B}(\mathbf{x},\mathbf{H})\mathds{1}_{\mathbb{T}_{d}^{\alpha(\mathbf{H}),B}}(\mathbf{x})$.
This modified associated kernel is such that, for any fixed bandwidth matrix $\mathbf{H}$,
\begin{equation}
 \widetilde{\theta}(\mathbf{x},\mathbf{H})=  \left\{\begin{array}{lll}\begin{aligned}
  \widetilde{\theta}_{I}(\mathbf{x},\mathbf{H}):  \:\:\: \mathbf{a}_{\widetilde{\theta}_I}(\mathbf{x},\mathbf{H}) = \mathbf{0}&   &   & \textrm{if}\: \mathbf{x}\in\mathbb{T}^{\alpha(\mathbf{H}),I}_d\\
  \widetilde{\theta}_{B}(\mathbf{x},\mathbf{H})~~~~~~~~~~~~~~~~~~~~\:\:\:\:\:\:&   &   &\textrm{if}\: \mathbf{x}\in\mathbb{T}^{\alpha(\mathbf{H}),B}_d
\end{aligned}
\end{array}\right.\label{systhetatild}
\end{equation}
must be continuous on $\mathbb{T}_d$ and constant on $\mathbb{T}^{\alpha(\mathbf{H}),B}_d$. 

\begin{Pro}\label{Proprs}
The function $K_{\widetilde{\theta}(\mathbf{x},\mathbf{H})}$ on its support $\mathbb{S}_{\widetilde\theta(\mathbf{x},\mathbf{H})}=\mathbb{S}_{\theta(\mathbf{x},\mathbf{H})}$ and obtained from (\ref{systhetatild}) is also a general associated kernel.
\end{Pro}

\textbf{Proof.} Since $K_{\theta(\mathbf{x}, \mathbf{H})}$ is a general associated kernel for all $\mathbf{x}\in\mathbb{T}_{d} = \mathbb{T}^{\alpha(\mathbf{H}),I}_{d} \cup \mathbb{T}^{\alpha(\mathbf{H}),B}_{d}$, one gets the first condition of Definition \ref{Defkern} because
$\mathbb{S}_{\theta(\mathbf{x}, \mathbf{H})}=\mathbb{S}_{\widetilde{\theta}(\mathbf{x}, \mathbf{H})}$.
According to Proposition~\ref{ProTyp-Ass} it follows that, for a given random variable $\mathcal{Z}_{\widetilde{\theta}(\mathbf{x}, \mathbf{H})}$ with pdf
$K_{\widetilde{\theta}(\mathbf{x}, \mathbf{H})}$, we obtain the last two conditions of Definition~\ref{Defkern} as
$$\mathbb{E}\left(\mathcal{Z}_{\widetilde{\theta}(\mathbf{x}, \mathbf{H})}\right) = \mathbf{x}+
\mathbf{a}_{\widetilde{\theta}}(\mathbf{x},\mathbf{H}) \textrm{ and } \mbox{Cov}\left(\mathcal{Z}_{\widetilde{\theta}(\mathbf{x},\mathbf{H})}\right) = \mathbf{B}_{\widetilde{\theta}}(\mathbf{x},\mathbf{H}).$$
Both quantities $\mathbf{a}_{\widetilde{\theta}}(\mathbf{x}, \mathbf{H})$ and $\mathbf{B}_{\widetilde{\theta}}(\mathbf{x}, \mathbf{H})$ tend, respectively, to $\mathbf{0}$ and $\mathbf{0}_d$ as $\mathbf{H}$ goes to $\mathbf{0}_d$. In particular, from (\ref{systhetatild}) we have
$\mathbf{a}_{\widetilde{\theta}_{I}}(\mathbf{x}, \mathbf{H})=\mathbf{0}.\blacksquare$

Similar to both (\ref{Asskern}) and (\ref{AssCKernest}), the \textit{modified associated kernel estimator} $\widetilde{f}_{n}$ using $K_{\widetilde{\theta}(\mathbf{x},\mathbf{H})}$ is then defined by
\begin{equation}
 \widetilde f_{n}(\mathbf{x})=\dfrac{1}{n} \displaystyle \sum_{i=1}^{n}K_{\widetilde{\theta}(\mathbf{x},\mathbf{H})} \left(\mathbf{X}_{i} \right).\label{Kernest2}
\end{equation}
The following result gives only in the interior $\mathbb{T}^{\alpha(\mathbf{H}),I}_{d}$ of $\mathbb{T}_{d}$  the pointwise expressions of the bias and the variance of $\widetilde f_{n}$. Of course, the corresponding expressions in the boundary regions $\mathbb{T}^{\alpha(\mathbf{H}),B}_{d}$ are tedious to write with respect to the $N_d$ (\ref{bounregions}) boundary situations from (\ref{systhetatild}).

\begin{Pro} \label{Prop}
 Let $\widehat{f}_{n}$ and  $\widetilde{f}_{n}$ be the multivariate associated  kernel estimators of $f$ defined in (\ref{AssCKernest}) and (\ref{Kernest2}) respectively. Then, for $\mathbf{x}\in\mathbb{T}^{\alpha(\mathbf{H}),I}_{d}$ as in (\ref{systhetatild}):
\begin{equation}
 \Bias\left\{\widetilde{f}_{n}(\mathbf{x})\right\}
=  \dfrac{1}{2}\trace\left(\mathbf{B}_{\widetilde{\theta}_I}(\mathbf{x}, \mathbf{H}) \nabla^{2}f\left(\mathbf{x}\right)\right) +  \oldstylenums{0}\left\{\trace\left(\mathbf{H}^{2}\right)\right\} \label{modbias} 
\end{equation}
and
\begin{equation*}
  \Var\left\{\widetilde{f}_{n}(\mathbf{x})\right\}  \simeq \Var\left\{\widehat{f}_{n}(\mathbf{x})\right\} \label{varm} ~as~ n\to\infty.
\end{equation*}
\end{Pro}

\textbf{Proof.} We have the first result (\ref{modbias}) by replacing in (\ref{biaisg}) $\widehat{f}_{n}$, $\mathbf{a}_{\theta}$ and $\mathbf{B}_{\theta}$ by $\widetilde{f}_{n}$, $\mathbf{a}_{\widetilde\theta_{I}}$ and $\mathbf{B}_{\widetilde\theta_{I}}$, respectively.
For the last result, considering (\ref{var}) it is sufficient to show that $$\left\|K_{\theta(\mathbf{x},\mathbf{H})}\right\|^{2}_{2}\simeq
\left\|K_{\widetilde\theta(\mathbf{x},\mathbf{H})}\right\|^{2}_{2}\textrm{ as }\mathbf{H}\rightarrow \mathbf{0}_d.$$
Since $K_{\theta(\mathbf{x},\mathbf{H})}$ and $K_{\widetilde\theta(\mathbf{x},\mathbf{H})}$ are general associated kernels of the same type $K_{\theta}$ with $r_{2}=r_{2}(K_{\theta})$ and $\widetilde r_{2}=\widetilde r_{2}(K_{\widetilde\theta})$, then there exists a common   largest positive real number $r^*_{2}=r^*_{2}(K_{\theta})$  such that
\begin{equation*}
\left\|K_{\theta(\mathbf{x},\mathbf{H})}\right\|^{2}_{2}\lesssim c_{2}(\mathbf{x})(\det\mathbf{H})^{-r^*_{2}}\textrm{ and } \left\|K_{\widetilde\theta(\mathbf{x},\mathbf{H})}\right\|^{2}_{2}\lesssim
\widetilde{c}_{2}(\mathbf{x})(\det\mathbf{H})^{-r^*_{2}}
\end{equation*}
with $0 < {c}_{2}(\mathbf{x}),\widetilde{c}_{2}(\mathbf{x}) < \infty$. Taking $c(\mathbf{x})=\sup\left\{c_{2}(\mathbf{x}), \widetilde{c}_{2}(\mathbf{x})\right\}$, we have $\left\|K_{\theta(\mathbf{x},\mathbf{H})}\right\|^{2}_{2}
\lesssim c(\mathbf{x})(\det\mathbf{H})^{-2r^*_{2}}$ and $\left\|K_{\widetilde\theta(\mathbf{x},\mathbf{H})}\right\|^{2}_{2}\lesssim c(\mathbf{x})(\det\mathbf{H})^{-2r^*_{2}}$. Since $c(\mathbf{x})/n(\det\mathbf{H})^{2r^*_{2}}=\oldstylenums{0}(n^{-1}(\det\mathbf{H})^{-2r^*_{2}})$ then
$\left\|K_{\widetilde\theta(\mathbf{x},\mathbf{H})}\right\|^{2}_{2}\simeq\left\|K_{\theta(\mathbf{x},\mathbf{H})}\right\|^{2}_{2}$.$\blacksquare$

Thus, we define the asymptotic expression of the MISE of $\widetilde{f}_{n}$ in the interior $\mathbb{T}^{\alpha(\mathbf{H}),I}_{d}$ as follows:
\begin{equation*}
\AMISE_{\widetilde{\theta}_I}\left(\widetilde{f}_{n}\right)  =  \int_{\mathbb{T}^{\alpha(\mathbf{H}),I}_{d}}\left(\left\{ \dfrac{1}{2}\trace\left( \mathbf{B}_{\widetilde{\theta}_I}(\mathbf{x}, \mathbf{H}) \nabla^{2}f\left(\mathbf{x}\right)\right) \right\}^{2}   + \frac{1}{n} \left\|K_{\theta(\mathbf{x}, \mathbf{H})}\right\|^{2}_{2}f(\mathbf{x})\right)d\bold{x}. \label{AMISEm}
\end{equation*}All the results in this section can be easily deduced for both cases of diagonal and Scott bandwidth matrices. The following section provides some detailed results for $d = 2$ with a bivariate beta kernel having correlation structure on $\mathbb{T}_{2} = \left[0, 1\right]\times\left[0, 1\right]$. The corresponding associated kernel is built from a technique due to \cite{S66} and using two independent univariate beta pdfs. See \cite{BL09}, \cite{KBJ10} for some other examples of bivariate types of kernels and \cite{KBJ00} in multivariate case.

\section{Bivariate beta kernel with correlation structure}
\label{sec:BivBetaKern}

This section presents the generalizable procedure of a bivariate beta kernel estimator from a bivariate beta pdf with correlation structure to its corresponding associated kernel. The standard associated kernel is built by a variant of the mode-dispersion method deduced from (\ref{eqtheta}). Then, we provide properties of both versions of the corresponding estimators.

\subsection{Type of bivariate beta kernel}
\label{ssec:Typeofbetabiv}

In order to control better the effects of correlation, we here consider a flexible type of bivariate beta kernel for which the correlation structure is introduced by \cite{S66}; see also \cite{L96}. Let us take two independent univariate beta distributions with pdfs
\begin{equation}
 g_{j}(t) = \frac{1}{\mathscr{B}(p_{j}, q_{j})}t^{p_{j}-1}(1 - t)^{q_{j}-1}\mathds{1}_{\left[0, 1\right]}(t), ~~~~j = 1, 2,\label{gunif}
\end{equation}
where $\mathscr{B}(p_{j}, q_{j}) = \int_{0}^{1}t^{p_{j}-1}(1 - t)^{q_{j} - 1}dt$ is the usual beta function with $p_{j} > 0$ and $q_{j} > 0$. Their means and variances are, respectively,
\begin{equation}
 \mu_{j}  = \frac{p_{j}}{p_{j} + q_{j}} = \mu_{j}(p_{j}, q_{j})   \quad \mbox{ and } \quad  \sigma_{j}^{2}  = \frac{p_{j}q_{j}}{(p_{j} + q_{j})^{2}(p_{j} + q_{j} + 1)} = \sigma_{j}^{2}(p_{j}, q_{j}). \label{MV}
\end{equation}
Also, $g_{j}$ are unimodal for $p_{j}\geq 1$, $q_{j}\geq 1$ and $(p_{j}, q_{j})\neq (1, 1)$, with mode and dispersion parameters:
\begin{equation}
m_{j}(p_{j}, q_{j})  = \frac{p_{j} - 1}{p_{j} + q_{j} - 2} \quad \mbox{ and } \quad    d_{j}  = \frac{1}{p_{j} + q_{j} - 2} = d_{j}(p_{j}, q_{j}). \label{MD}
\end{equation}

The corresponding pdf (or type of kernel) of the bivariate beta-Sarmanov with correlation and from $g_j$ of (\ref{gunif}) is then denoted by $g_{\theta}$ ($=BS_{\theta}$) and defined as:
\begin{equation}
 g_{\theta}(\mathbf{v}) = g_{1}(v_{1})g_{2}(v_{2})\left[1 + \rho\times\dfrac{v_{1} - \mu_{1}(p_{1}, q_{1})}{\sigma_{1}(p_{1}, q_{1})}\times\dfrac{v_{2} - \mu_{2}(p_{2}, q_{2})}{\sigma_{2}(p_{2}, q_{2})}\right]\mathds{1}_{\left[0, 1\right]\times\left[0, 1\right]}(\mathbf{v}), \label{gbiv}
\end{equation}
with $\mathbf{v} = (v_{1}, v_{2})^{\top}$ and $\theta := \theta(p_{1}, q_{1}, p_{2}, q_{2}, \rho) \in \Theta \subseteq \mathbb{R}^{5}$. Depending on $p_j$ and $q_j$, the correlation parameter $\rho=\rho(p_{1}, q_{1}, p_{2}, q_{2})$ belongs to the following interval 
\begin{equation}
\left[-\varepsilon, \varepsilon'\right] \subset [-1, 1], \label{rho}
\end{equation}
with nonnegative reals
\begin{equation*}
\varepsilon  = \left(\max_{v_1,v_2} \left\{\dfrac{v_{1} - \mu_{1}(p_{1}, q_{1})}{\sigma_{1}(p_{1}, q_{1})}\times\dfrac{v_{2} - \mu_{2}(p_{2}, q_{2})}{\sigma_{2}(p_{2}, q_{2})}\right\}\right)^{-1}\end{equation*}
and 
\begin{equation*}
 \varepsilon' = \left\lvert\left(\min_{v_1,v_2} \left\{\dfrac{v_{1} - \mu_{1}(p_{1}, q_{1})}{\sigma_{1}(p_{1}, q_{1})}\times\dfrac{v_{2} - \mu_{2}(p_{2}, q_{2})}{\sigma_{2}(p_{2}, q_{2})}\right\}\right)^{-1}\right\lvert.
\end{equation*}
Thus, the mean vector and covariance matrix of $g_{\theta}$ are, respectively,
\begin{equation*}
 \boldsymbol{\mu} = (\mu_{1}, \mu_{2})^{\top}   \quad \mbox{ and } \quad  \boldsymbol{\Sigma}  = \left(
                \begin{array}{cc}
                  \sigma_{1}^{2} &  \sigma_{1} \sigma_{2}\rho \\
                    \sigma_{1} \sigma_{2}\rho & \sigma_{2}^{2} \\
                \end{array}
    \right).
\end{equation*}
\begin{figure}[hbtp]
  \mbox{
\subfloat[(a0): $\rho = 0$]{\includegraphics[width=205pt,height=165pt,scale=0.95]{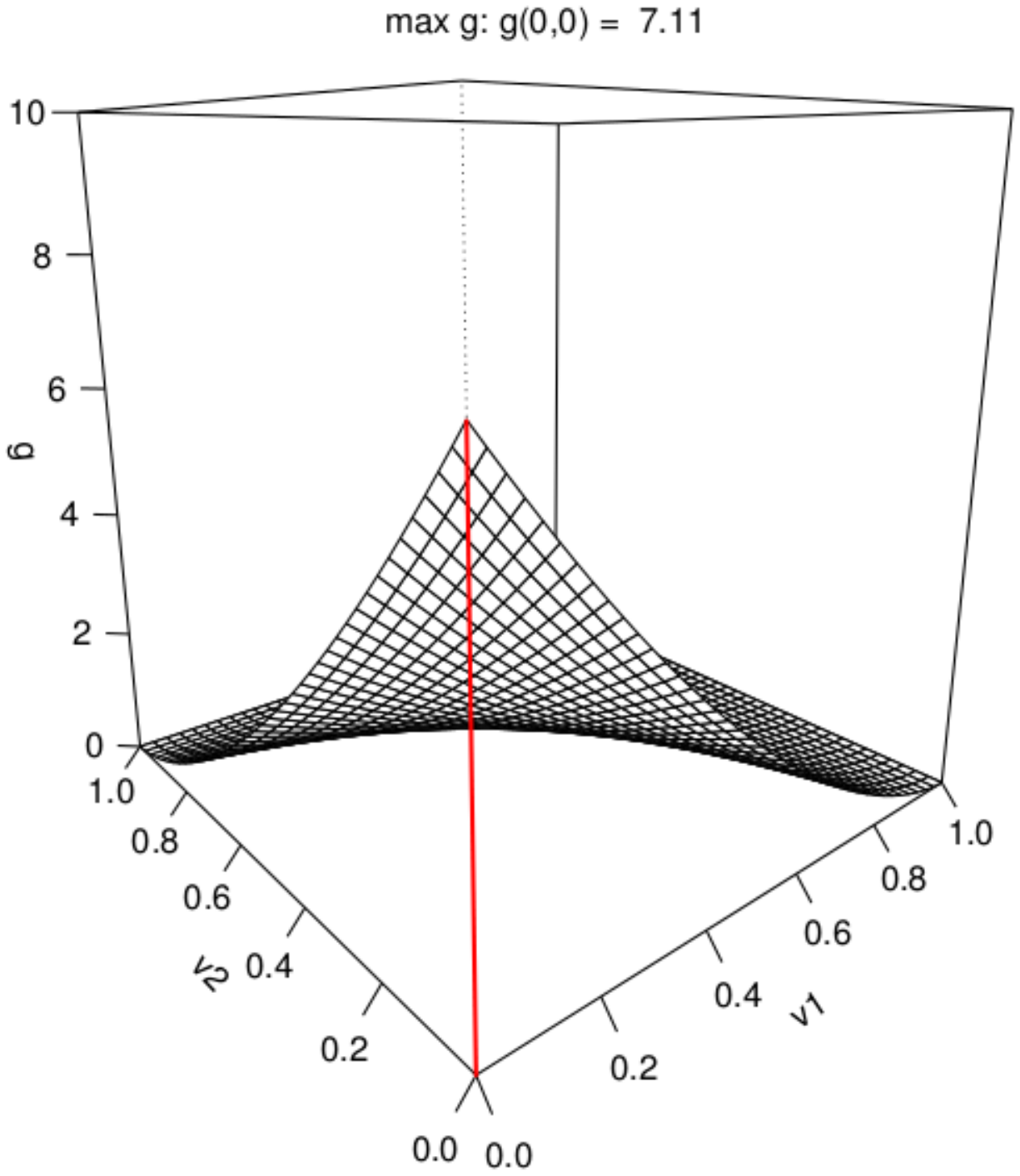} }\:\:\:
\subfloat[(a1): $\rho = 0.18$]{\includegraphics[width=205pt,height=165pt,scale=0.95]{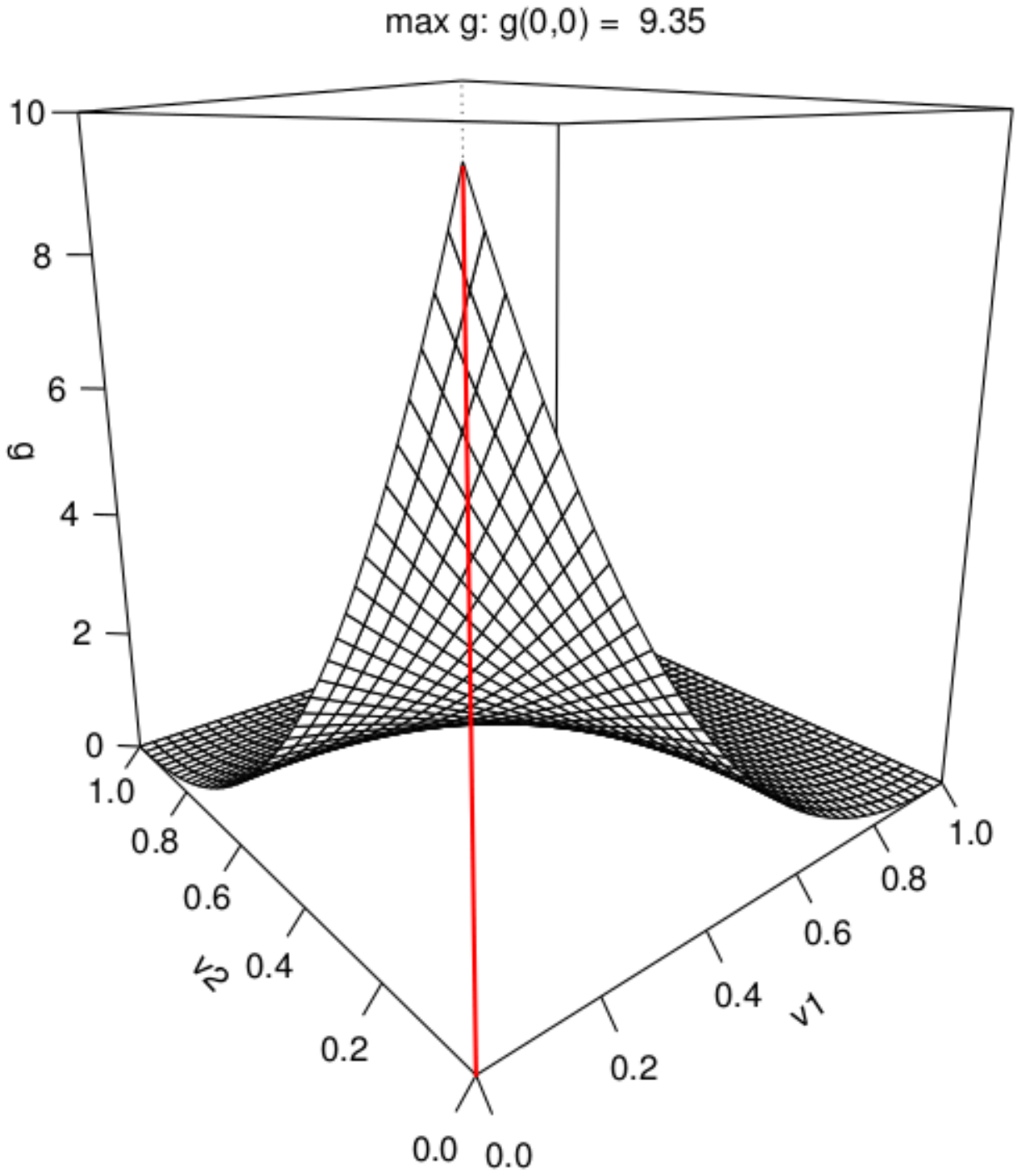} } 
}
  \mbox{
\subfloat[(b0): $\rho = 0$]{\includegraphics[width=205pt,height=165pt,scale=0.95]{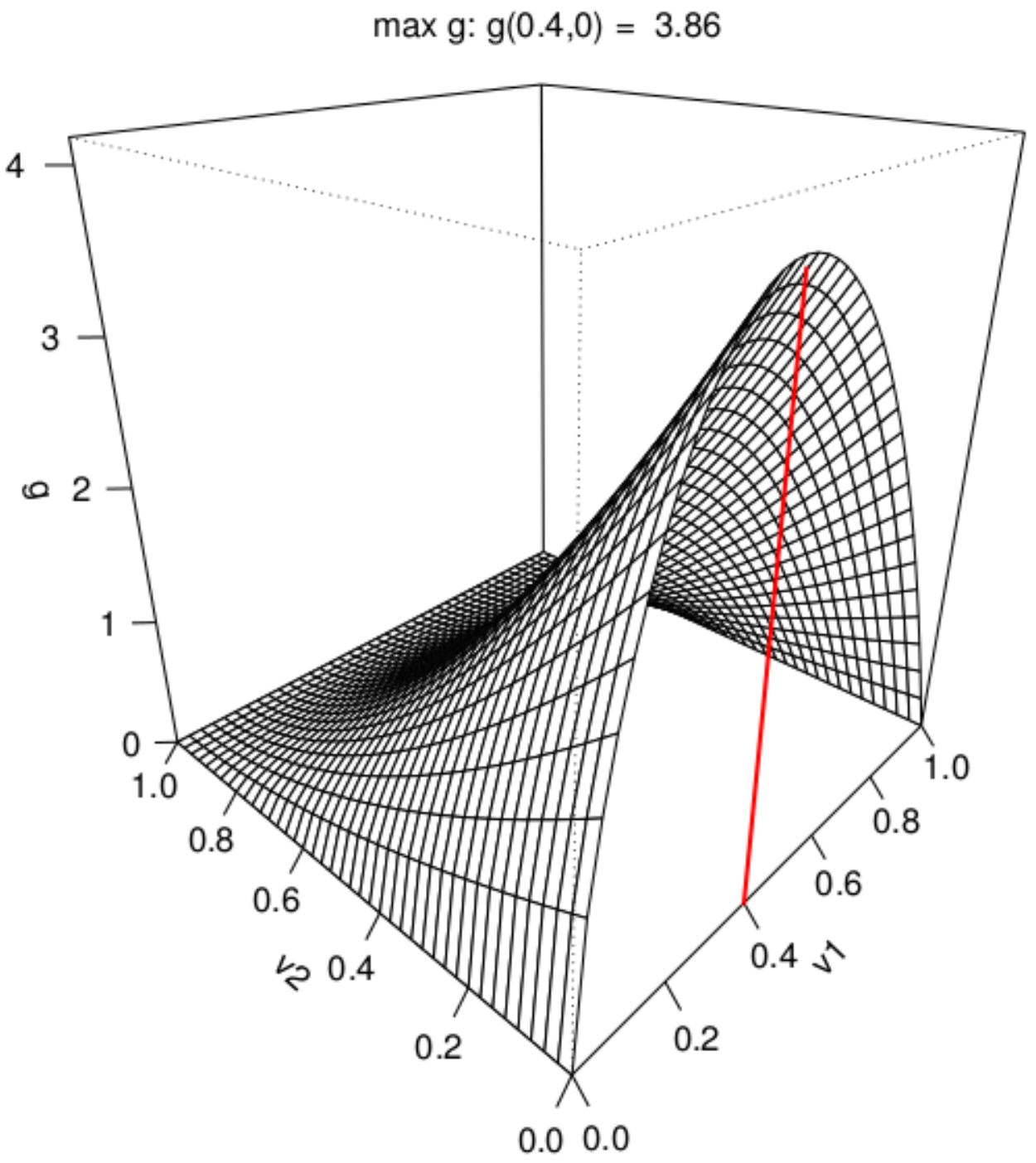} }\:\:\: 
\subfloat[(b1): $\rho = 0.12$]{\includegraphics[width=205pt,height=165pt,scale=0.95]{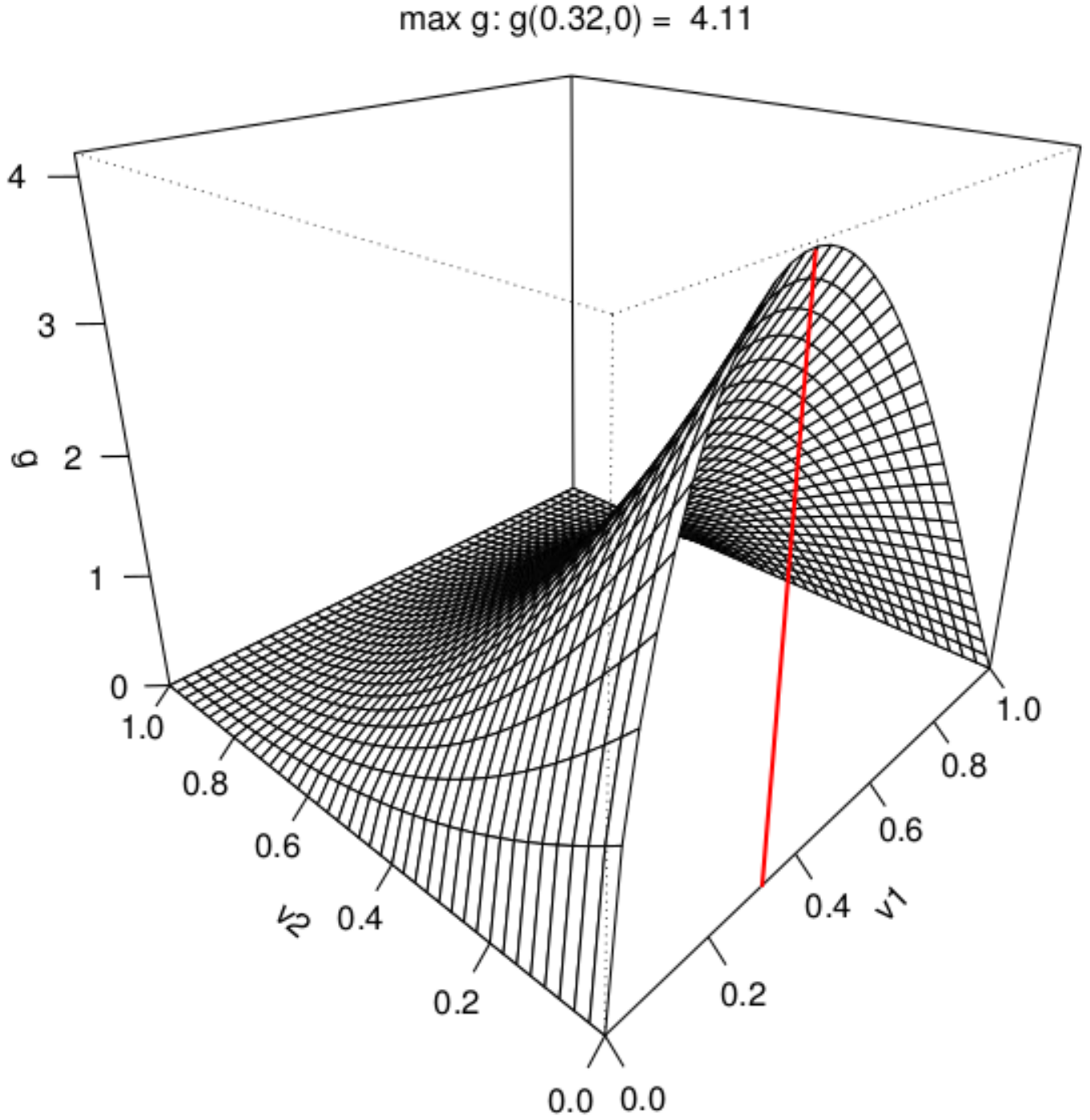} } 
}
  \mbox{
\subfloat[(c0): $\rho = 0$]{\includegraphics[width=205pt,height=165pt,scale=0.95]{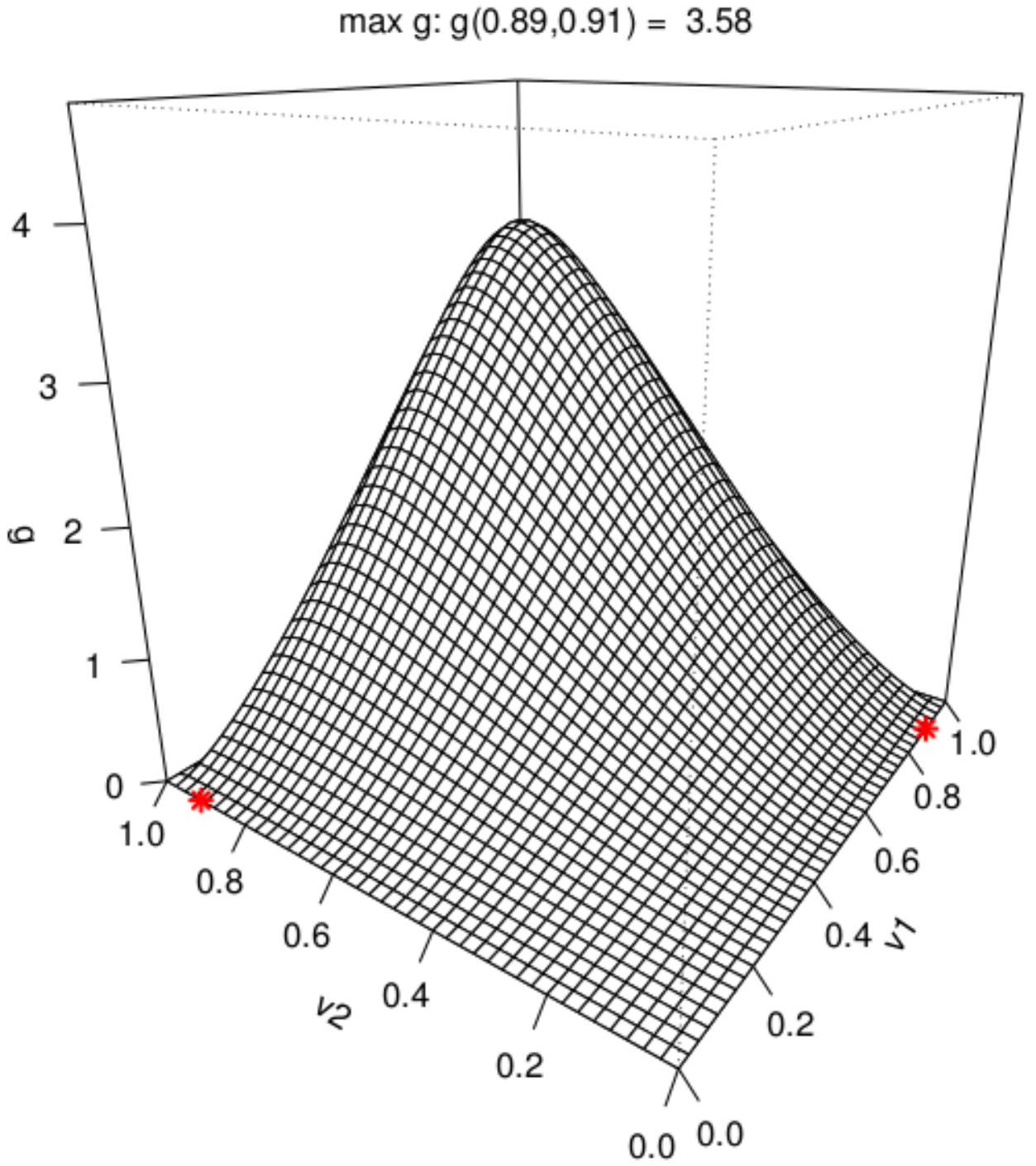} }\:\:\: 
\subfloat[(c1): $\rho = 0.19$]{\includegraphics[width=205pt,height=165pt,scale=0.95]{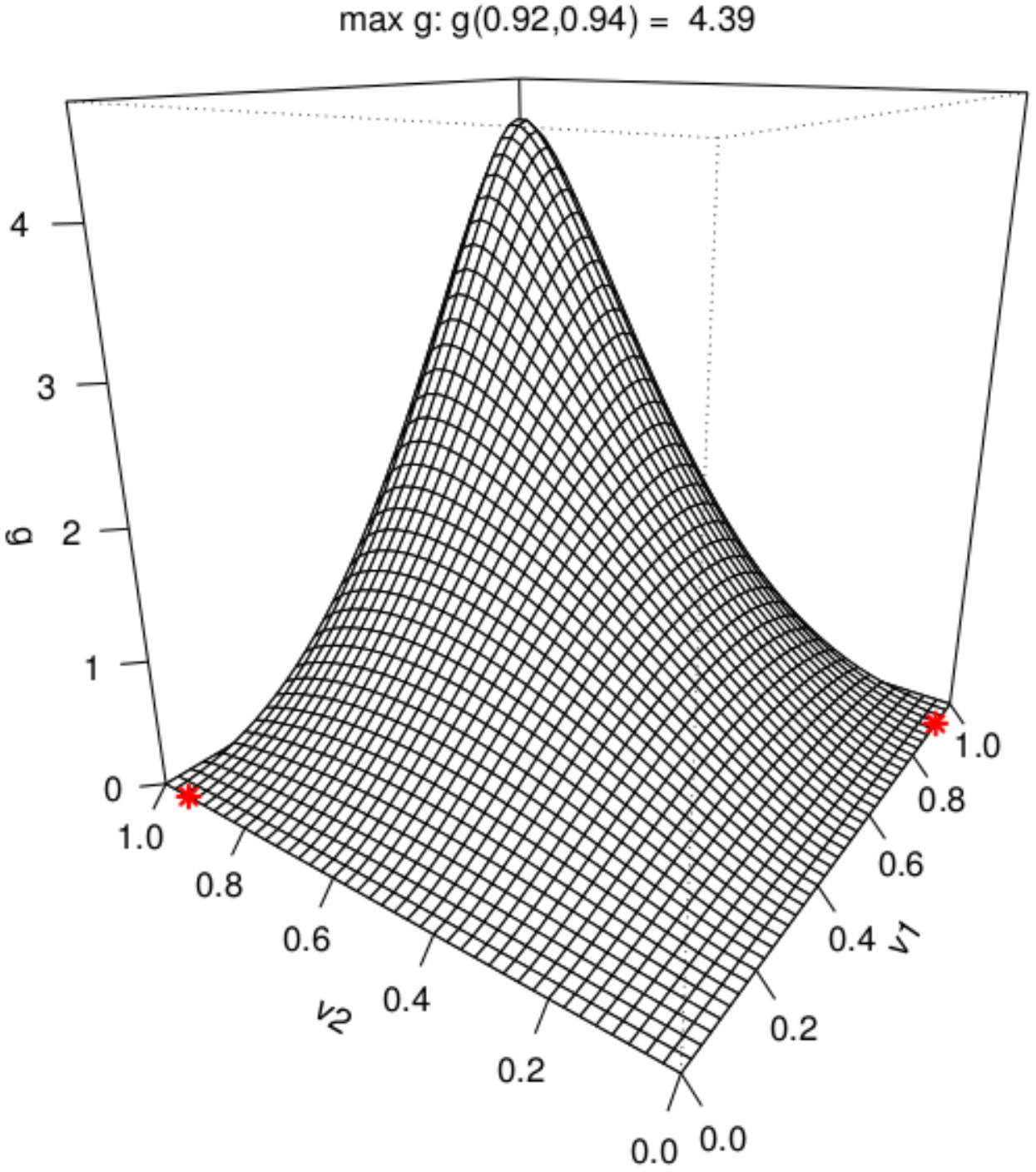} }
}
\caption{Some shapes of the bivariate beta-Sarmanov type (\ref{gbiv}) with different effects of correlations on the unimodality ((a): $p_{1} = p_{2} = 1, q_{1} = q_{2} = 8/3$; (b): $p_{1} = 5/3, p_{2} = 1, q_{1} = 2, q_{2} = 8/3$; (c): $p_{1} = 149/60, p_{2} = 151/60, q_{1} = 71/60, q_{2} = 23/20$).}
\label{figdist}
\end{figure}
The unimodality of $g_{\theta}$ in (\ref{gbiv}) also occurs for $p_{j} \geq 1$, $q_{j} \geq 1$ and $(p_{j}, q_{j}) \neq (1, 1)$ with $j = 1, 2$. However, the corresponding mode vector $\mathbf{m}_{\theta}=\mathbf{m}(p_{1},p_{2},q_{1},q_{2},\rho) =: \mathbf{m}_{\rho}$ of $g_{\theta}$ does not have an explicit expression; but, we numerically verified that this mode vector $\mathbf{m}_{\rho}$ is slightly shifted with respect to the modal vector (\ref{MD}) of the two independent margins that we denote by $\mathbf{m}_{0} = (m_{1}(p_{1}, q_{1}), m_{2}(p_{2}, q_{2}))^{\top}$ for $\rho = 0$. 

Figure~\ref{figdist} illustrates some different effects of both null and positive correlations (\ref{rho}) on the unimodality of (\ref{gbiv}). The negative correlations will show the opposite effects in terms of positions according to the modal vector $\mathbf{m}_{0}$ for $\rho = 0$. In other words, the correlation parameter $\rho$ enables the pdf $g_{\theta}$ to reach points which are inaccessible with the null correlation. The parameter values of both univariate beta (\ref{gunif}) used for Figure~\ref{figdist} produce the following intervals (\ref{rho}) of correlation:
\begin{itemize}
 \item $\rho \in \left[-0.100, 0.210\right]$ if $p_{1} = p_{2} = 1, q_{1} = q_{2} = 8/3$ for an angle;
 \item $\rho \in \left[-0.120, 0.143 \right]$ if $p_{1} = 5/3, p_{2} = 1, q_{1} = 2, q_{2} = 8/3$ for an edge; 
 \item $\rho \in \left[-0.008,  0.214\right]$ if $p_{1} = 149/60, p_{2} = 151/60, q_{1} = 71/60, q_{2} = 23/20$ for the interior.
\end{itemize}

Concerning the dispersion matrix $\mathbf{D}_{\rho}$ of the bivariate beta-Sarmanov type we consider
 \begin{equation}
  \mathbf{D}_{\rho} = \left(
                \begin{array}{cc}
                  d_{1} & \left(d_{1}d_{2}\right)^{1/2}\rho \\
                   \left(d_{1}d_{2}\right)^{1/2}\rho & d_{2} \\
                \end{array} \label{dispersion}
    \right),
\end{equation}
where $d_{1}$ and $d_{2}$ are the dispersion parameters (\ref{MD}) of margins and $\rho$ the correlation parameter. This dispersion matrix is the analogue of the covariance one. Since we do not have a closed expression of the modal vector $\mathbf{m}_{\rho}$, we cannot use the bivariate mode-dispersion method (\ref{eqtheta}) for a construction of the bivariate beta-Sarmanov kernel.

\subsection{Bivariate beta-Sarmanov kernel}
\label{ssec:BivbetaSarmanov}

From the previous section, the standard version of the bivariate beta-Sarmanov kernel is here constructed by using the modal vector $\mathbf{m}_{0}= (m_{1}(p_{1}, q_{1}), m_{2}(p_{2}, q_{2}))^{\top}$ of $\rho = 0$ instead of $\mathbf{m}_{\rho}$ as a variant of the mode-dispersion method (\ref{eqtheta}). This choice will be compensated in the bandwidth matrix $\mathbf{H}$ connected to the complete dispersion matrix $\mathbf{D}_{\rho}$ with correlation structure (\ref{dispersion}). 

Indeed, solving $\left(\theta(\mathbf{m}_{0}, \mathbf{D}_{\rho})\right)^{\top} = (\mathbf{x}, \vech\mathbf{H})^{\top}$ in the sense of $\mathbf{m}_{0} = \mathbf{x}$ and $\mathbf{D}_{\rho}=\mathbf{H}$ leads to the new reparametrization of $g_{\theta}$ of (\ref{gbiv}) from $\theta = \theta(p_{1}, q_{1}, p_{2}, q_{2}, \rho) \subseteq \mathbb{R}^{5}$ into
\begin{equation}
\theta(\mathbf{x}, \mathbf{H}) = \left(\frac{x_{1}}{h_{11}} + 1, \frac{1 - x_{1}}{h_{11}} + 1, \frac{x_{2}}{h_{22}} + 1, \frac{1 - x_{2}}{h_{22}} + 1, \frac{h_{12}}{(h_{11}h_{22})^{1/2}}\right)^{\top}, ~~\forall \mathbf{x} = \begin{pmatrix}x_{1} \\ x_{2}\end{pmatrix}, \mathbf{H}=\begin{pmatrix}h_{11}  & h_{12} \\ h_{12} & h_{22}\end{pmatrix}. \label{theta}
\end{equation}
Rewriting (\ref{MV}) in terms of (\ref{theta}), the means $\mu_{j}(p_{j}, q_{j})$ and variances $\sigma_{j}^{2}(p_{j}, q_{j})$ of the univariate beta pdfs become
\begin{equation*}
 \widetilde{\mu}_{j} = \frac{x_{j} + h_{jj}}{1 + 2h_{jj}} =   \widetilde{\mu}_{j}(x_{j}, h_{jj}) \quad \mbox{and}\quad \widetilde{\sigma}_{j}^{2} =  \frac{(x_{j} + h_{jj})(1 + h_{jj} - x_{j})}{(1 + 2h_{jj})^{2}(1+3h_{jj})}h_{jj} = \widetilde{\sigma}_{j}^{2}(x_{j}, h_{jj}).
\end{equation*}
Finally, the bivariate beta-Sarmanov kernel is defined as $BS_{\theta(\mathbf{x}, \mathbf{H})}:= g_{\theta(\mathbf{x}, \mathbf{H})}$ such that
\begin{eqnarray}
\label{betakern}
 BS_{\theta(\mathbf{x}, \mathbf{H})}(v_1,v_2) &=& \left(\frac{v_{1}^{x_{1}/h_{11}}(1-v_{1})^{(1 - x_{1})/h_{11}} }{\mathscr{B}(1 + x_{1}/h_{11}, 1 + (1 - x_{1})/h_{11})}\right)  \left(\frac{v_{2}^{x_{2}/h_{22}}(1-v_{2})^{(1 - x_{2})/h_{22}}}{\mathscr{B}(1 + x_{2} / h_{22}, 1 + (1 - x_{2})/h_{22})} \right)  \nonumber \\
&&\times \left(1 + h_{12}\times\dfrac{v_{1} - \widetilde{\mu}_{1}(x_{1}, h_{11})}{h_{11}^{1/2}\widetilde{\sigma}_{1}(x_{1}, h_{11})}\times\dfrac{v_{2} - \widetilde{\mu}_{2}(x_{2}, h_{22})}{h_{22}^{1/2}\widetilde{\sigma}_{2}(x_{2}, h_{22})}\right)\mathds{1}_{\left[0, 1\right]^2}(v_1,v_2),
\end{eqnarray}
with the constraints
\begin{equation}
h_{12} \in \left[-\beta, \beta^{\prime}\right] \cap \left(-\sqrt{h_{11}h_{22}}\:,\sqrt{h_{11}h_{22}}\right), \label{h12}
\end{equation}
\begin{equation*}
\beta  = \left(\max_{v_1,v_2} \left\{\dfrac{v_{1} - \widetilde{\mu}_{1}(x_{1}, h_{11})}{h_{11}^{1/2}\widetilde{\sigma}_{1}(x_{1}, h_{11})}\times\dfrac{v_{2} - \widetilde{\mu}_{2}(x_{2}, h_{22})}{h_{22}^{1/2}\widetilde{\sigma}_{2}(x_{2}, h_{22})}\right\}\right)^{-1}
\end{equation*}
 and
\begin{equation*}
\beta^{\prime}  = \left\lvert\left(\min_{v_1,v_2} \left\{\dfrac{v_{1} - \widetilde{\mu}_{1}(x_{1}, h_{11})}{h_{11}^{1/2}\widetilde{\sigma}_{1}(x_{1}, h_{11})}\times\dfrac{v_{2} - \widetilde{\mu}_{2}(x_{2}, h_{22})}{h_{22}^{1/2}\widetilde{\sigma}_{2}(x_{2}, h_{22})}\right\}\right)^{-1}\right\lvert.
\end{equation*}

The first interval $\left[-\beta, \beta^{\prime}\right]$ of (\ref{h12}) is the equivalent $\left[-\varepsilon, \varepsilon'\right]$ in (\ref{rho}) and the second one $(-\sqrt{h_{11}h_{22}}\:,\sqrt{h_{11}h_{22}})$ of (\ref{h12}) is due to the constraints from the bandwidth matrix $\mathbf{H}$, which is symmetric and positive definite. In practice, one often has: 
$\left[-\beta, \beta^{\prime}\right]\subset (-\sqrt{h_{11}h_{22}}\:,\sqrt{h_{11}h_{22}})$. Of course, 
the beta-Sarmanov kernel $BS_{\theta(\mathbf{x}, \mathbf{H})}$ satisfies Definition~\ref{Defkern} of associated kernel:
\begin{small}
\begin{eqnarray}
\mathbb{S}_{\theta(\mathbf{x}, \mathbf{H})} = \left[0, 1\right]\times\left[0, 1\right],~~ \mathbf{a}_{\theta}(\mathbf{x}, \mathbf{H}) = (a_{\theta 1}, a_{\theta 2})^{\top}~~ \mbox{with} ~ a_{\theta j} = \frac{(1 - 2x_{j})h_{jj}}{1 + 2h_{jj}} = a_{\theta j}(x_{j}, h_{jj}), \nonumber \\
  \mathbf{B}_{\theta}(\mathbf{x},  \mathbf{H}) = (b_{\theta ij})_{i, j = 1, 2} ~\mbox{with }
 b_{\theta jj} = \widetilde{\sigma}_{j}^{2}(x_{j}, h_{jj}) ~ \mbox{and} ~ b_{\theta 12} = \frac{h_{12}}{\sqrt{h_{11}h_{22}}}\widetilde{\sigma}_{1}(x_{1}, h_{11})\widetilde{\sigma}_{2}(x_{2}, h_{22}). \label{condBS}
\end{eqnarray}
\end{small}
Table~\ref{Eqvfigkern} shows some effects of the correlation parameter $h_{12}$ in $BS_{\theta(\mathbf{x}, \mathbf{H})}$ on the modal vector and maximum values, which are obtained by using corresponding values of $(p_{1}, q_{1}, p_{2}, q_{2}, \rho)$ for Figure~\ref{figdist}. 

\begin{table}[h!]
\begin{center}
\begin{tabular}{lcccl}
\hline \hline
Position &  $\mathbf{x}$ &Interval of $h_{12}$ & $\mathbf{H}$ & Maximum value of $BS_{\theta}$ \\
\hline
\multirow{3}{*}{\textit{Angle}}  &\multirow{3}{*}{(0, 0) }& \multirow{3}{*}{[-0.040, 0.128]}&
$\mathbf{Diag}\left(0.600,   0.600\right)$
&  $BS_{\theta(\mathbf{x}, \mathbf{H})}(0, 0)$ = 7.11 \\
 & & & $\begin{pmatrix}
0.600 & 0.128 \\
0.128& 0.600
\end{pmatrix}$&$BS_{\theta(\mathbf{x}, \mathbf{H})}(0,0)$ = 9.77\\ \hline
 \multirow{3}{*}{\parbox{1.5cm}{\textit{Edge} } } &\multirow{3}{*}{(0.40, 0.00)} & \multirow{3}{*}{[-0.050, 0.104]}&
$\mathbf{Diag}\left(0.600,   0.600\right)$
& $BS_{\theta(\mathbf{x}, \mathbf{H})}(0.40,0.00)$ = 3.86 \\
 
 & & &$\begin{pmatrix}
0.600 & 0.104 \\
0.104 & 0.600
\end{pmatrix}$&$BS_{\theta(\mathbf{x}, \mathbf{H})}(0.32,0.00)$ = 4.11\\ \hline
 \multirow{3}{*}{\textit{Interior}}  & \multirow{3}{*}{(0.89, 0.91)} &  \multirow{3}{*}{[-0.060, 0.128]}&
$\mathbf{Diag}\left(0.612,   0.600\right)$
& $BS_{\theta(\mathbf{x}, \mathbf{H})}(0.89,0.91)$ = 3.58 \\
 & & &$\begin{pmatrix}
0.612 & 0.123 \\
0.123 & 0.600
\end{pmatrix}$&$BS_{\theta(\mathbf{x}, \mathbf{H})}(0.92,0.94)$ = 4.46\\
 \hline \hline
\end{tabular}
\caption{The corresponding values of Figure~\ref{figdist} for the bivariate beta-Sarmanov kernel (\ref{betakern}).}
\label{Eqvfigkern}
\end{center}
\end{table}

\subsection{Bivariate beta-Sarmanov kernel estimators}
\label{ssec:Betasarmanovkernel}

\subsubsection{Standard version of the estimator}
\label{sssec:standbeta}

In this particular case, the beta-Sarmanov kernel estimator
\begin{equation}
  \widehat{f}_{n}(\mathbf{x}) = \dfrac{1}{n} \displaystyle \sum_{i=1}^{n}BS_{\theta(\mathbf{x}, \mathbf{H})}(\mathbf{X}_{i}),~~ \forall  \mathbf{x} \in \left[0, 1\right]\times\left[0, 1\right] \label{est}
\end{equation}
also satisfies Proposition~\ref{Pro2}. Table~\ref{Tablnormconst} allows to observe the effect of correlation (\ref{h12}) on the total mass  $\Lambda(n; \mathbf{H}, BS_{\theta})\neq 1$ by using four samples of simulated data.
\begin{table}[!h]
\begin{center}
\begin{tabular}{llllllll}
\hline \hline
$\Lambda(n, \mathbf{H}, BS_{\theta})$& sample 1& sample 2 & sample 3&sample 4  \\
\hline
$h_{12} = $-0.0003 &1.034019 & 0.9954256& 1.002369& 1.025671 \\
$h_{12} = 0 $&1.034078 & 0.9955653 &1.002418 & 1.025731 \\
$h_{12} = 0.0004$ & 1.034157 &0.9957517&  1.002483 & 1.025811 \\
 \hline \hline
\end{tabular}
\caption{Some values of $\Lambda(n; \mathbf{H}, BS_{\theta})$ for $h_{11} = 0.10$, $h_{22} = 0.07$ and $n=1000$.} \label{Tablnormconst}
\end{center}
\end{table}
Fixing $\mathbf{x} = (x_{1}, x_{2})^{\top}$ in $[0, 1]\times [0, 1]$, the pointwise bias is written as
\begin{eqnarray*}
 \Bias\left\{\widehat{f}_{n}(\mathbf{x})\right\} &=&
  a_{\theta 1}\frac{\partial f }{\partial x_{1}}(\mathbf{x}) + a_{\theta 2}\frac{\partial f }{\partial x_{2}}(\mathbf{x}) +
\frac{1}{2} \left\{ (a_{\theta 1}^{2} + b_{\theta 11})\frac{\partial^{2} f}{\partial x_{1}\partial x_{1}}(\mathbf{x}) + 2(a_{\theta 1}a_{\theta 2} + b_{\theta 12})\frac{\partial^{2} f}{\partial x_{1}\partial x_{2}}(\mathbf{x})\right. \nonumber\\
&& \left. + ~(a_{\theta 2}^{2} + b_{\theta 22})\frac{\partial^{2} f}{\partial x_{2}\partial x_{2}}(\mathbf{x}) \right\} + \oldstylenums{0}(h_{11}^{2} + 2h_{12}^{2} + h_{22}^{2});  \label{biais}
\end{eqnarray*} 
and, the pointwise variance is
\begin{equation*}
  \Var\left\{\widehat{f}_{n}(\mathbf{x})\right\}  = \frac{1}{n} \left\|BS_{\theta(\mathbf{x}, \mathbf{H})}\right\|^{2}_{2}f(\mathbf{x}) + \mbox{  } \oldstylenums{0}\left(n^{-1}\;(\mbox{det}\mathbf{H})^{-r_{2} }\right),
\end{equation*}
with
\begin{eqnarray*}
\left\|BS_{\theta(\mathbf{x}, \mathbf{H})}\right\|^{2}_{2} &=& \frac{\mathscr{B}\left(1 + 2x_{1}/h_{11}, 1 + (1 - x_{1})/h_{11}\right)\mathscr{B}\left(1 + 2x_{2}/h_{22}, 1 + (1 - x_{2})/h_{22}\right) }{\left\{\mathscr{B}\left(1 + x_{1}/h_{11}, 1 + (1 - x_{1})/h_{11}\right)\mathscr{B}\left(1 + x_{2}/h_{22}, 1 + (1 - x_{2})/h_{22}\right)\right\}^{2}} \nonumber \\
&& \times~\left\{  1 + \frac{(2x_{1} - 1)(2x_{2} - 1)}{2(1 + h_{11})(1 + h_{22})} \left(\frac{(x_{1} + h_{11})^{-1}(1 + 3h_{11})(x_{2} + h_{22})^{-1}(1 + 3h_{22})}{(1 - x_{1} + h_{11})(1 - x_{2} + h_{22})}\right)^{1/2}h_{12} \right. \nonumber \\
&&  \left. +~  \left(\frac{2x_{1} + h_{11}}{h_{11}(1 + 2h_{11})^{-1}} + \frac{2(x_{1} + h_{11})(1 + h_{11})}{h_{11}^{2}(2 + 3h_{11})^{-1}}\right) 
 \left(\frac{2x_{2} + h_{22}}{h_{22}(1 + 2h_{22})^{-1}} + \frac{2(x_{2} + h_{22})(1 + h_{22})}{h_{22}^{2}(2 + 3h_{22})^{-1}}\right)\right. \nonumber \\
&& \left. \times \left(\frac{(1 + h_{11})^{-1}(1 + 3h_{11})(1 + h_{22})^{-1}(1 + 3h_{22})}{4(1 - x_{1} + h_{11})(2 + 3h_{11})(1 - x_{2} + h_{22})(2 + 3h_{22})}\right) h_{12}^{2} \right\}. \label{normbs}
\end{eqnarray*}
Using (\ref{condBS}) the AMISE (\ref{amiseg}) becomes here
\begin{eqnarray*} \label{amise}
\AMISE\left(\widehat{f}_{n}\right)  &=&   \int_{\left[0, 1\right]\times\left[0, 1\right]}\left(\left[ a_{\theta 1}\frac{\partial f }{\partial x_{1}}(\mathbf{x}) + a_{\theta 2}\frac{\partial f }{\partial x_{2}}(\mathbf{x}) +
\frac{1}{2} \left\{ (a_{\theta 1}^{2} + b_{\theta 11})\frac{\partial^{2} f}{\partial x_{1}^2}(\mathbf{x}) \right. \right.\right.  \nonumber \\
&&\left. \left. \left.  + 2(a_{\theta 1}a_{\theta 2} + b_{\theta 12})\frac{\partial^{2} f}{\partial x_{1}\partial x_{2}}(\mathbf{x}) + (a_{\theta 2}^{2} + b_{\theta 22})\frac{\partial^{2} f}{\partial x_{2}^2}(\mathbf{x}) \right\} \right]^{2} + \frac{1}{n} \left\|BS_{\theta(\mathbf{x}, \mathbf{H})}\right\|^{2}_{2}f(\mathbf{x})\right)d\bold{x}.
\end{eqnarray*}

The bandwidth matrix is selected  by the LSCV method (\ref{Hcv}) on the set $\mathcal{D}$ of diagonal bandwidth matrices. Concerning full and Scott cases, this LSCV method is used under $\mathcal{H}_{1}$ and $\mathcal{S}_{1}$, respectively, subsets of $\mathcal{H}$ and $\mathcal{S}$ verifying  the constraint of the beta-Sarmanov kernel (\ref{h12}). Their algorithms are described below and used for numerical studies in Section~\ref{sec:Simulations}.

\subsubsection*{\textit{Algorithms of LSCV method (\ref{Hcv}) for three forms of bandwidth matrices in two dimensions ($d=2$)}}\label{par:Algorithme}

\begin{enumerate}

\item[A1.] Full bandwidth matrices. 
\begin{enumerate}[1.]
\item Choose two intervals $H_{11}$ and $H_{22}$ related to $h_{11}$ and $h_{22}$, respectively.
\item For $\delta = 1, \ldots, \ell(H_{11})$ and $\gamma = 1, \ldots, \ell(H_{22})$,
\begin{enumerate}[(a)]
 \item Compute the interval $H_{12}[\delta,\gamma]$ related to $h_{12}$ from constraints (\ref{h12}); 
 \item For $\lambda = 1, \ldots, \ell(H_{12}[\delta,\gamma])$, \\
  Compose the full bandwidth matrix $\mathbf{H}(\delta,\gamma,\lambda):=\left(h_{ij}(\delta,\gamma,\lambda)\right)_{i,j=1,2}$ with $h_{11}(\delta,\gamma,\lambda)=H_{11}(\delta)$, $h_{22}(\delta,\gamma,\lambda)=H_{22}(\gamma)$ and 
  $h_{12}(\delta,\gamma,\lambda)=H_{12}[\delta,\gamma](\lambda)$.
   \end{enumerate}
  \item Apply LSCV method on the set $\mathcal{H}_1$ of all full bandwidth matrices $\mathbf{H}(\delta,\gamma,\lambda)$.
\end{enumerate}

\item[A2.] Scott bandwidth matrices.  
\begin{enumerate}[1.]
\item Choose an interval $H$ related to $h$ and a fixed bandwidth matrix $\mathbf{H}_{0} = \left(h_{0ij}\right)_{i,j=1,2}$.
\item For $\zeta = 1, \ldots, \ell(H)$, 
\begin{enumerate}[(a)]
 \item Compute the interval $H_{012}[\zeta]$ related to $h_{012}$ from constraints (\ref{h12});
 \item For $\kappa = 1, \ldots, \ell(H_{012}[\zeta])$, \\
  Compose the given bandwidth matrix $\mathbf{H}_{0}(\zeta,\kappa):=\left(h_{0ij}(\zeta,\kappa)\right)_{i,j=1,2}$ with $h_{011}(\zeta,\kappa)=h_{011}$, $h_{022}(\zeta,\kappa)=h_{022}$ and 
  $h_{012}(\zeta,\kappa)=H_{012}[\zeta](\kappa)$;
  \item Compose then the Scott bandwidth matrix $\mathbf{H}(\zeta,\kappa) := H(\zeta)\times \mathbf{H}_{0}(\zeta,\kappa)$.
   \end{enumerate}
  \item Apply LSCV method on the set $\mathcal{S}_1$ of all Scott bandwidth matrices $\mathbf{H}(\zeta,\kappa)$.
\end{enumerate}

\item[A3.] Diagonal bandwidth matrices. 
\begin{enumerate}[1.]
\item Choose two intervals $H_{11}$ and $H_{22}$ related to $h_{11}$ and $h_{22}$, respectively.
\item For $\delta = 1, \ldots, \ell(H_{11})$ and $\gamma = 1, \ldots, \ell(H_{22})$, \\
  Compose the diagonal bandwidth matrix $\mathbf{H}(\delta,\gamma):=\mathbf{Diag}\left(H_{11}(\delta),H_{22}(\gamma)\right)$.
\item Apply LSCV method on the set $\mathcal{D}$ of all diagonal bandwidth matrices $\mathbf{H}(\delta,\gamma)$.
\end{enumerate}
\end{enumerate}
Let us conclude these algorithms by the following precisions. For a given interval $I$, the notation $\ell(I)$ is the total number of subdivisions of $I$ and $I(\eta)$ denotes the real value at the subdivision $\eta$ of $I$. Also, for practical uses of (A1) and (A3), both intervals $H_{11}$ and $H_{22}$ are generally chosen to be $(0,1)$. In the case of the Scott bandwidth matrix (A2), we retain the interval $H=(0,2)$ and the fixed bandwidth matrix $\mathbf{H}_{0} = \widehat{\mathbf{\Sigma}}$, where $\widehat{\mathbf{\Sigma}}$ is the sample covariance matrix. See Figure \ref{Hcvfig} for graphical illustrations.

\subsubsection{Modified version of the estimator}
\label{sssec:Modversionbeta}

Being large in the standard version, the pointwise bias of the beta-Sarmanov kernel estimator (\ref{est}) must be reduced. Following the algorithm of Section~\ref{ssec:Modifiedversion} and without numerical illustration in this paper, the first step  divides  $[0, 1]\times[0, 1]$ in nine subregions of order $\alpha({\mathbf{H}}) = (\alpha_{1}(h_{11}), \alpha_{2}(h_{22}))^{\top}$ with $\alpha_{j}(h_{jj}) > 0$ for $j = 1,2$:
 \begin{enumerate}[a.]
\item only one interior subregion denoted as 
$\mathbb{T}^{\alpha(\mathbf{H}),I}_{2} = (\alpha_{1}(h_{11}), 1 - \alpha_{1}(h_{11})) \times (\alpha_{2}(h_{22}), 1 - \alpha_{2}(h_{22}))$;
 \item eight boundary subregions divided in two parts as
 \begin{enumerate}[(i)]
       \item four angle subregions denoted by
          \begin{eqnarray*}\mathbb{T}_{2}^{\alpha(\mathbf{H}),A}
                         &=& [0, \alpha_{1}(h_{11})] \times [0, \alpha_{2}(h_{22})]  \cup [1 - \alpha_{1}(h_{11}), 1] \times [0, \alpha_{2}(h_{22})] \\
                         && ~\cup [0, \alpha_{1}(h_{11})] \times [1 - \alpha_{2}(h_{22}), 1] \cup [1 - \alpha_{2}(h_{22}), 1] \times [1 - \alpha_{2}(h_{22}), 1],
                         \end{eqnarray*}
        \item four edge subregions denoted by
               \begin{eqnarray*}\mathbb{T}^{\alpha(\mathbf{H}),E}_{2}
                         &=& (\alpha_{1}(h_{11}), 1 - \alpha_{1}(h_{11})) \times [1 - \alpha_{2}(h_{22}), 1]  \cup [0, \alpha_{1}(h_{11})] \times (\alpha_{2}(h_{22}), 1 - \alpha_{2}(h_{22}))   \\
              && ~\cup (\alpha_{1}(h_{11}), 1 - \alpha_{1}(h_{11})) \times [0, \alpha_{2}(h_{22})] \\
              &&~\cup  [1 - \alpha_{1}(h_{11}), 1] \times (\alpha_{2}(h_{22}), 1 - \alpha_{2}(h_{22})).
                         \end{eqnarray*}                          
 \end{enumerate}
\end{enumerate}

As for the second step, we consider the three functions 
$\psi_{1}, \psi_{2}: [0, 1] \rightarrow  [0, 1]$  and  $\psi_{3}:  \mathcal{H} \rightarrow \mathbb{R}$ such that
\begin{equation}
\psi_{j}(z_{j}) = \alpha_{j}(h_{jj})\{z_{j} - \alpha_{j}(h_{jj}) + 1\}, ~~ \forall z_{j} \in [0, 1], \;j = 1, 2\mbox{ and } \psi_{3}(\mathbf{H}) = \frac{h_{12}}{\sqrt{h_{11}h_{22}}}.\label{phi3}
\end{equation}
Each axis of $[0, 1] \times [0, 1]$ has one interior subregion $(\alpha_{j}(h_{jj}), 1 - \alpha_{j}(h_{jj}))$ and two boundary regions $[0, \alpha_{j}(h_{jj})]$ and $[1 - \alpha_{1}(h_{jj}), 1]$ for $j = 1, 2$. Thus, from (\ref{systhetatild}) with $d = 1$, one gets the new parametrization of each margin beta kernel used in $BS_{\theta(\mathbf{x}, \mathbf{H})}$ of (\ref{betakern}) with $\mathbf{x} = (x_{1}, x_{2})^{\top}$: for $j = 1, 2$,
\begin{equation}
    \left\{
      \begin{array}{ll}
       \left(\frac{\psi_{j}(x_{j})}{h_{jj}}, \frac{x_{j}}{h_{jj}}\right)  & \hbox{} \text{ if } x_{j} \in [0, \alpha_{j}(h_{jj})], \\
\left(\frac{x_{j}}{h_{jj}}, \frac{1 - x_{j}}{h_{jj}}\right)  & \hbox{} \text{ if } x_{j} \in (\alpha_{j}(h_{jj}), 1 - \alpha_{j}(h_{jj})),\\
\left(\frac{1 - x_{j}}{h_{jj}}, \frac{\psi_{j}(1 - x_{j})}{h_{jj}}\right)  & \hbox{} \text{ if } x_{j} \in  [1 - \alpha_{i}(h_{jj}), 1].
      \end{array}
\right. \label{param}
\end{equation}
Therefore, using $\psi_{3}$ of (\ref{phi3}) and by combination of (\ref{param}) for each of the nine subregions of $\left[0, 1\right]\times\left[0, 1\right]$, then  $\widetilde{\theta}$ is expressed by
 \begin{equation}
\widetilde{\theta}(\mathbf{x}, \mathbf{H}) = \left\{
      \begin{array}{ll}
       \left(\frac{\psi_{1}(x_{1})}{h_{11}}, \frac{x_{1}}{h_{11}}, \frac{\psi_{2}(x_{2})}{h_{22}}, \frac{x_{2}}{h_{22}}, \frac{h_{12}}{(h_{11}h_{22})^{1/2}}\right)^{\top}  & \hbox{} \text{ if } \mathbf{x} \in [0, \alpha_{1}(h_{11})] \times [0, \alpha_{2}(h_{22})]\\
\left(\frac{\psi_{1}(x_{1})}{h_{11}}, \frac{x_{1}}{h_{11}},\frac{1 - x_{2}}{h_{22}}, \frac{\psi_{2}(1 - x_{2})}{h_{22}}, \frac{h_{12}}{(h_{11}h_{22})^{1/2}}\right)^{\top}  & \hbox{} \text{ if } \mathbf{x}\in [0, \alpha_{1}(h_{11})] \times [1 - \alpha_{2}(h_{22}), 1]\\
     \left(\frac{1 - x_{1}}{h_{11}}, \frac{\psi_{1}(1 - x_{1})}{h_{11}}, \frac{\psi_{2}(x_{2})}{h_{22}}, \frac{x_{2}}{h_{22}}, \frac{h_{12}}{(h_{11}h_{22})^{1/2}}\right)^{\top} & \hbox{} \text{ if } \mathbf{x} \in [1 - \alpha_{1}(h_{11}), 1] \times [0, \alpha_{2}(h_{11})]\\
  \left(\frac{1 - x_{1}}{h_{11}}, \frac{\psi_{1}(1 - x_{1})}{h_{11}}, \frac{1 - x_{2}}{h_{22}}, \frac{\psi_{2}(1 - x_{2})}{h_{22}}, \frac{h_{12}}{(h_{11}h_{22})^{1/2}}\right)^{\top} & \hbox{} \text{ if } \mathbf{x} \in [1 - \alpha_{1}(h_{11}), 1] \times [1 - \alpha_{2}(h_{22}), 1] \\

\left(\frac{x_{1}}{h_{11}}, \frac{1 - x_{1}}{h_{11}}, \frac{x_{2}}{h_{22}}, \frac{1 - x_{2}}{h_{22}}, \frac{h_{12}}{(h_{11}h_{22})^{1/2}}\right)= \widetilde{\theta_I}(\mathbf{x}, \mathbf{H}) & \hbox{} \text{ if } \mathbf{x} \in (\alpha_{1}(h_{11}), 1 - \alpha_{1}(h_{11})) \times (\alpha_{2}(h_{22}), 1 - \alpha_{2}(h_{22}))\\

\left(\frac{\psi_{1}(x_{1})}{h_{11}}, \frac{x_{1}}{h_{11}}, \frac{x_{2}}{h_{22}}, \frac{1 - x_{2}}{h_{22}}, \frac{h_{12}}{(h_{11}h_{22})^{1/2}}\right)^{\top}  & \hbox{} \text{ if } \mathbf{x} \in [0, \alpha_{1}(h_{11})] \times (\alpha_{2}(h_{22}), 1 - \alpha_{2}(h_{22}))\\

        \left(\frac{x_{1}}{h_{11}}, \frac{1 - x_{1}}{h_{11}}, \frac{\psi_{2}(x_{2})}{h_{22}}, \frac{1 - x_{2}}{h_{22}}, \frac{h_{12}}{(h_{11}h_{22})^{1/2}}\right)^{\top}   & \hbox{} \text{ if } \mathbf{x} \in (\alpha_{1}(h_{11}), 1 - \alpha_{1}(h_{11})) \times [0, \alpha_{2}(h_{22})]\\

\left(\frac{x_{1}}{h_{11}}, \frac{1 - x_{1}}{h_{11}}, \frac{1 - x_{2}}{h_{22}}, \frac{\psi_{2}(1 - x_{2})}{h_{22}}, \frac{h_{12}}{(h_{11}h_{22})^{1/2}}\right)^{\top} & \hbox{} \text{ if } \mathbf{x} \in (\alpha_{1}(h_{11}), 1 - \alpha_{1}(h_{11})) \times [1 - \alpha_{2}(h_{22}), 1]\\

  \left(\frac{1 - x_{1}}{h_{11}}, \frac{\psi_{1}(1 - x_{1})}{h_{11}}, \frac{x_{2}}{h_{22}}, \frac{1 - x_{2}}{h_{22}}, \frac{h_{12}}{(h_{11}h_{22})^{1/2}}\right)^{\top} & \hbox{} \text{ if } \mathbf{x} \in [1 - \alpha_{1}(h_{11}), 1]\times (\alpha_{2}(h_{22}), 1 - \alpha_{2}(h_{22})).

      \end{array}
    \right. \label{tilde}
\end{equation} 
Notice that the last component $h_{12}(h_{11}h_{22})^{-1/2}$ of $\widetilde{\theta}$  does not change in all subregions of the support. This is because of the choice of $\mathbf{m}_{0}$ related to the construction of the beta-Sarmanov kernel (\ref{betakern}). According to Proposition~\ref{Proprs} the modified beta-Sarmanov kernel obtained from (\ref{tilde}) and denoted by $BS_{\widetilde{\theta}(\mathbf{x}, \mathbf{H})}$  is  an associated kernel with  $\mathbb{S}_{\widetilde{\theta}(\mathbf{x}, \mathbf{H})} = \left[0, 1\right]\times\left[0, 1\right]$,
\begin{equation}
\mathbf{a}_{\widetilde{\theta}}(\mathbf{x}, \mathbf{H})
= \left\{
      \begin{array}{ll}

       \left(\frac{(1 - x_{1})\psi_{1}(x_{1}) + x_{1}^{2}}{x_{1} + \psi_{1}(x_{1})}, \frac{(1 - x_{2})\psi_{2}(x_{2}) + x_{2}^{2}}{x_{2} + \psi_{2}(x_{2})}\right)^{\top}  & \hbox{} \small{\text{ if } \mathbf{x} \in [0, \alpha_{1}(h_{11})] \times [0, \alpha_{2}(h_{22})]}\\

\left(\frac{(1 - x_{1})\psi_{1}(x_{1}) + x_{1}^{2}}{x_{1} + \psi_{1}(x_{1})}, 0\right)^{\top}  & \hbox{} \small{\text{ if } \mathbf{x} \in [0, \alpha_{1}(h_{11})] \times (\alpha_{2}(h_{22}), 1 - \alpha_{2}(h_{22}))}\\

\left(\frac{(1 - x_{1})\psi_{1}(x_{1}) + x_{1}^{2}}{x_{1} + \psi_{1}(x_{1})}, \frac{(1 - x_{2})\{x_ {2} - \psi_{2}(1 - x_{2})\}}{1 - x_{2} + \psi_{2}(1 - x_{2})}\right)^{\top}  & \hbox{} \small{\text{ if } \mathbf{x} \in [0, \alpha_{1}(h_{11})] \times [1 - \alpha_{2}(h_{22}), 1]}\\

        \left(0, \frac{(1 - x_{2})\psi_{2}(x_{2}) + x_{2}^{2}}{x_{2} + \psi_{2}(x_{2})}\right)^{\top} & \hbox{} \small{\text{ if } \mathbf{x} \in (\alpha_{1}(h_{11}), 1 - \alpha_{1}(h_{11})) \times [0, \alpha_{2}(h_{22})]}\\

\left(0, 0\right)^{\top} & \hbox{} \small{\text{ if } \mathbf{x} \in (\alpha_{1}(h_{11}), 1 - \alpha_{1}(h_{11})) \times (\alpha_{2}(h_{22}), 1 - \alpha_{2}(h_{22}))}\\

\left(0, \frac{(1 - x_{2})\{x_ {2} - \psi_{2}(1 - x_{2})\}}{1 - x_{2} + \psi_{2}(1 - x_{2})}\right)^{\top} & \hbox{} \small{ \text{ if } \mathbf{x} \in (\alpha_{1}(h_{11}), 1 - \alpha_{1}(h_{11})) \times [1 - \alpha_{2}(h_{22}), 1]}\\

        \left(\frac{(1 - x_{1})\{x_ {1} - \psi_{1}(1 - x_{1})\}}{1 - x_{1} + \psi_{1}(1 - x_{1})}, \frac{(1 - x_{2})\psi_{2}(x_{2}) + x_{2}^{2}}{x_{2} + \psi_{2}(x_{2})}\right)^{\top} & \hbox{} \small{ \text{ if } \mathbf{x} \in [1 - \alpha_{1}(h_{11}), 1] \times [0, \alpha_{2}(h_{11})]}\\

  \left(\frac{(1 - x_{1})\{x_ {1} - \psi_{1}(1 - x_{1})\}}{1 - x_{1} + \psi_{1}(1 - x_{1})}, 0\right)^{\top} & \hbox{} \small{\text{ if } \mathbf{x} \in [1 - \alpha_{1}(h_{11}), 1]\times (\alpha_{2}(h_{22}), 1 - \alpha_{2}(h_{22}))}\\

  \left(\frac{(1 - x_{1})\{x_ {1} - \psi_{1}(1 - x_{1})\}}{1 - x_{1} + \psi_{1}(1 - x_{1})}, \frac{(1 - x_{2})\{x_ {2} - \psi_{2}(1 - x_{2})\}}{1 - x_{2} + \psi_{2}(1 - x_{2})}\right)^{\top} & \hbox{} \small{\text{ if } \mathbf{x} \in [1 - \alpha_{1}(h_{11}), 1] \times [1 - \alpha_{2}(h_{22}), 1]}

      \end{array}
    \right.
\end{equation}
and, $\mathbf{B}_{\widetilde{\theta}}(\mathbf{x},  \mathbf{H}) = (b_{\widetilde{\theta} ij})_{i, j = 1, 2}$ such that $ b_{\widetilde{\theta} 12} = \left( h_{12}/(h_{11}h_{22})^{1/2}\right)b_{\widetilde{\theta}11}b_{\widetilde{\theta} 22}$ and 
$\left(b_{\widetilde{\theta}11},b_{\widetilde{\theta} 22}\right)$ is detailed as
\begin{eqnarray*}
&&\left\{
      \begin{array}{ll}

\left(\frac{h_{11}x_{1}\psi_{1}(x_{1})\{x_{1} + \psi_{1}(x_{1})\}^{-2}}{\{x_{1} + \psi_1(x_{1}) + h_{11}\}} , \frac{h_{22}x_{2}\psi_{2}(x_{2})\{x_{2} + \psi_{2}(x_{2})\}^{-2}}{\{x_{2} + \psi_2(x_{2}) + h_{22}\}} \right)  & \hbox{} \small{\text{if } \mathbf{x} \in [0, \alpha_{1}(h_{11})] \times [0, \alpha_{2}(h_{22})]}\\

\left(\frac{h_{11}x_{1}\psi_{1}(x_{1})\{x_{1} + \psi_{1}(x_{1})\}^{-2}}{\{x_{1} + \psi_1(x_{1}) + h_{11}\}} , \frac{h_ {22}x_{2}(1 - x_{2})}{1 + h_{22}}\right)  & \hbox{} \small{\text{if } \mathbf{x} \in [0, \alpha_{1}(h_{11})] \times (\alpha_{2}(h_{22}), 1 - \alpha_{2}(h_{22}))}\\

\left(\frac{h_{11}x_{1}\psi_{1}(x_{1})\{x_{1} + \psi_{1}(x_{1})\}^{-2}}{\{x_{1} + \psi_1(x_{1}) + h_{11}\}} , \frac{h_{22}(1 - x_{2})\{1 - x_{2} + \psi_{2}(1 - x_{2})\}^{-2}}{\{\psi_{2}(1 - x_{2})\}^{-1}\{1 - x_{2} + \psi_{2}(1 - x_{2}) + h_{22}\}}\right)  & \hbox{} \small{ \text{if } \mathbf{x} \in [0, \alpha_{1}(h_{11})] \times [1 - \alpha_{2}(h_{22}), 1]}\\

\left(\frac{h_ {11}x_{1}(1 - x_{1})}{1 + h_{11}}, \frac{h_{22}x_{2}\psi_{2}(x_{2})\{x_{2} + \psi_{2}(x_{2})\}^{-2}}{\{x_{2} + \psi_2(x_{2}) + h_{22}\}} \right) & \hbox{} \small{\text{if } \mathbf{x} \in (\alpha_{1}(h_{11}), 1 - \alpha_{1}(h_{11})) \times [0, \alpha_{2}(h_{22})]}\\

\left(\frac{h_ {11}x_{1}(1 - x_{i})}{1 + h_{11}}, \frac{h_ {22}x_{2}(1 - x_{2})}{1 + h_{22}}\right)& \hbox{} \small{\text{if } \mathbf{x} \in \mathbb{T}^{\alpha(\mathbf{H}),I}_{2}}\\

\left(\frac{h_ {11}x_{1}(1 - x_{1})}{1 + h_{11}}, \frac{h_{22}(1 - x_{2})\{1 - x_ {2} + \psi_{2}(1 - x_{2})\}^{-2}}{\{\psi_{2}(1 - x_{2})\}^{-1}\{1 - x_{2} + \psi_{2}(1 - x_{2}) + h_{22}\}}\right) & \hbox{} \small{\text{ if } \mathbf{x} \in (\alpha_{1}(h_{11}), 1 - \alpha_{1}(h_{11})) \times [1 - \alpha_{2}(h_{22}), 1]}\\

\left(\frac{h_{11}(1 - x_{1})\{1 - x_{1} + \psi_{1}(1 - x_{1})\}^{-2}}{\{\psi_{1}(1 - x_{1})\}^{-1}\{1 - x_{1} + \psi_{1}(1 - x_{1}) + h_{11}\}}, \frac{h_{22}x_{2}\psi_{2}(x_{2})\{x_{2} + \psi_{2}(x_{2})\}^{-2}}{\{x_{2} + \psi_2(x_{2}) + h_{22}\}} \right) & \hbox{} \small{ \text{if } \mathbf{x} \in [1 - \alpha_{1}(h_{11}), 1] \times [0, \alpha_{2}(h_{11})]}\\

\left(\frac{h_{11}(1 - x_{1})\{1 - x_ {1} + \psi_{1}(1 - x_{1})\}^{-2}}{\{\psi_{1}(1 - x_{1})\}^{-1}\{1 - x_{1} + \psi_{1}(1 - x_{1}) + h_{11}\}}, \frac{h_ {22}x_{2}(1 - x_{2})}{1 + h_{22}}\right) & \hbox{} \small{\text{ if } \mathbf{x} \in [1 - \alpha_{1}(h_{11}), 1]\times (\alpha_{2}(h_{22}), 1 - \alpha_{2}(h_{22}))}\\

\left(\frac{h_{11}(1 - x_{1})\{1 - x_{1} + \psi_{1}(1 - x_{1})\}^{-2}}{\{\psi_{1}(1 - x_{1})\}^{-1}\{1 - x_{1} + \psi_{1}(1 - x_{1}) + h_{11}\}}, \frac{h_{22}(1 - x_{2})\{1 - x_ {2} + \psi_{2}(1 - x_{2})\}^{-2}}{\{\psi_{2}(1 - x_{2})\}^{-1}\{1 - x_{2} + \psi_{2}(1 - x_{2}) + h_{22}\}}\right) & \hbox{} \small{\text{if } \mathbf{x} \in [1 - \alpha_{1}(h_{11}), 1] \times [1 - \alpha_{2}(h_{22}), 1]}.
      \end{array}
    \right. \label{Bii}
\end{eqnarray*}

The corresponding modified beta-Sarmanov kernel estimator
\begin{equation}
  \widetilde{f}_{n}(\mathbf{x}) = \dfrac{1}{n} \displaystyle \sum_{i=1}^{n}BS_{\widetilde{\theta}(\mathbf{x}, \mathbf{H})}(X_{i}),~~ \forall  \mathbf{x} \in \left[0, 1\right]\times\left[0, 1\right] \label{estmod}
\end{equation}
has, for all $\mathbf{x} \in  \mathbb{T}^{\alpha(\mathbf{H}),I}_{2} = (\alpha_{1}(h_{11}), 1 - \alpha_{1}(h_{11})) \times (\alpha_{2}(h_{22}), 1 - \alpha_{2}(h_{22}))$,
\begin{eqnarray}
 \Bias\left\{\widetilde{f}_{n}(\mathbf{x})\right\}  &=& \frac{1}{2} \left\{b_{\widetilde{\theta}_{I}11}\frac{\partial^{2} f}{\partial x_{1}^2}(\mathbf{x})  + 2b_{\widetilde{\theta}_{I}12}\frac{\partial^{2} f}{\partial x_{2}\partial x_{2}}(\mathbf{x})  + b_{\widetilde{\theta}_{I}22}\frac{\partial^{2} f}{\partial x_{2}^2}(\mathbf{x}) \right\} + \oldstylenums{0}(h_{11}^{2} + 2h_{12}^{2} + h_{22}^{2}) \nonumber
\end{eqnarray}
and
\begin{equation*}
  \Var\left\{\widetilde{f}_{n}(\mathbf{x})\right\}  \simeq \Var\left\{\widehat{f}_{n}(\mathbf{x})\right\}  ~as~ n\to\infty.
\end{equation*}Thus, the asymptotic expression of the MISE of $\widetilde{f}_{n}$ on $\mathbb{T}^{\alpha(\mathbf{H}),I}_{2}$ is given by
\begin{eqnarray*}
\AMISE_{\widetilde{\theta}_I}\left(\widetilde{f}_{n}\right)  &=&   \int_{\mathbb{T}^{\alpha(\mathbf{H}),I}_{2}}\left(\left[ \frac{1}{2} \left\{b_{\widetilde{\theta}_{I}11}\frac{\partial^{2} f}{\partial x_{1}^2}(\mathbf{x})  + 2b_{\widetilde{\theta}_{I}12}\frac{\partial^{2} f}{\partial x_{2}\partial x_{2}}(\mathbf{x})  + b_{\widetilde{\theta}_{I}22}\frac{\partial^{2} f}{\partial x_{2}^2}(\mathbf{x}) \right\} \right]^{2}  \right. \nonumber \\
&&\left.+ \frac{1}{n} \left\|BS_{\theta(\mathbf{x}, \mathbf{H})}\right\|^{2}_{2}f(\mathbf{x})\right)d\bold{x}. 
\end{eqnarray*}
The modified beta-Sarmanov kernel also depend on the choice of scalars $\alpha_{j}(h_{jj})$ for $j = 1, 2$. The user can set the values of  $\alpha_{j}(h_{jj})$ according to his practical objective. For example in univariate case, \cite{{C99},{C00}} took $\alpha_{j}(h_{jj}) = 2h_{jj}$.
From (\ref{theta}) to (\ref{estmod}) and when $h_{12} = 0$, we have the same formulas for multiple associated kernels of  \cite{BR10}.  Also, similar results can be obtained for the Scott bandwidth matrices. However, numerical illustrations are so long and tedious tasks; see, e.g., \cite{HiruSakudo13} for $d=1$.
 
\section{Simulation studies and real data analysis}
\label{sec:Simulations}

In this section, we compare the performance of the three forms of bandwidth matrices of  Table~\ref{parameterstable}. The optimal bandwidth matrix is chosen by LSCV method (\ref{Hcv}) using the algorithms A1, A2 and A3 given at the end of Section~\ref{sssec:standbeta} and their indications. All computations were done on the computational resource\footnote{Dell Poweredge R900, Processor Xeon X7350, 2.93 GHz, 32 Go RAM} of Laboratoire de Math\'ematiques de Besan\c{c}on by using the R software; see \citet{R13}. The comparisons will be done using the standard version of the beta-Sarmanov kernel estimator (\ref{est}) through simulations studies and an illustration on real dataset.

\subsection{Simulation studies}
\label{ssec:Simulations studies}

\begin{figure}[htbp!]

 \begin{center}
  \mbox{
\subfloat[(A)]{\includegraphics[width=210pt,height=155pt,scale=0.95]{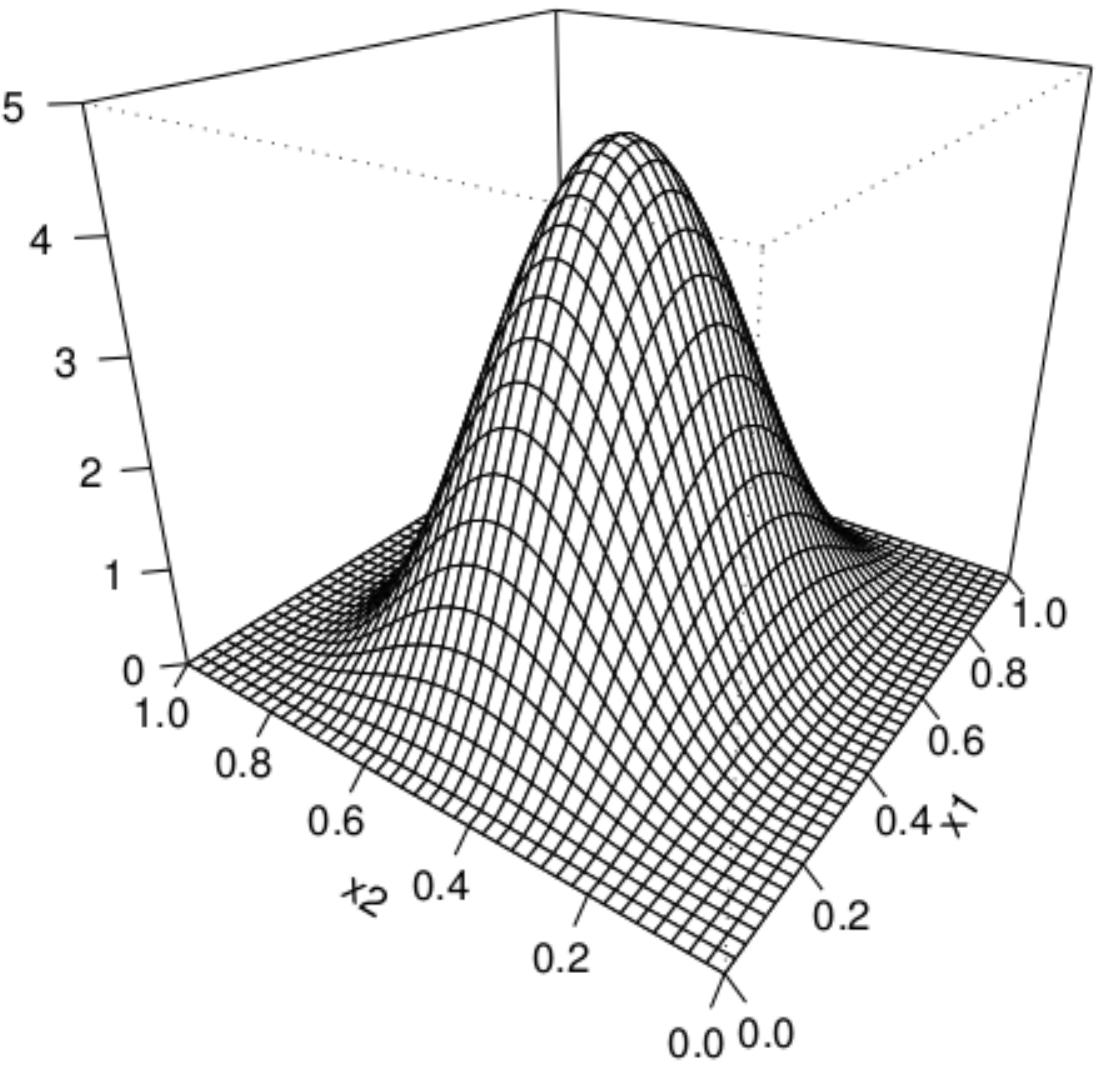} } \:
\subfloat[(B)]{\includegraphics[width=210pt,height=155pt,scale=0.95]{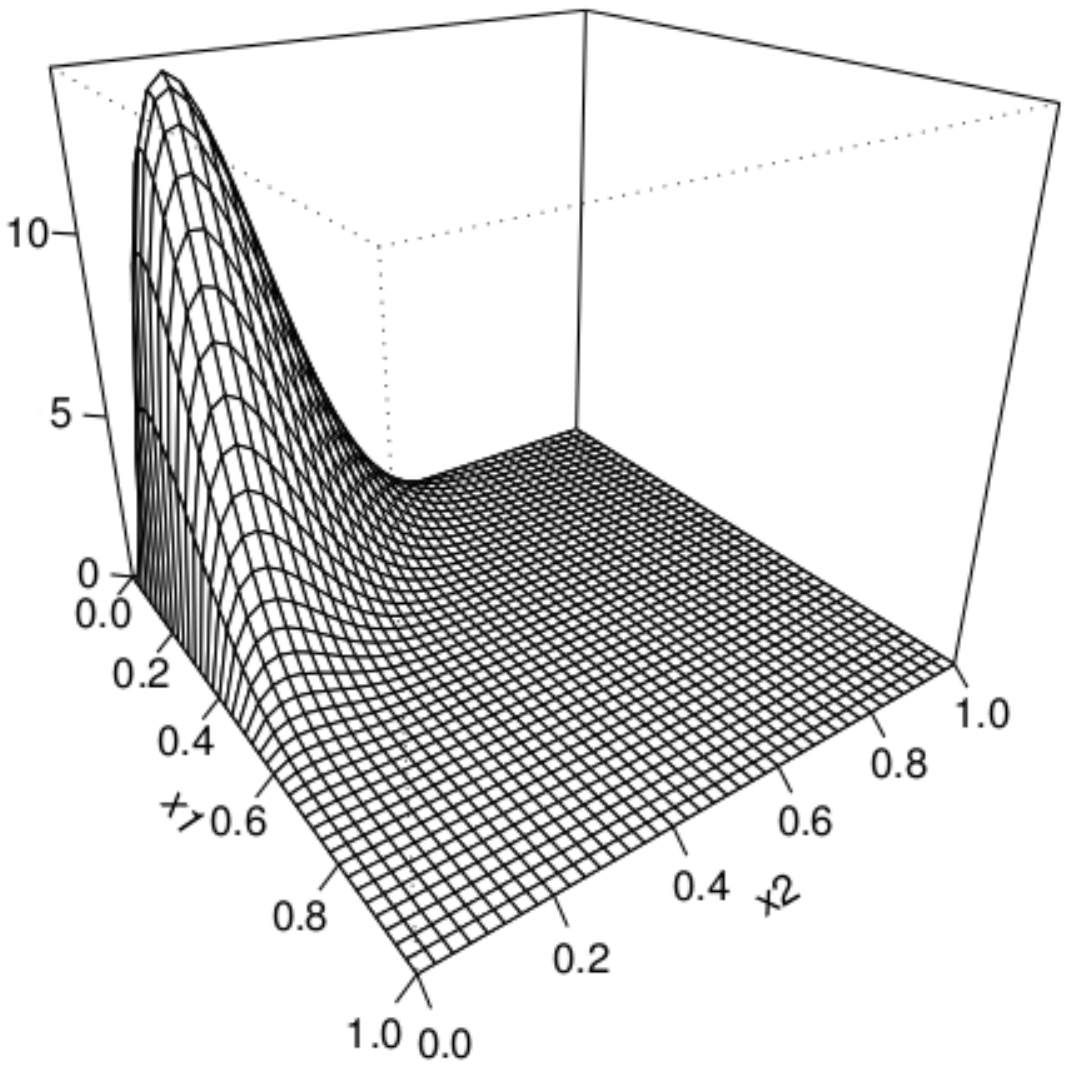} }  
 } 
   \mbox{
\subfloat[(C)]{\includegraphics[width=210pt,height=155pt,scale=0.95]{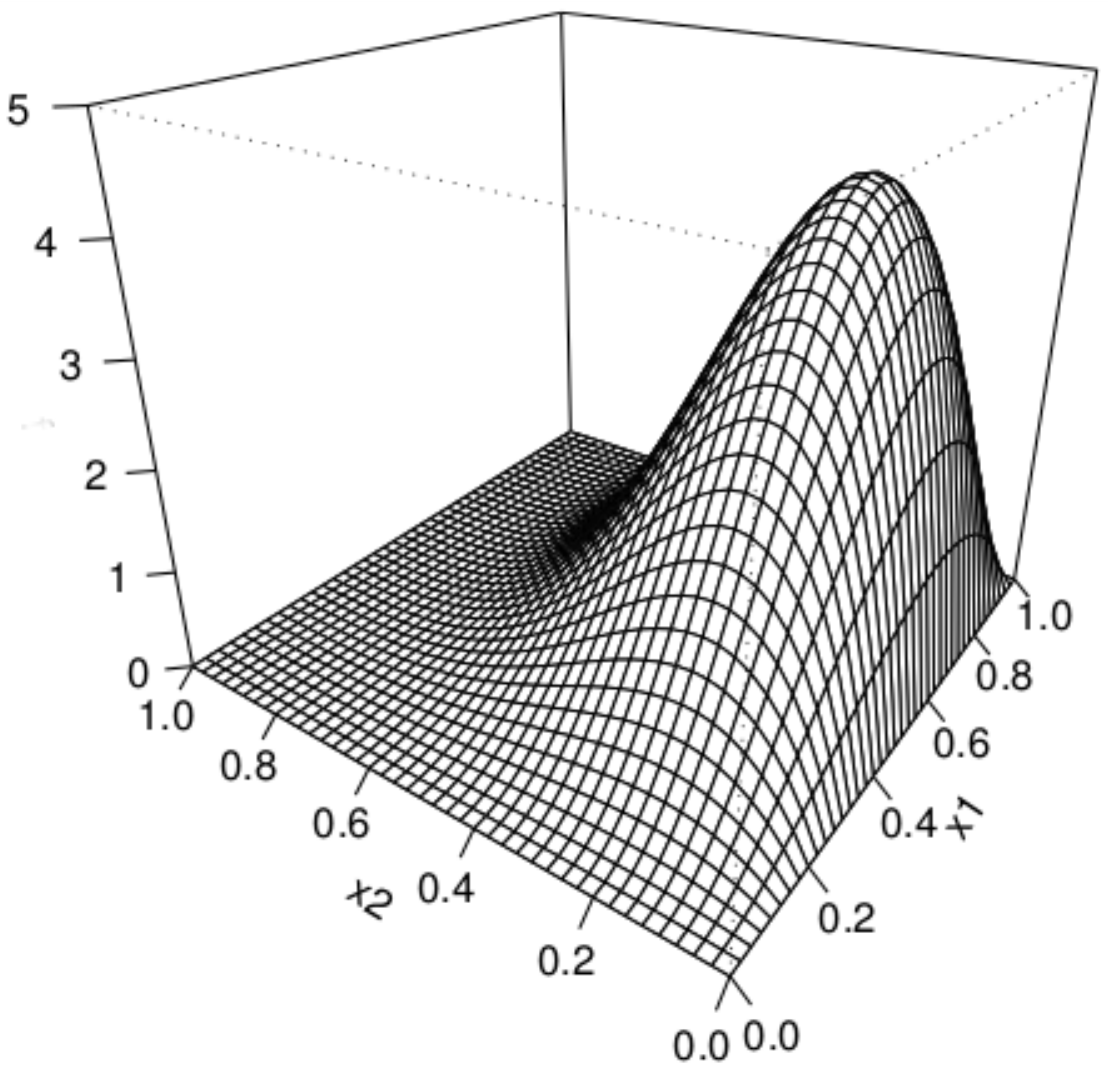} } \:
\subfloat[(D)]{\includegraphics[width=210pt,height=155pt,scale=0.95]{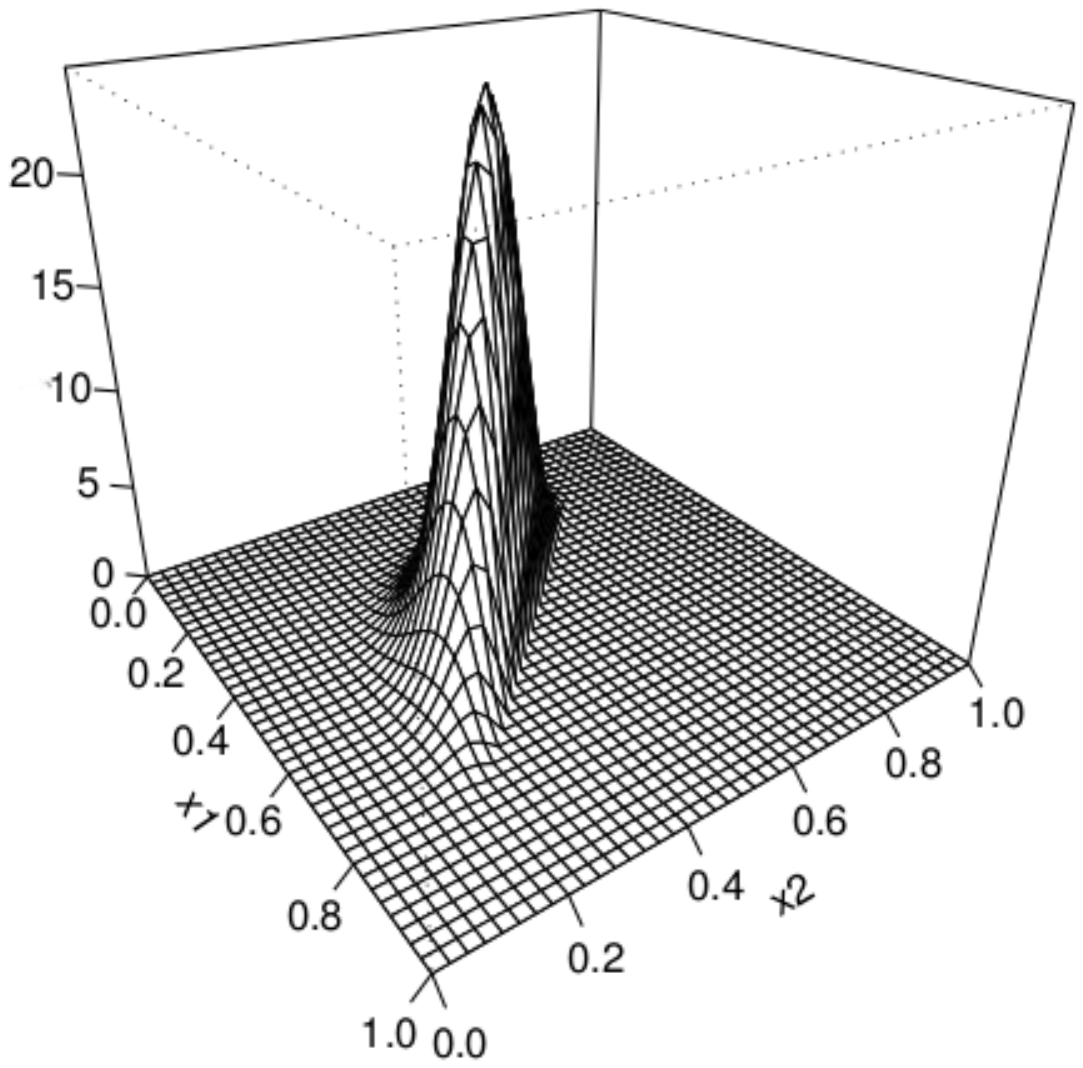} }
}

\mbox{
\subfloat[(E)]{\includegraphics[width=210pt,height=155pt,scale=0.95]{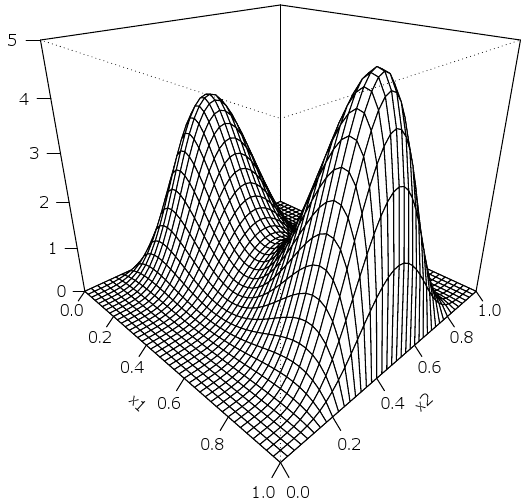} } \:
\subfloat[(F)]{\includegraphics[width=210pt,height=155pt,scale=0.95]{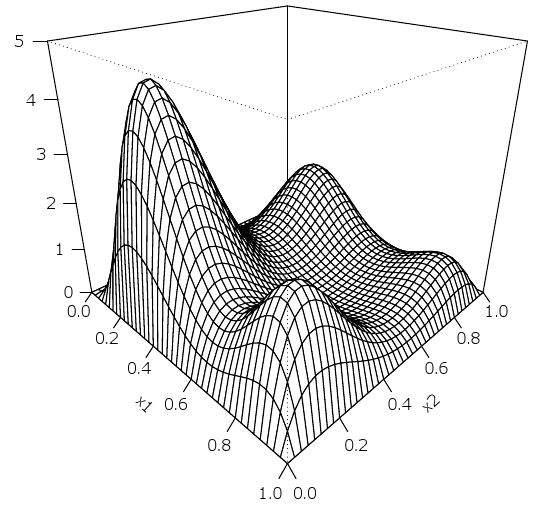} }
}
 \end{center}
\caption{Six plots of density functions defined at the begining of Section~\ref{ssec:Simulations studies} and considered for simulations.}\label{simdist}
\end{figure}

We consider six target densities with supports included in $[0,1]\times[0,1]$ and labeled A, B, C, D, E and F respectively. They have different correlation structure and some local modes. The plots for these densities are given in Figure~\ref{simdist}.
\begin{itemize}
 \item Density A is the bivariate beta density without correlation $\rho = 0$ such that $(p_{1}, q_{1}) = (3,3)$ and $(p_{2},q_{2}) = (5,5)$ as parameters values in univariate beta density (\ref{gunif}), respectively;
\item  density B is the bivariate Dirichlet density
\begin{equation*}
 f(v_{1}, v_{2}) = \frac{\Gamma(\alpha_{1} + \alpha_{2} + \alpha_{3})}{\Gamma(\alpha_{1})\Gamma(\alpha_{2})\Gamma( \alpha_{3})}v_{1}^{\alpha_{1} - 1}v_{2}^{\alpha_{2} - 1}(1 - v_{1} - v_{2})^{\alpha_{3} - 1}, \mbox{ } ~~~\mbox{ } v_{1}, v_{2} \geq 0, v_{1} + v_{2} \leq 1,
\end{equation*}
where $\Gamma(\cdot)$ is the classical gamma function, with parameters values $\alpha_{1} = \alpha_{2} = 2$, $\alpha_{3} =7$ and, therefore, the moderate value of $\rho=-(\alpha_{1}\alpha_{2})^{1/2}(\alpha_{1}+\alpha_{3})^{-1/2}(\alpha_{2}+\alpha_{3})^{-1/2}= -0.2222$; 
\item  density C is the bivariate beta density without correlation $\rho = 0$ such that $(p_{1},q_{1})=(3, 2)$ and $(p_{2}, q_{2}) = (2, 5)$ in (\ref{gunif}), respectively;
\item density D is the bivariate Dirichlet as the density B but with $\alpha_{1} = \alpha_{2} = 10$,  $\alpha_{3} = 3$ and then $\rho = -0.7692$.
\item density E is the bivariate density without correlation $\rho = 0$ defined as follows:
\begin{equation*}
 f(v_{1}, v_{2}) = \left[(3/7)g_{1}(v_{1}) + (4/7)g_{2}(v_{1})\right]\times g_{3}(v_{2})\, 
\end{equation*}such that $g_{1}$, $g_{2}$ and $g_{3}$ are univariate beta densities 
(\ref{gunif}) with parameters values $(p_{1}, q_{1}) = (2,7)$ and $(p_{2},q_{2}) = (7,2)$ and $(p_{3}, q_{3}) = (6,6)$.

\item density F is the bivariate density without correlation $\rho = 0$:
\begin{equation*}
 f(v_{1}, v_{2}) = \left[(8/11)g_{1}(v_{1}) + (3/11)g_{2}(v_{1})\right]\times \left[ (5/7)g_{3}(v_{2}) + (2/7)g_{4}(v_{2})\right]\, 
\end{equation*}such that $g_{1}$, $g_{2}$, $g_{3}$ and $g_{4}$ are univariate beta densities (\ref{gunif}) with parameters values $(p_{1}, q_{1}) = (3.5,7)$ and $(p_{2},q_{2}) = (7,3.5)$, $(p_{3}, q_{3}) = (7,2)$ and $(p_{4}, q_{4}) = (2,7)$.
\end{itemize}

Table~\ref{Timehcv} presents the execution times needed for computing the LSCV method for each of the three types of bandwidth matrix with respect to only one replication of sample sizes $n = 100$ and $500$ for each of the target densities A, B, C, D, E and F of Figure~\ref{simdist}. For $n=100$ the computational times of the LSCV method for full bandwidth matrix are longer than both the Scott and diagonal bandwidth matrices, which have the almost identical Central Processing Unit (CPU) times. The constraints (\ref{h12}) of the correlation of Sarmanov induce that execution times of the Scott bandwidth matrix are relatively a bit longer than for the diagonal ones. Let us note that for full bandwidth matrices, the execution times become very large when the number of observations is large as seen in Table~\ref{Timehcv}; however, these CPU times  can be considerably reduced by parallelism processing, in particular for for the full LSCV method. These constraints (\ref{h12}) reflect the difficulty of finding the appropriate bandwidth matrix with correlation structure by LSCV method (\ref{Hcv}). Hence, we are able to infer that the Scott and diagonal LSCV method do not impose excessive computation burdens; and, the Scott procedure takes into account the structure of (null and moderate) correlations in the sample. 

\begin{table}[!htbp]
\begin{center}
\begin{tabular}{rlrrrrrr}
\hline \hline
 \multicolumn{1}{r}{$n$} &$\mathbf{H}$ & \multicolumn{1}{c}{A}& \multicolumn{1}{c}{B}& \multicolumn{1}{c}{C}&\multicolumn{1}{c}{D} &\multicolumn{1}{c}{E} &\multicolumn{1}{c}{F}  \\
\hline
\multirow{3}{*}{$100$}  & Full   & 2.8310 &2.7570&  2.8793 & 2.7931&27940&2.8362  \\
            &Scott &0.2173 & 0.2251& 0.2207& 0.2202&0.2212&0.2204 \\
            & Diagonal &0.1436 & 0.1456& 0.1526& 0.1536&1.1531&1.1446  \\            
 \hline \hline
\end{tabular}
\caption{Typical  CPU times (in hours) for one replication  of LSCV method (\ref{Hcv}) by using the algorithms A1, A2 and A3 given at the end of Section~\ref{sssec:standbeta}.} \label{Timehcv}
\end{center}
\end{table}

We now examine the efficiency of various bandwidth matrices in Table~\ref{parameterstable} via $$\widehat{ISE}=\dfrac{1}{N}\displaystyle{\sum}_{m=1}^{N}\displaystyle{\int}_{\left[0,1\right]\times\left[0,1\right]}\left\lbrace\widehat{f}_{n}(\mathbf{x})- f(\mathbf{x})\right\rbrace^{2} d\bold{x},$$ where $N =100$ is the  number of replications. Table~\ref{Err} shows some expected values of $\widehat{ISE}$ for the three forms of bandwidth matrices with respect  to the densities A, B, C, D, E and F of Figure~\ref{simdist}, and according only to the sample size $n = 100$ because of excess of computational times (see Table~\ref{Timehcv}). 
Globally, the full and Scott bandwidth matrices with correlation structure perform better in terms of the quality of smoothing than the diagonal one without correlation. Even if the correlation is almost non-existent in the sample (e.g. models  A and C of Table~\ref{Err}), we attend the good behavior of the full and Scott bandwidth matrices. Also, these bandwidth matrices with correlation structure suit for   multimodal target densities  (e.g. E and F of Table~\ref{Err}). For moderate correlation (e.g. models B of Table~\ref{Err}) we can recommend the light version of bandwidth matrices with correlation structure which is the Scott one. As for strong correlation in the sample (e.g. models D of Table~\ref{Err}) it is not preferable to use the Scott bandwidth matrix.
\begin{table}[!htbp]
\begin{center}
\begin{tabular}{clll}
\hline
\hline
Models &\multicolumn{1}{c}{Full}  &\multicolumn{1}{c}{Scott} &\multicolumn{1}{c}{Diagonal}  \\\hline

    A& 0.2554(0.2253) & 0.1887(0.0890) &0.2963(0.2121)  \\

   B & 0.8571(0.2422) & 0.5999(0.2582) &0.9743(0.5973) \\

  C& 0.1985(0.0737) & 0.2150(0.0828) &0.2384(0.1740) \\

    D  & 1.1481(0.0937) & 9.1188(0.7466) &1.6420(0.5735)  \\

    E  & 0.2786(0.0785) & 0.2498(0.0985) &0.5056(0.1526)  \\

    F  &0.6025(0.1122) &0.6791(0.1475)  & 0.8025(0.1254) \\
 \hline\hline
\end{tabular}
\caption{Expected values (and their standard deviation) of $\widehat{ISE}$ with $N=100$ replications of sample sizes $n = 100$ using three types of bandwidth matrices $\widehat{\mathbf{H}}$ (\ref{Hcv}) for each of four models of Figure~\ref{simdist}.} \label{Err}
\end{center}
\end{table}

Finally, Tables~\ref{Timehcv} and \ref{Err} indicate that the choice of the Scott bandwidth matrix using LSCV method is a good alternative to the full one in the purpose of preserving a correlation structure in the bandwidth matrix for an associated kernel estimator.
 
\subsection{Real data analysis}
\label{ssec:Illustration}

\begin{figure}[htbp!]
\begin{center}
  \mbox{
\subfloat[(o) ]{\includegraphics[width=195pt,height=195pt,scale=0.95]{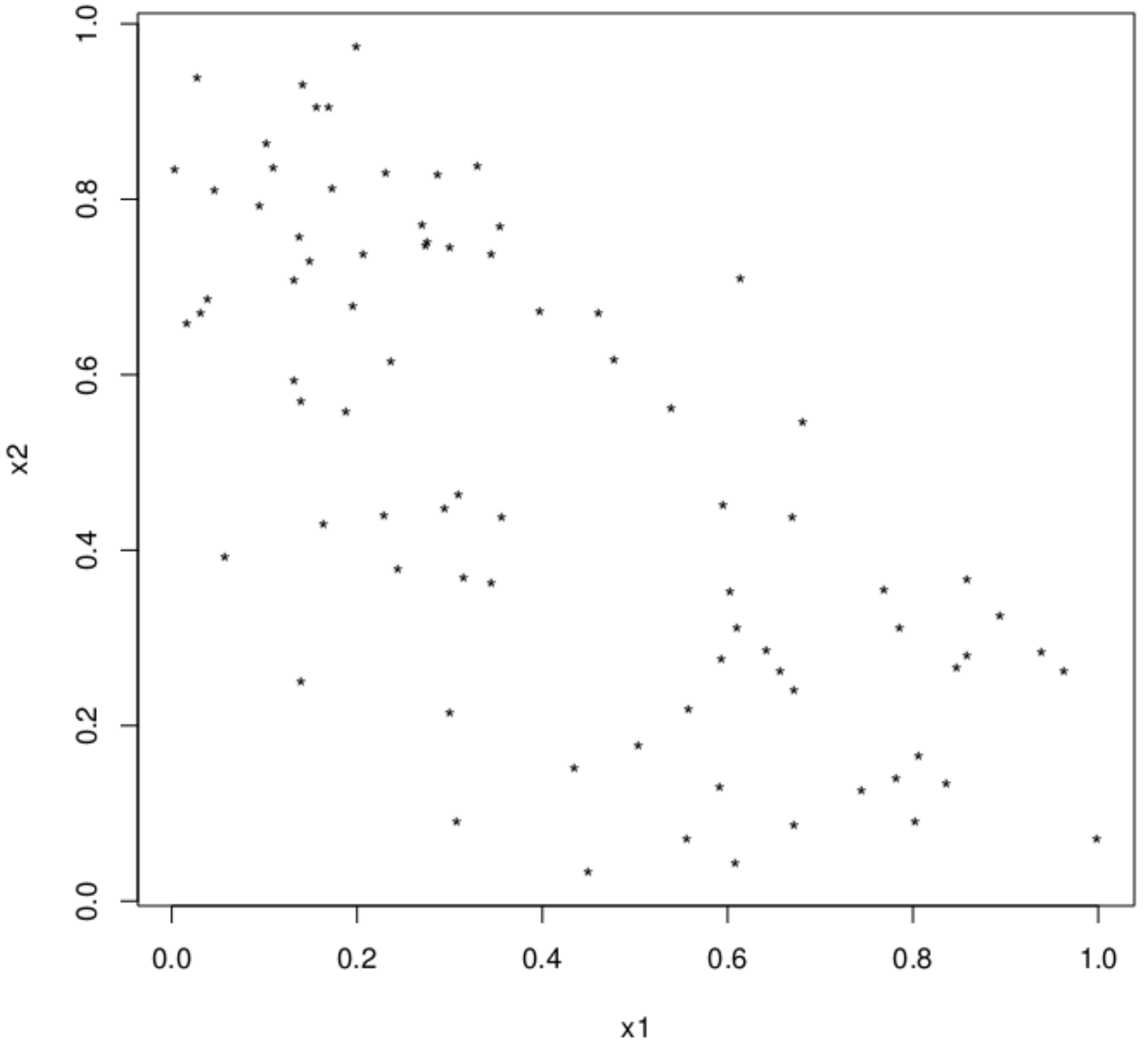} } \:
\subfloat[(a)]{\includegraphics[width=215pt,height=200pt,scale=0.95]{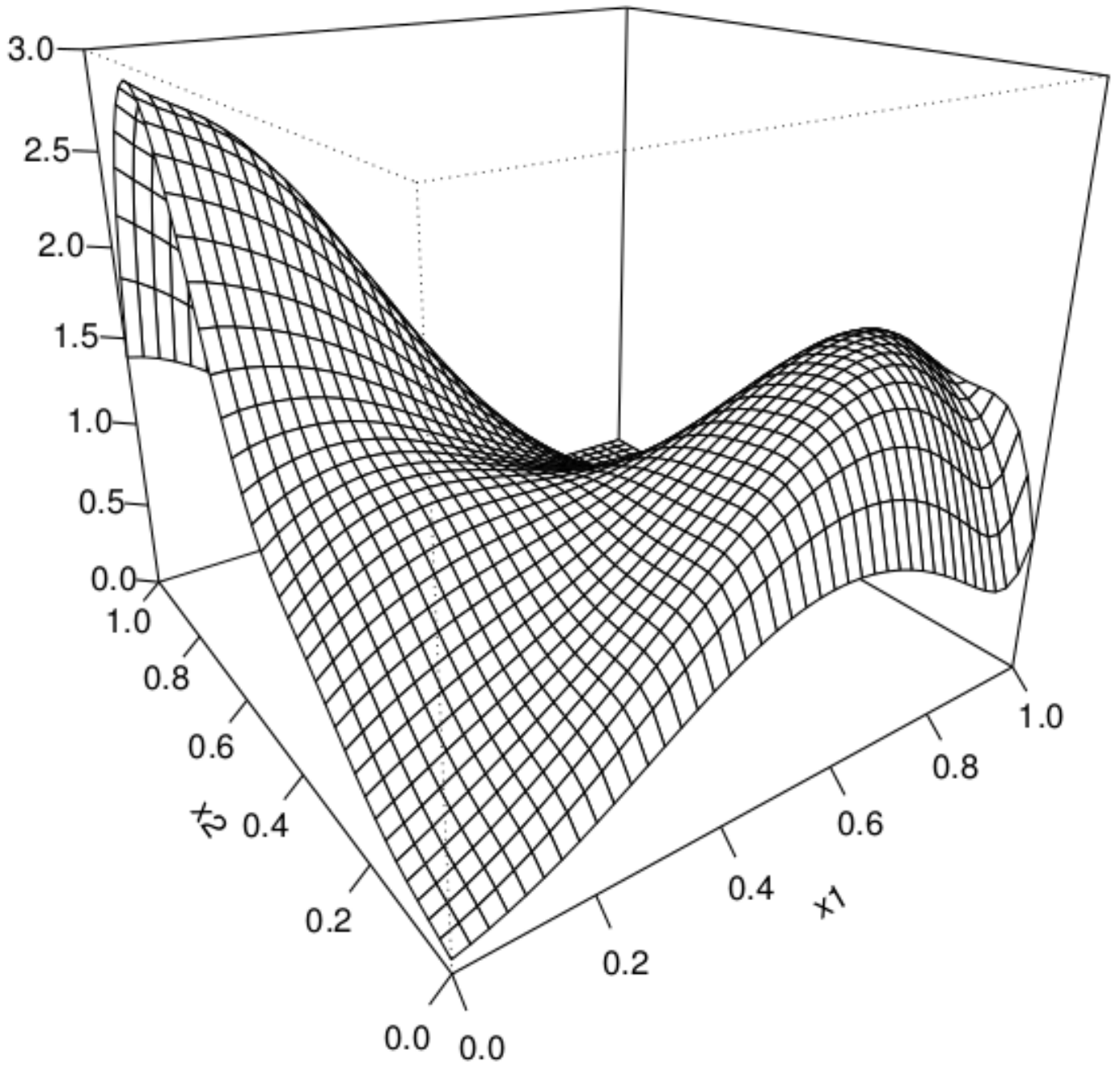} }  } 
  \mbox{
\subfloat[(b) ]{\includegraphics[width=215pt,height=200pt,scale=0.95]{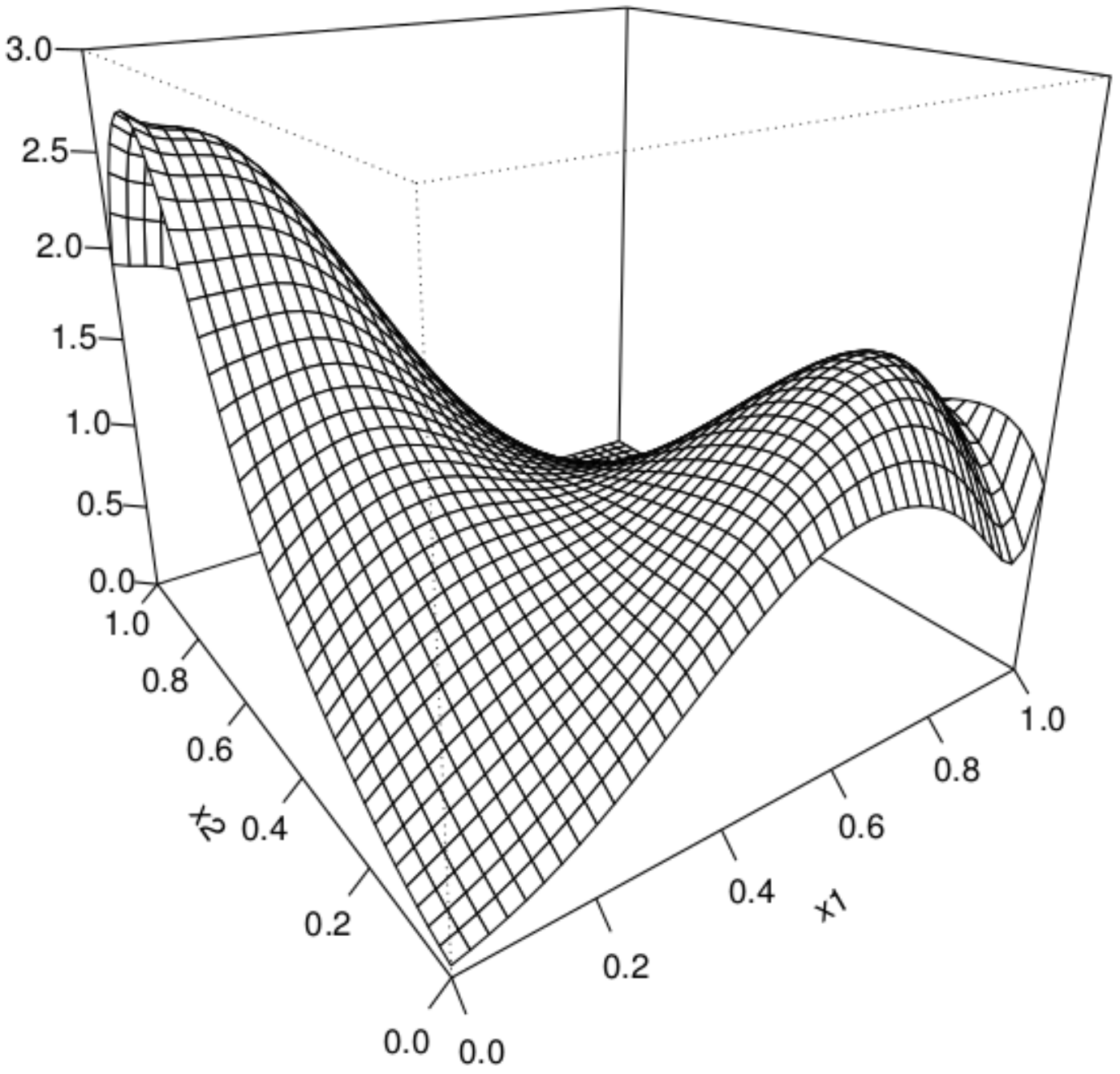} } \:
\subfloat[(c)]{\includegraphics[width=215pt,height=200pt,scale=0.95]{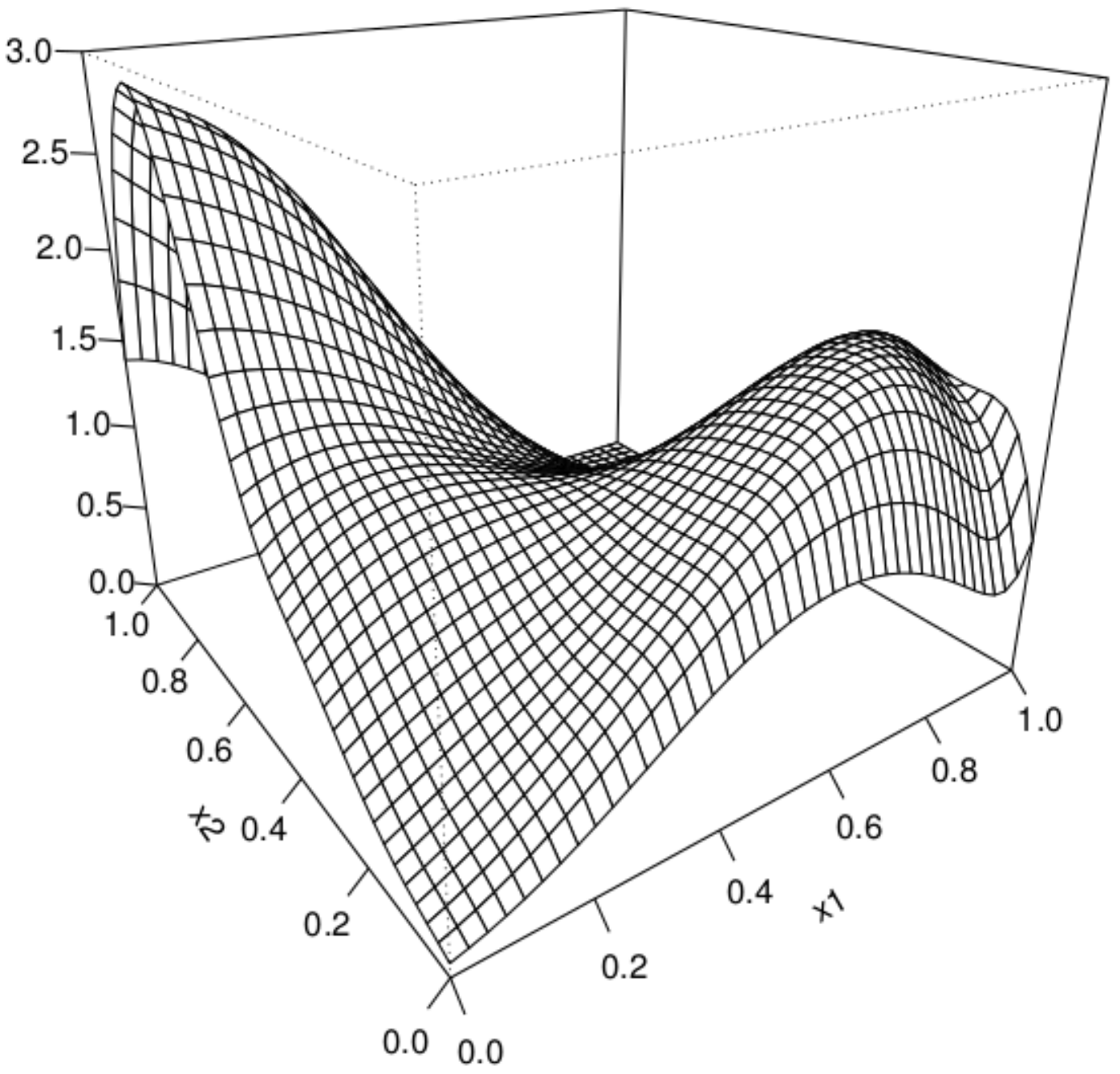} }  }
 \end{center}
	\caption{Graphical representations (o) of the real dataset with $n =80$, $\overline{x}_{1} = 0.4915$, $\overline{x}_{2} = 0.4126 $, $\widehat{\sigma}_{1} = 0.2757$, $\widehat{\sigma}_{2} = 0.2720$, $\widehat{\rho} = -0.6949$ and its corresponding smoothings according to the choice of bandwidth matrix $\widehat{\mathbf{H}}$ by LSCV: (a) full, (b) Scott and (c) diagonal.}\label{EstFull-Scottfig}
\end{figure}

We applied the standard version of the beta-Sarmanov estimator (\ref{est}) on paired rates data set according to the three types of bandwidth matrices in Table~\ref{parameterstable}. The dataset of sample size $n=80$ in Graph~(o) of Figure~\ref{EstFull-Scottfig} has been provided by Francial G. ~\cite{L13} during his last stay in Burkina Faso. It represents the popular ratings, for the two first consecutive ballots of the same electoral mandate of five years, of a political figure in different departments. Note that, for both elections, the prominent politician was finally elected in the second round. We are here interested to the opposite behavior of peoples during first rounds of both elections since the second round is governed by political alliances, which do not constitute a reference for the own popularity of a candidate. 

Indeed, for many African countries, political elections are generally fought on tribal ethnic origins and partisan interests. The data displays opposed viewpoints between the results of the first round of the first election ($x_{1}$) and the first round of the second one ($x_{2}$). In fact, the first election saw the candidate program mostly adopted by its clan (tribe and allied tribes) and rejected by the others. A few years later, facing the social discontent, he assumed a new political program which ends to another consultation of the people in a time $x_{2}$. This reversal made him loose the support of his supporters but he received the backing up of formers opponents. It is thus noticeable that its own popularity is not much different between the first round of both elections; the first gave an average $\overline{x}_{1} = 0.4915$ and the second $\overline{x}_{2} = 0.4126$. However, there is a significant negative correlation in the dataset: $\widehat{\rho} = -0.6949$. Unlike overall trends $\overline{x}_{1}$, $\overline{x}_{2}$ and $\widehat{\rho}$, the empirical distribution Graph~(o) of Figure~\ref{EstFull-Scottfig} gives more details of the electoral situation department by department. Hence, we need a nonparametric smoothing of this joint distribution by using associated kernels.

\begin{figure}[htbp!] 
 \begin{center}
  \mbox{
\subfloat[(A1) ]{\includegraphics[width=150pt,height=150pt,scale=0.97]{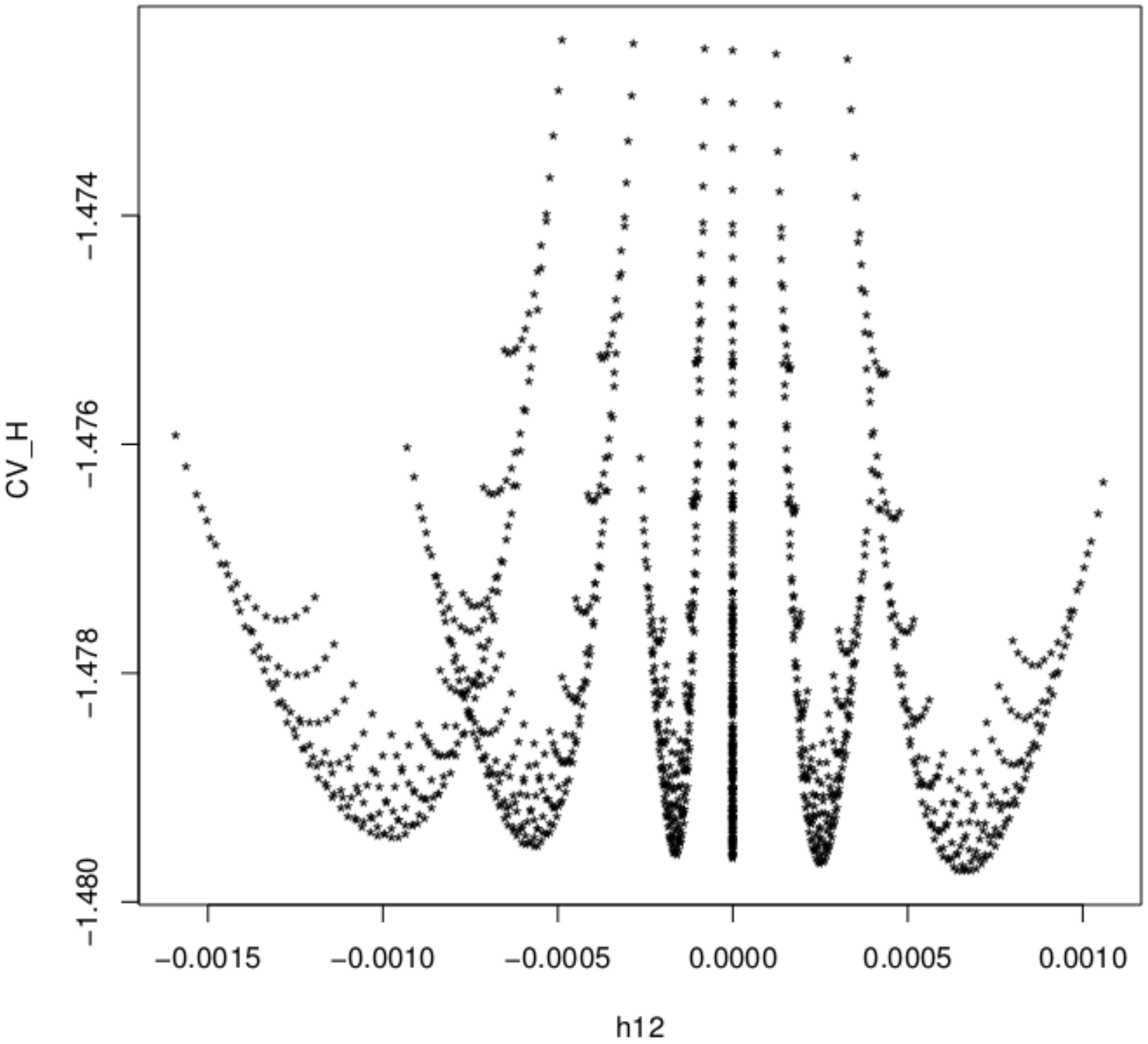} } \:
\subfloat[(A2) ]{\includegraphics[width=150pt,height=155pt,scale=0.97]{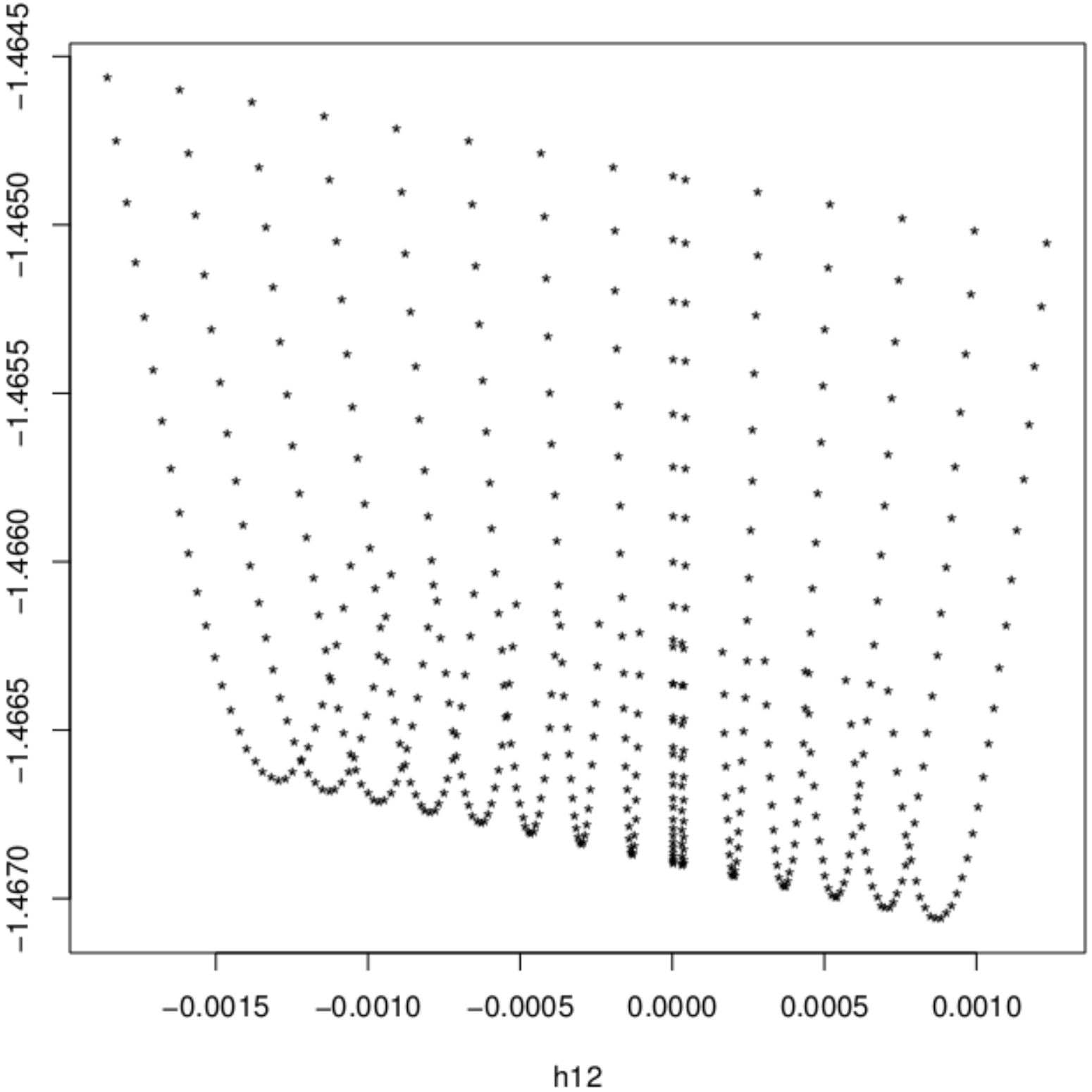} } \:
  \subfloat[(A3)]{\includegraphics[width=150pt,height=150pt,scale=0.97]{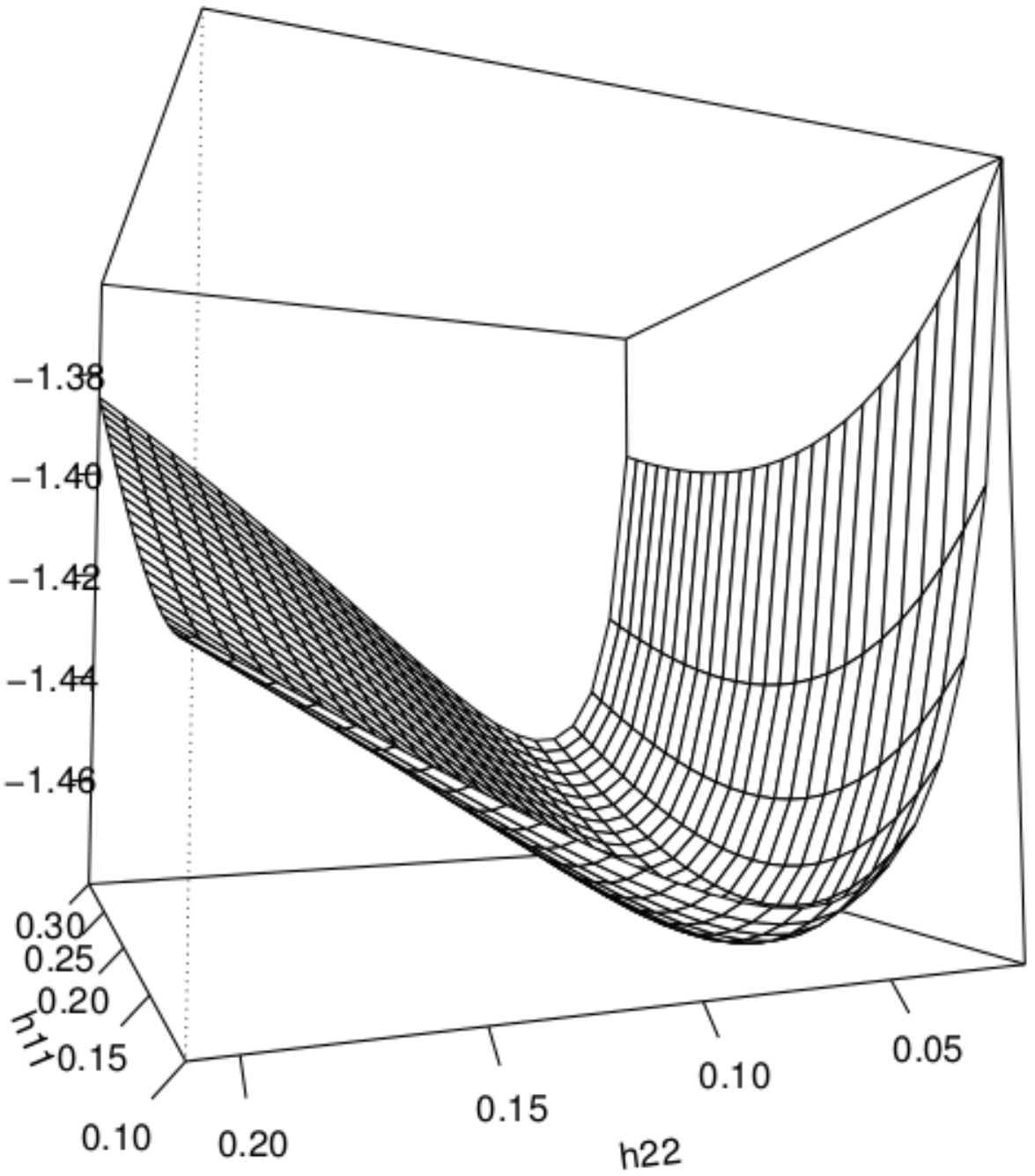} } } 
 \end{center}
\caption{Some plots of $\mathbf{H} \mapsto LSCV(\mathbf{H})$ from algorithms at the end of Section~\ref{sssec:standbeta} for real data in Graph (o) of Figure~\ref{EstFull-Scottfig}: (A1) full $\widehat{\mathbf{H}}$, (A2) Scott $\widehat{\mathbf{H}}$ and (A3) diagonal $\widehat{\mathbf{H}}$.}
\label{Hcvfig}
\end{figure}

In order to smooth the joint empirical distribution of these paired data, we apply the beta-Sarmanov kernel estimator in its standard version (\ref{est}). Figure~\ref{Hcvfig} shows the results of the LSCV algorithms A1, A2 and A3 with the ratings dataset; see (\ref{Hcv}) and the end of Section~\ref{sssec:standbeta}. The computation time of the LSCV is in the same trend as in Table~\ref{Timehcv} for $n=100$. To simplify the presentations in Figure~\ref{Hcvfig} for both the Scott and full bandwidth matrices, we only plot $h_{12} \mapsto LSCV(\mathbf{H})$ for some values of $h_{11}$ and $h_{22}$. In all cases we observe that there is a global minimum. 
The obtained optimal bandwidth matrices are
$$\widehat{\mathbf{H}}_{\rm{full}} = \begin{pmatrix}0.160000  & 0.000655 \\ 0.000655 & 0.057500\end{pmatrix},$$ $$\widehat{\mathbf{H}}_{\rm{Scott}} = 1.405\times \begin{pmatrix}0.076039  & 0.000622 \\ 0.000622 & 0.073994\end{pmatrix} = \begin{pmatrix}0.106836  & 0.000874 \\ 0.000874 & 0.103962\end{pmatrix} $$ 
and
$\widehat{\mathbf{H}}_{\rm{diagonal}} =  \mathbf{Diag}\left(0.160000,   0.057500\right)$. The resulting estimates are displayed in Figure~\ref{EstFull-Scottfig}. 

From simulation studies of previous section, the full bandwidth matrix provides the reference smoothing which is appropriated for correlated data; see Graph (a) of Figure~\ref{EstFull-Scottfig}.  In record time, the Scott bandwidth matrix $\widehat{\mathbf{H}}_{\rm{Scott}}$ delivers  similar smoothing in Graph (b) of Figure~\ref{EstFull-Scottfig} as the full and diagonal ones (see respectively Graphs (a) and (c) of Figure~\ref{EstFull-Scottfig}). In conclusion, we found anywhere the shape of a ``carpet flying'' in balance, smoothing thus the joint empirical distribution of the electoral situation of the candidate. This balance situation makes him to win in the second round of both elections.

\section{Summary and final remarks}
\label{sec:Conclusion}

We have presented general associated kernels (with or without correlation structure) that varying their shape according to the target point along the support. Excluding the classical associated kernels, the local adaptability of these associated kernels (depending on the target point $\mathbf{x}$ and the bandwidth matrix $\mathbf{H}$) means that they are free of boundary bias but have a slightly bias  different. Furthermore, the forms of bandwidth matrices used in the case with correlation structure have both theoretical and practical significances. Under the criterion of cross-validation, we therefore recommend the Scott bandwidth matrix which is more workable than the full one. A method of construction, called {\it multivariate mode dispersion} method, for these kernels are introduced. Also, we have proposed an algorithm of bias-reductions of their corresponding associated kernel estimators. An extension to discrete multivariate associated kernels is obviously possible. Similarly, a work is in progress on nonparametric multiple regression composed by continuous and discrete univariate associated kernels (e.g. Proposition~\ref{Proprod}).

Constructed by the correlation structure of \cite{S66} and by a variant of mode dispersion method, the bivariate beta-Sarmanov kernel estimator is completely study with the optimal bandwidth matrix chosen by cross-validation.  This technique can be extended in multivariate case for different type of kernels which are continuous and also discrete. In fact, from two or more univariate independent pdfs or probability mass functions, the correlation structure of \cite{S66} and a variant of mode dispersion method can always allow to build a multivariate type of kernel with correlation structure; and, therefore, produced the corresponding multivariate Sarmanov kernel. In terms of execution times from using the cross-validation method, we advise to use the Scott bandwidth matrix because of its flexibility and efficiency for no very strong correlation in the dataset. 

Simulation experiments and analysis of a real dataset  provide insight into the behavior of the bandwidth matrix for small and moderate sample sizes. Tables~\ref{Timehcv} and ~\ref{Err} and Figure~\ref{EstFull-Scottfig} can be conceptually summarized as follows. As expected, the full bandwidth matrix is frequently better than the others. An alternative with correlation structure has been proposed for the cross-validation technique: the Scott bandwidth matrix which has a comparable gain in execution times than the diagonal one; also, it produces better results than the diagonal in most cases. So, we recommend the Scott bandwidth matrix for multivariate use with the cross-validation technique. Further research in this direction are in progress, especially on the choice of optimal bandwidth matrix by using Bayesian approaches; see, e.g., \citet{ZAK14}. 


\section{Appendix}
\label{sec:Appendix}

This section is dedicated to the proof of Proposition~\ref{Pro3}. Indeed, using successively (\ref{EstversKern}) and Taylor's formula around $\mathbb{E}\left(\mathcal{Z}_{\theta(\mathbf{x}, \mathbf{H})}\right)$ and then $\mathbf{x}$, and also the invariance under cyclic permutations of the operator $\trace$, the result (\ref{biaisg}) is shown by
\begin{small}
\begin{eqnarray*}
 \Bias\left\{\widehat{f}_{n}(\mathbf{x})\right\} &=& \mathbb{E}\left\{ f\left(\mathcal{Z}_{\theta(\mathbf{x}, \mathbf{H})}\right)\right\} - f(\mathbf{x}) \nonumber \\
&=& f\left(\mathbb{E}\left(\mathcal{Z}_{\theta(\mathbf{x}, \mathbf{H})}\right)\right) +  \dfrac{1}{2}\mathbb{E}\left[\trace\left\{\left(\mathcal{Z}_{\theta(\mathbf{x}, \mathbf{H})} - \mathbb{E}\left(\mathcal{Z}_{\theta(\mathbf{x}, \mathbf{H})}\right)\right)^{\top}\nabla^{2}f\left(\mathbb{E} \left(\mathcal{Z}_{\theta(\mathbf{x}, \mathbf{H})}\right)\right)\right. \right. \\
&& \left. \left. \left(\mathcal{Z}_{\theta(\mathbf{x}, \mathbf{H})} - \mathbb{E}\left(\mathcal{Z}_{\theta(\mathbf{x}, \mathbf{H})}\right)\right)\right\}\right] - f(\mathbf{x}) + \oldstylenums{0}\left[\mathbb{E}\left\{\left(\mathcal{Z}_{\theta(\mathbf{x}, \mathbf{H})} - \mathbb{E}\left(\mathcal{Z}_{\theta(\mathbf{x}, \mathbf{H})}\right)\right)^{\top}\left(\mathcal{Z}_{\theta(\mathbf{x}, \mathbf{H})} - \mathbb{E}\left(\mathcal{Z}_{\theta(\mathbf{x}, \mathbf{H})}\right)\right)\right\}\right] \\
&=&f\left(\mathbb{E}\left(\mathcal{Z}_{\theta(\mathbf{x}, \mathbf{H})}\right)\right) - f(\mathbf{x}) +  \dfrac{1}{2}\trace\left[\nabla^{2}f\left(\mathbb{E} \left(\mathcal{Z}_{\theta(\mathbf{x}, \mathbf{H})}\right)\right)\mathbb{E}\left\{\left(\mathcal{Z}_{\theta(\mathbf{x}, \mathbf{H})} - \mathbb{E}\left(\mathcal{Z}_{\theta(\mathbf{x}, \mathbf{H})}\right)\right)^{\top}\right. \right. \\
&& \left. \left. \left(\mathcal{Z}_{\theta(\mathbf{x}, \mathbf{H})} - \mathbb{E}\left(\mathcal{Z}_{\theta(\mathbf{x}, \mathbf{H})}\right)\right)\right\}\right]+ \oldstylenums{0}\left[\trace\left(\mathbb{E}\left\{\left(\mathcal{Z}_{\theta(\mathbf{x}, \mathbf{H})} - \mathbb{E}\left(\mathcal{Z}_{\theta(\mathbf{x}, \mathbf{H})}\right)\right)\left(\mathcal{Z}_{\theta(\mathbf{x}, \mathbf{H})} - \mathbb{E}\left(\mathcal{Z}_{\theta(x, \mathbf{H})}\right)\right)^{\top}\right\}\right)\right] \\
&=&f\left(\mathbb{E}\left(\mathcal{Z}_{\theta(\mathbf{x}, \mathbf{H})}\right)\right) - f(\mathbf{x}) +  \dfrac{1}{2}\trace\left[\nabla^{2}f\left(\mathbb{E} \left(\mathcal{Z}_{\theta(\mathbf{x}, \mathbf{H})}\right)\right)\mbox{Cov}\left(\mathcal{Z}_{\theta(\mathbf{x}, \mathbf{H})}\right)\right] + \oldstylenums{0}\left[ \trace\left\{\mbox{Cov}\left(\mathcal{Z}_{\theta(\mathbf{x}, \mathbf{H})}\right)\right\}\right]\\
&=&f\left(\mathbf{x} +  \mathbf{a}_{\theta}(\mathbf{x}, \mathbf{H})\right) - f(\mathbf{x}) +  \dfrac{1}{2}\trace\left[\nabla^{2}f\left(\mathbf{x} + \mathbf{a}_{\theta}(\mathbf{x}, \mathbf{H})\right)\mathbf{B}_{\theta}(\mathbf{x}, \mathbf{H})\right] + \oldstylenums{0}\left[ \trace\left\{\mathbf{B}_{\theta}(\mathbf{x}, \mathbf{H})\right\}\right]\\
&=& \mathbf{a}^{\top}_{\theta}(\mathbf{x}, \mathbf{H})\nabla f\left(\mathbf{x}\right) + \dfrac{1}{2}\trace\left[\left\{\mathbf{a}_{\theta}(\mathbf{x}, \mathbf{H})\mathbf{a}^{\top}_{\theta}(\mathbf{x}, \mathbf{H}) + \mathbf{B}_{\theta}(\mathbf{x}, \mathbf{H}) \right\}\nabla^{2}f\left(\mathbf{x}\right)\right] + \oldstylenums{0}\left\{\trace\left(\mathbf{H}^{2}\right)\right\}.
\end{eqnarray*}
\end{small}
In fact, the rest $\oldstylenums{0}\left\{\trace\left(\mathbf{H}^2\right)\right\}$ comes from  $\trace(\mathbf{B}_{\theta}(\mathbf{x}, \mathbf{H})) =  O\left(\trace\mathbf{H}^{2}\right)$ deduced from Proposition~\ref{Pro1} of classical associated kernels and 
\begin{eqnarray*}
\oldstylenums{0}\left[\mathbb{E}\left\{\left(\mathcal{Z}_{\theta(\mathbf{x}, \mathbf{H})} - \mathbb{E}\left(\mathcal{Z}_{\theta(\mathbf{x}, \mathbf{H})}\right)\right)^{\top}\left(\mathcal{Z}_{\theta(\mathbf{x}, \mathbf{H})} - \mathbb{E}\left(\mathcal{Z}_{\theta(\mathbf{x}, \mathbf{H})}\right)\right)\right\}\right] &=& \\  \mathbb{E}\left\{\oldstylenums{0}_{p}\left(\mathcal{Z}_{\theta(\mathbf{x}, \mathbf{H})} - \mathbb{E}\left(\mathcal{Z}_{\theta(\mathbf{x}, \mathbf{H})}\right)\right)^{\top}\left(\mathcal{Z}_{\theta(\mathbf{x}, \mathbf{H})} - \mathbb{E}\left(\mathcal{Z}_{\theta(\mathbf{x}, \mathbf{H})}\right)\right)\right\}, 
\end{eqnarray*} where $\oldstylenums{0}_{p}(\cdot)$ is the probability rate of convergence. 

Concerning the variance (\ref{var}) one first has
\begin{eqnarray*}
   \Var\left\{\widehat{f}_{n}(\mathbf{x})\right\} &=& \frac{1}{n} \mathbb{E}\left\{K^{2}_{\theta(\mathbf{x}, \mathbf{H})}\left(\mathbf{X}_{1}\right)\right\} -  \frac{1}{n} \left[\mathbb{E}\left\{K_{\theta(\mathbf{x}, \mathbf{H})}\left(\mathbf{X}_{1}\right)\right\}\right]^{2} \\
&=& \dfrac{1}{n} \int_{\mathbb{S}_{\theta(\mathbf{x}, \mathbf{H})} \cap \mathbb{T}_{d}} K^{2}_{\theta(\mathbf{x}, \mathbf{H})} (\mathbf{u})f(\mathbf{u})d\bold{u} - \frac{1}{n} \left[\mathbb{E}\left\{K_{\theta(\mathbf{x}, \mathbf{H})}\left(\mathbf{X}_{1}\right)\right\}\right]^{2} \\
&=& I_{1} - I_{2}.
\end{eqnarray*}
From (\ref{EstversKern}) and (\ref{biaisg}), one has the following behavior of the second term 
$$I_{2} := (1/n)\left[\mathbb{E}\left\{K_{\theta(\mathbf{x}, \mathbf{H})}\left(\mathbf{X}_{1}\right)\right\}\right]^{2} \simeq (1/n)f^{2}\left(\mathbf{x}\right) \simeq O\left(1/n\right)$$
since $f$ is bounded for all $\mathbf{x} \in \mathbb{T}_{d}$. By using Taylor's expansion around of $\mathbf{x}$, the first term $I_{1} :=   (1/n) \int_{\mathbb{S}_{\theta(\mathbf{x}, \mathbf{H})} \cap \mathbb{T}_{d}} K^{2}_{\theta(\mathbf{x}, \mathbf{H})} (\mathbf{u})f(\mathbf{u})d\bold{u}$ becomes
\begin{equation*}
  I_{1}  =  \dfrac{1}{n} f(\mathbf{x})\int_{\mathbb{S}_{\theta(\mathbf{x}, \mathbf{H})} \cap \mathbb{T}_{d}} K^{2}_{\theta(\mathbf{x}, \mathbf{H})} (\mathbf{u})d\bold{u} + R\left(\mathbf{x}, \mathbf{H}\right),
\end{equation*}
 with
\begin{eqnarray*}
 R\left(\mathbf{x}, \mathbf{H}\right) &=& \dfrac{1}{n} \int_{\mathbb{S}_{\theta(\mathbf{x}, \mathbf{H})} \cap \mathbb{T}_{d}} K^{2}_{\theta(\mathbf{x}, \mathbf{H})} (\mathbf{u})\left[\left(\mathbf{u} - \mathbf{x}\right)^{\top}\nabla f(\mathbf{x}) + \dfrac{1}{2} \left(\mathbf{u}- \mathbf{x}\right)^{\top} \nabla^{2} f (\mathbf{x}) \left(\mathbf{u}- \mathbf{x}\right) \right. \\
  &&\left. +~ \oldstylenums{0}\left\{\left(\mathbf{u} - \mathbf{x}\right)^{\top}\left(\mathbf{u}- \mathbf{x}\right)\right\} \right]d\bold{u}.
\end{eqnarray*}
A similar argument from \citet[Lemma]{C99} with $f$ bounded on $\mathbb{T}_{d}$ shows the existence of $r_{2}$ and then the condition $||K_{\theta(\mathbf{x}, \mathbf{H})}||^{2}_{2}  \lesssim  c_{2}(\mathbf{x})(\det\mathbf{H})^{-r_{2}}$ leads  successively to

\begin{eqnarray*}
0 \;\:\leq \;\:  R\left(\mathbf{x}, \mathbf{H}\right)  &\lesssim& \dfrac{1}{n\;(\det\mathbf{H})^{r_{2}}} \int_{\mathbb{S}_{\theta(\mathbf{x}, \mathbf{H})} \cap \mathbb{T}_{d}} c_{2}(\mathbf{x})\left\{\left(\mathbf{u} - \mathbf{x}\right)^{\top}\nabla f(\mathbf{x}) + \dfrac{1}{2} \left(\mathbf{u} - \mathbf{x}\right)^{\top} \nabla^{2} f (\mathbf{x}) \left(\mathbf{u} - \mathbf{x}\right)\right\}d\bold{u}\\ &&\simeq \oldstylenums{0} \left\{n^{-1}(\det\mathbf{H})^{-r_{2}}\right\}.\blacksquare
\end{eqnarray*}

\section*{Acknowledgements}
We sincerely thank Francial G. Libengu\'e for useful discussions and for the dataset of illustration.

\section*{References}

\end{document}